\documentclass[11pt]{article}   	                 
\usepackage{amsmath,amsthm, amssymb,amsfonts,amscd,mathabx}

\usepackage[colorlinks=true,linkcolor=blue,citecolor=red]{hyperref}
\usepackage[pdftex]{graphicx}                 
\usepackage{wrapfig}                               
\usepackage[font=footnotesize]{caption}                     
\usepackage{float}
\usepackage[symbol]{footmisc}
\counterwithin{figure}{section}

\let\OLDitemize\itemize
\renewcommand\itemize{\OLDitemize\addtolength{\itemsep}{-5pt}}

\numberwithin{equation}{section}
\newtheorem{theorem}{Theorem}[section]
\newtheorem{conjecture}[theorem]{Conjecture}
\newtheorem{lemma}[theorem]{Lemma}

\newtheorem{definition}[theorem]{Definition}
\newtheorem{corollary}[theorem]{Corollary}

\newtheorem{remark}[theorem]{Remark}
{\begin{trivlist}\item[]\textbf{Proof#1 }}%
{\qed\end{trivlist}}

 {\begin{trivlist}\item[]\textbf{Acknowledgments }}{\end{trivlist}}

\newcommand{\R}{\mathbb{R}}
\newcommand{\C}{\mathbb{C}}
\newcommand{\N}{\mathbb{N}}
\newcommand{\Z}{\mathbb{Z}}

\newcommand{\T}{\mathbb{T}}

\newcommand{\Rmnum}[1]{\uppercase\expandafter{\romannumeral #1\relax}}

\newcommand{\rmO}{\mathrm{O}}

\newcommand{\rmd}{\mathrm{d}}
\newcommand{\rme}{\mathrm{e}}
\newcommand{\rmi}{\mathrm{i}}

\renewcommand{\Re}{\mathrm{Re}\,}
\renewcommand{\Im}{\mathrm{Im}\,}

\renewcommand{\leq}{\leqslant}
\renewcommand{\geq}{\geqslant}

\def\im{\mathop{\mathrm{\,Im}\,}}

\def\eps{\varepsilon}

\newfam\bifam
\font\tenbi=cmmib10 scaled \magstep1 \font\sevenbi=cmmib10 at 11pt
\font\fivebi=cmmib10 at 6pt \textfont\bifam = \tenbi
\scriptfont\bifam = \sevenbi \scriptscriptfont\bifam= \fivebi

\begin{document}

\begin{center}
{\fontsize{14}{14}\fontfamily{cmr}\fontseries{b}\selectfont{
Selection mechanisms in front invasion\footnote[4]{
The authors gratefully acknowledge partial support through NSF grants DMS-2510541 \& DMS-2202714 (M.A.), DMS-2406623 (M.H.),  and DMS-2205663 \& DMS-2506837 (AS).
}}}\\[0.2in]
Montie Avery$^{1}$, Matt Holzer$^2$, and  Arnd Scheel$^{3}$\\[0.1in]
\textit{\footnotesize
$^1$Department of Mathematics, Emory University, 400 Dowman Drive, Atlanta, GA 30322, USA \\[0.02in]
$^2$Department of Mathematical Sciences, George Mason University, Fairfax, VA 22030, USA\\[0.01in]
$^3$School of Mathematics, University of Minnesota, Minneapolis, 206 Church St SE, MN 55414, USA
}
\end{center}
\vspace*{.4in}
\begin{abstract}
\noindent
We review progress on questions related to front propagation into unstable states and point out open problems in the area. We  strive to highlight different theoretical perspectives and challenges while also addressing more practical questions with examples and guides to computational methods. Throughout we take a dynamical systems point of view that focuses on the ability of invasion processes to act as a selection mechanism in complex systems.
\end{abstract}


%
%

 \newpage
\tableofcontents
\section{Introduction}

When studying dynamical systems, one is concerned with the long-time behavior of solutions to at times complex models arising in the sciences, in engineering, or in mathematics. In the simplest cases, one seeks to establish convergence to a stable asymptotic state. More interesting behavior arises when a trivial asymptotic state is unstable. One then seeks to describe the fate of perturbations of this unstable state, following for instance its unstable manifold. In large systems, this unstable manifold is usually very high-dimensional and there is in general no simple description of the fate of all solutions in this unstable manifold.
Spatially extended systems are an example of such complex systems. In a translation-invariant system, instability of a trivial state can manifest itself from a small localized perturbation at any arbitrary location, thus contributing to a high-dimensional unstable manifold.

In many contexts, one is specifically interested in  such \emph{spatially localized} perturbations and would like to track their growth temporally and spatially. Monitoring the spatial leading edge of the perturbations, one then tries to characterize the evolution of the disturbance through the evolution of a propagating front solution. Despite the complexity arising through the high dimensionality of the unstable manifold, one then often observes spatio-temporal behavior that is largely independent of the initial perturbation: the growth of disturbances, mediated by fronts, propagates  at selected speeds and selects distinguished states in its wake.

In experimental or naturally observed settings, unstable states typically arise from a sudden parameter change or quench, that renders a previously stable, trivial equilibrium unstable. Alternatively, a mode of instability was simply not present in the system and is newly introduced, as is the case in many problems in ecology and biology, where invasive species spread in an environment that was in a stable ecological equilibrium until the introduction of the new species. Similarly in spirit, introduction of a new reactant may trigger instabilities in combustion, or mutated cells may displace healthy tissue in cancer growth through cell proliferation and migration.

A slightly different motivation arises when an instability is subject to an advection mechanism. Close to the instability threshold, perturbations grow in amplitude but are advected out of a fixed observational window. {\color{black} Far from the onset of instability, perturbations grow in any fixed window.} This {\color{black} former, \emph{convective}, instability changes into the latter, \emph{absolute}, instability} only at a secondary threshold, when the growth of perturbations in the fixed window of observation is marginal, neither exponentially growing nor decaying. Again, the interest is in determining the speed of invasion, particularly the critical parameter value at which the speed is zero in the laboratory frame.

Mathematical approaches toward a description of front invasion phenomena originate in early work by Fisher~\cite{fisher} and Kolmogorov, Petrov, and Piscounov~\cite{kpp}, who studied scalar parabolic equations; see also \cite{aronson,uchiyama}. Independently, instabilities in plasmas and fluids led to the development of methods to track instabilities pointwise in  linear equations~\cite{sturrock,bers1983handbook,briggs} with far-ranging implications for growth and invasion problems in the physical sciences. Dee and Langer~\cite{deelanger} predicted patterns in growth processes based on these linear criteria, which were more recently matched with a nonlinear front to obtain corrections for speeds and wavenumbers by Ebert and van Saarloos~\cite{EbertvanSaarloos}. Those results expand on logarithmic corrections found earlier in Bramson's probabilistic approach~\cite{Bramson1,Bramson2}, which in turn was yet more recently reproduced and refined using PDE methods in order preserving systems~\cite{Comparison1, Comparison2, NRRkppasy, Graham}.

The work we describe here starts with the observation that most mathematical results on front invasion hinge on the presence of comparison principles, while much of the interest in the sciences has centered on systems without any ordering structure; we refer to~\cite{vanSaarloosReview} for an extensive review of this latter point of view. We therefore emphasize here mathematical tools in front invasion that do not rely on comparison principles, taking  a dynamical systems view point. We hope that such a more general and yet rigorous  mathematical approach can add clarity to our understanding of front invasion, and also detect novel phenomena, making assumptions explicit and pointing to situations where those fail.

The purpose of this review is to demonstrate how recent work has made progress on both aspects, bringing rigor to the selection of fronts, and pointing towards novel phenomena. Importantly, we feel that there are indeed a multitude of interesting mathematical questions that arise in this context, and we hope that this review will stimulate interested researchers to contribute to research on front selection. On the other hand, we hope that the mathematical clarity that comes with a more rigorous approach to front selection will also help more applied researchers use predictions and concepts in their work and contribute to the development of effective numerical algorithms. This review therefore includes in each topic, next to a  description of theory and open problems, a practical guide to applying the results, numerically or analytically.

The remainder of this introduction contains a selective illustration of both the theory in order-preserving systems and challenges and phenomena that arise in pattern-forming systems.

%
%

\subsection{The simplest example: existence, positivity, and stability}\label{s : FKPPintro}

\begin{figure}
\includegraphics[height=1.6in]{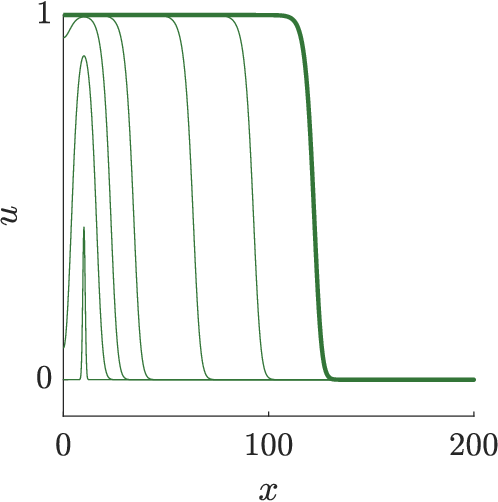}\hfill
\includegraphics[height=1.6in]{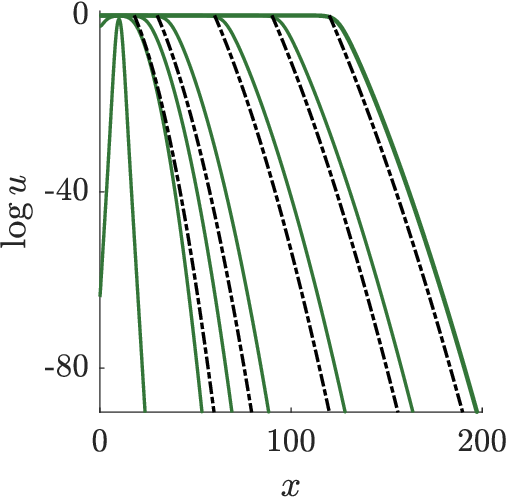}\hfill
\includegraphics[height=1.6in]{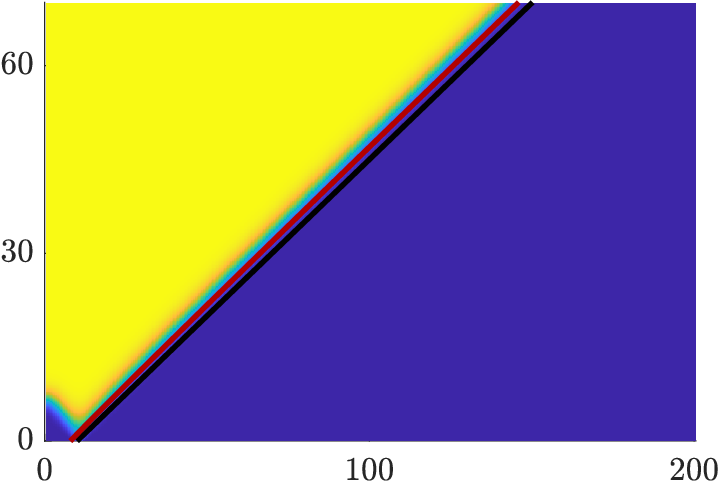}
\caption{Solution to~\eqref{e:fkpp} with small Gaussian initial condition near $x=10$, snap shots of profiles (left) and logarithms of profiles (center) at times $t=0,5,9,15,30,45,60$. Also shown in the center are parabolas corresponding to Gaussians solving the linear equation (black dashed). The space-time plot (right) shows spreading with speed 2 (black line) and even better agreement when including logarithmic corrections {\color{black}$x=2t-\frac{3}{2}\log t$} (red). }
\label{f:kppspacetime}
\end{figure}

Kolmogorov, Petrovsky, and Piscounov~\cite{kpp} and, independently, Fisher~\cite{fisher} studied the equation
\begin{equation}
u_t=u_{xx}+f(u), \qquad x\in\R,\ t>0,\qquad  u=u(t,x)\in\R,
\label{e:fkpp}
\end{equation}
\textcolor{black}{in the particular case where $ f(u)=u(1-u)$.} Here subscripts denote partial derivatives. In fact, it was proven in~\cite{kpp} that \emph{step-like initial conditions} of the specific form $u(0,x)=0$ for $x>0$ and $u(0,x)=1$ for $x<0$ propagate with \emph{spreading speed} $c_*=2$, stated for instance in the form that for all $x$,
\begin{equation}\label{e:spsp}
\lim_{t\to\infty}u(t,x-ct)=\left\{\begin{array}{ll}
0,& \text{ if }c>c_*,\\
1,& \text{ if } c<c_*;
\end{array}\right.
\end{equation}
see Fig.~\ref{f:kppspacetime}. 
This spreading speed depends on the specific form of $f$. For $f(u)=u(1-u)(u+a)$ with $a>0$, the Nagumo equation, \eqref{e:spsp} still holds for step-like initial conditions, but now with
\begin{equation}
    c_*=({1+2a})/{\sqrt{2}},\text{ when }a<1/2,\qquad 
    c_*=2\sqrt{f'(0)}= 2\sqrt{a},\text{ when }a\geq 1/2.
\end{equation}
We next explain a variety of aspects of this result. \textcolor{black}{We phrase our discussion in terms of a general nonlinearity $f(u)$ with $f(0)=0$ and $f'(0)>0$ keeping in mind the specific examples $f(u)=u(1-u)$ (FKPP) and $f(u)=u(1-u)(u+a)$ with $a>0$  (Nagumo). }

\begin{figure}
\centering\includegraphics[scale=1]{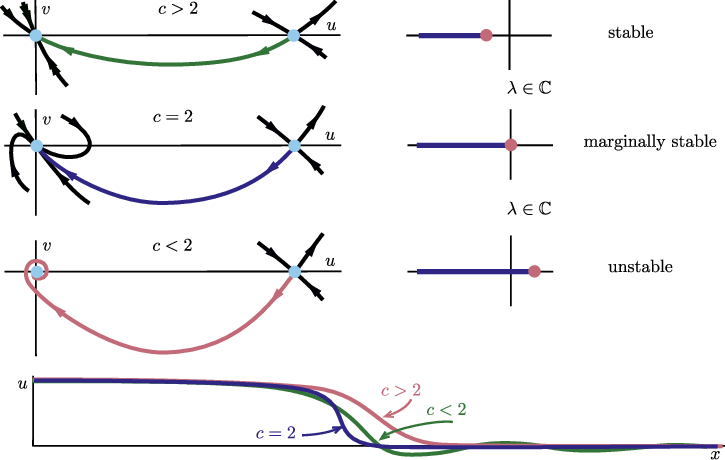}
\caption{Phase portrait sketches in the $u-u_x$-plane of~\eqref{e:kpptw}, $f(u)=u(1-u)$ (left), with $c$ decreasing from top to bottom. Also shown are spectra of the linearized operator in $L^p_{\mathrm{exp},c/2}$ (right); see~\eqref{e:Lpweighted} for function spaces. The marginally stable front is in fact the steepest front as illustrated in the bottom panel with plot sketches of $u_*(x;c)$.}
\label{f:kppphaseportrait}
\end{figure}

\begin{figure}
\centering\includegraphics[scale=1]{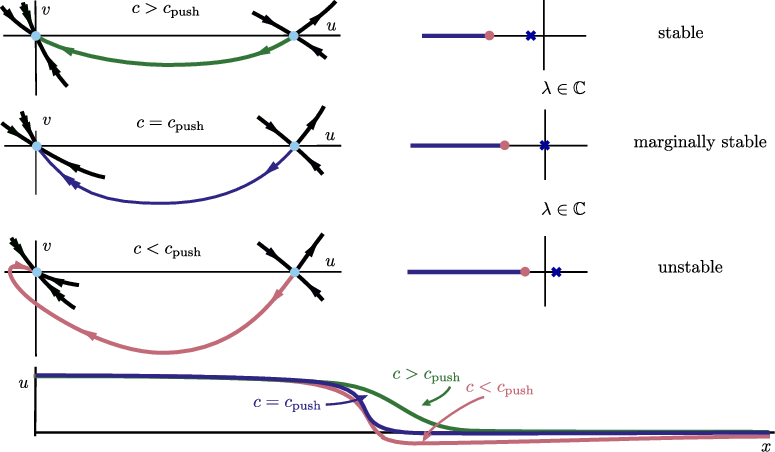}
\caption{Phase portrait sketches in the $u-u_x$-plane of~\eqref{e:kpptw}, $f(u)=u(1-u)(u+a)$ for $a<\frac{1}{2}$ (left), with $c$ decreasing from top to bottom. Also shown are spectra of the linearized operator in $L^p_{\mathrm{exp},c/2}$ (right); see~\eqref{e:Lpweighted} for function spaces. The marginally stable front is in fact the steepest front as illustrated in the bottom panel with plot sketches of $u_*(x;c)$.} \label{f:nagumophaseportrait}
\end{figure}

\textbf{Existence of traveling waves.} Propagating front solutions can be sought in the form $u(x-ct)$, which yields an ordinary differential equation in the comoving spatial variable $\xi=x-ct$,
\begin{equation}\label{e:kpptw}
u''+cu'+f(u)=0, \qquad \text{or }\quad  \left\{\begin{array}{ll}u'&=v\\v'&=-cv-f(u)\end{array}\right. .
\end{equation}
The phase plane analysis of this system \textcolor{black}{ in the case $f(u)=u(1-u)$} shows that there are precisely two equilibria: $p_-=(1,0)$ is a saddle and $p_+=(0,0)$ is a sink when $c>0$. A variety of shooting arguments can then show that there exists, for each $c>0$, a unique heteroclinic solution $(u_*,v_*)(\xi;c)$ with $(u_*,v_*)(\xi;c)\to p_\pm$ for $\xi\to \pm\infty$; see Fig.~\ref{f:kppphaseportrait}.
Note that the equilibrium $p_+$ is a node for $c>2$ and a focus for $c<2$. For $c>2$, the heteroclinic does not lie in the strong stable manifold of the equilibrium $p_+$. The former fact depends only on the linearization, that is, on the fact that $f'(0)=1$. The latter fact strongly depends on the shape of the nonlinearity, \textcolor{black}{exploiting for instance sublinear growth when $f(u)=u(1-u)$: for $u>0$, $f(u)< f'(0)u$. } In fact, when studying the Nagumo equation
\begin{equation}\label{e:nagumo}
u_t=u_{xx}+f(u),\qquad f(u)=u(1-u)(u+a),\quad a>0,
\end{equation}
the heteroclinic belongs to the strong stable manifold of $p_+$ when $a<1/2$ and $c=(1+2a)/\sqrt{2}$,  in which case the solution is explicitly given through $u_*(\xi)=(1+\exp(\xi/\sqrt{2}))^{-1}$; see Fig.~\ref{f:nagumophaseportrait}.

Following the nomenclature in \cite{stokes1976two}, we loosely refer to the case when the heteroclinic is never contained in the strong stable manifold as the \emph{pulled} case, with the front at speed $c=2$ as the pulled front, and to the other case as \emph{pushed} with the front contained in the strong stable manifold as the \emph{pushed front}; we will give {\color{black} alternative} descriptions below and {\color{black} precise} definitions in \S\ref{s:nlmsp}.

\textbf{Steepest fronts.} In order to obtain the exponential rate of decay of fronts as $x\to\infty$, we solve the stationary dispersion relation $\nu^2+c\nu+f'(0)=0$ and find the roots
\begin{equation}\label{e:nu}
\nu_\pm(c)=\frac{1}{2}\left(-c\pm \sqrt{c^2-4f'(0)} \right), \qquad \nu_-(c)<\nu_+(c).
\end{equation}
If the solution is not contained in the strong stable manifold as noted above in the case $f(u)=u(1-u)$, then the steepest decay occurs at
\begin{equation}\label{e:steep}
c_\mathrm{steep}=\mathrm{argmin}_c\{ \Re\nu_+(c)\}=2 \sqrt{f'(0)},
\end{equation}
which happens to be the speed at which step-function like initial conditions spread. If the traveling wave lies in the strong stable manifold for some $c_0>2$, we find
{\color{black} \begin{equation}\label{e:steeppush}
\nu_-(c_0)<2 \sqrt{f'(0)}\leq  \Re\nu_+(c),
\end{equation}
}
so that the front at this speed is steeper than all fronts not contained in the strong stable manifold. One can show \textcolor{black}{in the specific case of the cubic nonlinearity but also for general $f(u)$, }exploiting again the structure of the traveling-wave equation, that there exists at most one such value $c_0$ for which the front is contained in the strong stable manifold. It turns out that this front gives the spreading speed of step-function like initial conditions~\cite{HadelerRothe}.

In summary, we find that in this simple example of FKPP and Nagumo equations, but in fact also for general $f(u)$, we always find that the spreading speed of step-function like initial conditions is given by the speed of the \emph{steepest front}.

\textbf{Positivity.}
The solution $u_*(x;c)$ is strictly positive for all values $c\geq 2$. This can be seen directly for this specific nonlinearity designing invariant regions in phase space and exploiting sublinear growth of $f$ for $u>0$. Inspecting the phase portrait, positivity then follows from
\begin{itemize}
\item[(i)] the fact that the asymptotics of the solution are monotone for $c\geq 2$; and
\item[(ii)] the fact that the unstable manifold of $p_-$ does not intersect the strong stable manifold of $p_+$.
\end{itemize}
Clearly, (i) ceases to hold when the linearization at $p_+$ possesses complex eigenvalues. Failure of (ii)  would typically lead to a ``flip'' of the heteroclinic to the other side of the stable manifold. This phenomenon in fact occurs in~\eqref{e:kpptw}  when $f(u)=u(1-u)(u+a)$ with $a>1/2$, the pushed case discussed above.
\textcolor{black}{In all those scenarios,} the spreading speed of step-function like initial conditions is given by the \emph{slowest positive front}.

\textbf{Stability --- monotonicity.}
All fronts are equilibria in a comoving frame, that is, time-independent solutions to
\[
u_t=u_{\xi\xi}+cu_\xi + f(u).
\]
In order to understand if those solutions are indeed observable, one would first ask if small perturbations to a front profile stay bounded or even decay as $t\to\infty$. For monotone profiles, one may now simply consider initial conditions wedged between translates of the front, $u_*(\xi+\delta;c)<u(t=0;\xi)<u_*(\xi-\delta;c)$. The parabolic comparison principle then guarantees that this inequality holds for all $t>0$, so that the solution stays uniformly close to $u_*(\xi;c)$.

\textbf{Stability --- linearization.}
A different perspective is gained when thinking of the stability problem perturbatively. One therefore sets $u=u_*(\xi;c)+v(t,\xi)$ to find the equation
\[
v_t=\mathcal{L_*}v + N(v),
\]
with
\[
\mathcal{L}_*v=v_{\xi\xi}+cv_\xi + f'(u_*(\xi))v,\qquad\quad  N(v,\xi)=f(u_*(\xi)+v)-f(u_*(\xi))-f'(u_*(\xi))v=\rmO(|v|^2).
\]
One then wishes to guarantee that the spectrum of the elliptic operator $\mathcal{L}_*$ is negative, so that solutions to the linear equation decay exponentially, allowing for control of the nonlinear terms. For this, one necessarily needs to introduce weighted spaces: the essential spectrum of $\mathcal{L}_*$ contains the essential spectrum of the linearization at the leading edge, $\mathcal{L}_+v=v_{\xi\xi}+cv_\xi + f'(0)v$ which one readily finds to be unstable since $f'(0)>0$. In fact, Fourier transform shows that $\mathrm{spec}(\mathcal{L}_+)=\{-k^2+c\rmi k +f'(0) \, | \, k\in\R\}$ so that $\max\Re\,\mathrm{spec}(\mathcal{L}_+)=f'(0)$, in translation invariant norms such as $L^p(\R)$.

We define weighted $L^p$-spaces, $1\leq p\leq \infty$ for a weight $\omega(x)>0$ smooth, $\omega(x)=\rme^{\eta_+x},\ x>1$, $\omega(x)=\rme^{\eta_-x},\ x<-1$,
\begin{equation}
L^p_{\mathrm{exp},\eta_-,\eta_+}=\{u\in L^p_\mathrm{loc}(\R)|\|u\|_{L^p_{\mathrm{exp},\eta_-,\eta_+}}<\infty\},\qquad
\|u\|_{L^p_{\mathrm{exp},\eta_-,\eta_+}}=\|\omega u\|_{L^p},\qquad  L^p_{\mathrm{exp},\eta}=L^p_{\mathrm{exp},\eta,\eta}. \label{e:Lpweighted}
\end{equation}
Optimizing exponential weights, we wish to minimize $\max_{\eta_\pm}\{\Re\,\mathrm{spec}(\mathcal{L}_*)\}$, which in this simple example is achieved through the choice $\eta_-=\eta_+=c/2$. Fig.~\ref{f:nagumophaseportrait} shows how  the spectra of $\mathcal{L}_*$ change in this optimal weight as the speed varies. Coincidentally, we find
\begin{center}
$\Re\,\mathrm{spec}\,\mathcal{L}_*\leq 0$ \qquad $\Longleftrightarrow$ \qquad $u_*(\xi;c)$ is positive.
\end{center}
\textcolor{black}{We observe that the exponential decay rate of $\rme^{\mathcal{L}_*t}$, $|\max\Re\mathrm{spec}\,(\mathcal{L}_*)|,$ \emph{decreases} as $c$ decreases and vanishes at a critical wave speed $c_*$, with $c_*=2$ when $a>1/2$ and $c_*=(1+2a)/\sqrt{2}$ when $a<1/2$, which coincidentally is the spreading speed of a step-function initial condition.}


\textbf{Selection --- definition.}
\textcolor{black}{Slightly expanding on the characterization of spreading speeds as the asymptotic behavior of step-like initial conditions made explicit in \eqref{e:spsp}, we say a family of fronts $\{u(\xi+\xi_0;c)\mid \xi_0\in\R\}$, with speed $c$ is \emph{selected} if it has, in  suitable norms, an open basin of attraction that includes initial conditions that vanish in $x>0$. See Def.~\ref{def: selected front} for a precise statement.}

\textbf{Selection versus stability.} {\color{black}
It is typically the case that stable fronts, with an open, non-empty basin of attraction, exist for a continuum of speeds. However, many of these fronts will only attract initial data with well-prepared, slowly decaying tails. Such initial data are non-physical for invasion problems starting from localized disturbances, hence the distinction of selected fronts as those which attract steep, physically relevant data.
}

\textbf{Selection criteria: marginal positivity,  steepest decent, and marginal stability.}
In the scalar equations we have considered thus far, in particular in the FKPP and the Nagumo equation, we found that, somewhat coincidentally, selected fronts are the steepest fronts, they are the slowest positive fronts, and they are marginally stable  (in optimal exponential weights). In many situations, ordering structures have been shown to imply that the slowest positive front is in fact selected. On the other hand, a long-standing conjecture that the marginally stable front is in fact selected has recently been established absent any ordering principle~\cite{as1,avery2}. We shall revisit how these selection criteria relate to each other in the discussion,~\S\ref{s:dis}.

\textbf{Pulled and pushed fronts --- definition via marginal stability.}
\textcolor{black}{Based on the characterization of selection via linear marginal stability, we can then define a pulled front $u_*(x;c_*)$ as a front where marginal stability is due to essential spectrum of the linearization (in optimally weighted spaces) crossing the imaginary axis, while for a pushed front $u_*(x;c_*)$  point spectrum crosses the imaginary axis, as the speed $c$ decreases past $c_*$; see Defs.~\ref{def: rigid pulled} and~\ref{def: rigid pushed}, respectively, for details, and the right panels of Figs.~\ref{f:kppphaseportrait} and~\ref{f:nagumophaseportrait}, respectively, for illustrations of spectra.}

\begin{figure}
\centering\includegraphics[scale=1]{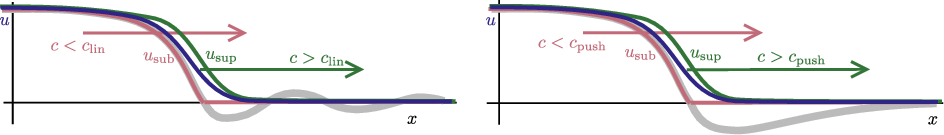}
\caption{Sub- and super-solutions that trap  solutions between speeds $c_*\pm \eps$ in the pulled (left, $c_*=c_\mathrm{lin}$) and in the pushed case (right, $c_*=c_\mathrm{push}$).}\label{f:subsuper}
\end{figure}

\textbf{Selection of the slowest positive front --- sub- and super-solutions.}
Exploiting the comparison principle, one can see very quickly and intuitively why the slowest positive front is in fact the selected front. In the simplest setting, consider initial conditions wedged between a sub- and a super-solution, given, for an arbitrary small $\eps>0$, through
\begin{equation}\label{e:subsuper}
u_\mathrm{sup}(x)=u_*(x;c_*+\eps),\qquad u_\mathrm{sub}(x)=
\left\{\begin{array}{ll}
u_*(x;c_*-\eps),& x<x_*,\\
0,& x\geq x_*,
\end{array}
\right.
\end{equation}
where $x_*=\inf\{x|u_*(x;c_*-\eps)=0\}$ is the left most zero of the non-positive front. One verifies that $u_\mathrm{sup}$ and $u_\mathrm{sub}$ are indeed super- and sub-solutions that guarantee that an invasion speed defined for instance through
\begin{equation}\label{e:cinv}
c_\mathrm{inv}=\lim_{t\to\infty} \frac{\sup\{x\mid u(t,x)>\delta\}}{t},
\end{equation}
is the speed $c_*$ of the slowest positive front; see Fig.~\ref{f:subsuper} for an illustration.

\textbf{Selection of the slowest stable front --- diffusive stability.}
A different perspective on selection is gained when considering classes of initial conditions that arise through cutting off the tail of the front, such as $u_0(x)=u_*(x;c)$ for $x<L$ and $u_0(x)=0$ for $x>L$ for some $L\gg 1$. Since $u_*(L;c)=\rmO(\rme^{-\eta L})$ for some $\eta$, one may hope that this initial condition could be viewed as a small perturbation of the front. Attempting to control the evolution of this perturbation of the initial condition as $t\to\infty$, one would then expect to identify the correct speed $c$ of the selected front.

In the case of a pushed front, the front tail decays with rate $\rme^{\nu_-x}$, with $\nu_-<-\frac{c}{2}$, while we can choose an exponential weight that enforces decay with rate $\rme^{-\frac{c}{2} x}$ in order to stabilize the essential spectrum. As a consequence, a cutoff is a small perturbation $\rmO(\rme^{(\nu_-+\frac{c}{2})L})$ in this weighted norm. Since the linearized equation gives exponential decay up to a simple eigenvalue associated with translations, one can readily control the nonlinear evolution of disturbances up to a fixed, small shift in the position of the front; see for instance~\cite{Sattinger}.

For a pulled front, this reasoning fails as the essential spectrum of the front touches the origin. The cutoff induces a perturbation that is \emph{not small} in the optimally weighted space, in fact growing linearly in $x$, while the linear equation exhibits only algebraic decay $t^{-3/2}$. We will describe in more detail in \S\ref{s:nlmsp} how one can nevertheless control this perturbation by writing the solution as a superposition of a front profile shifted logarithmically and a Gaussian wave packet in the leading edge; see \S\ref{s:nlmsp} for an outline,~\cite{EbertvanSaarloos} for matched asymptotics, and~\cite{as1,avery2} for proofs.

We emphasize that the perturbative perspective, based on linear information, yields \emph{local} selection results, in a suitably defined neighborhood of a given invasion front, while the approach based on comparison principles yields, when available, more global stability results, including for instance actual Heaviside initial conditions.

\begin{figure}
\includegraphics[height=.8in]{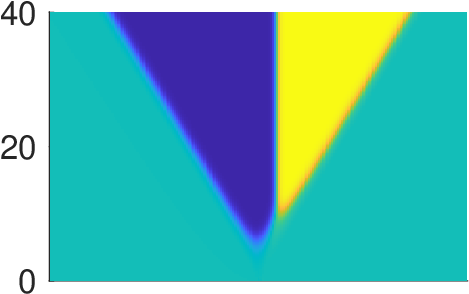}\hfill
\includegraphics[height=.8in]{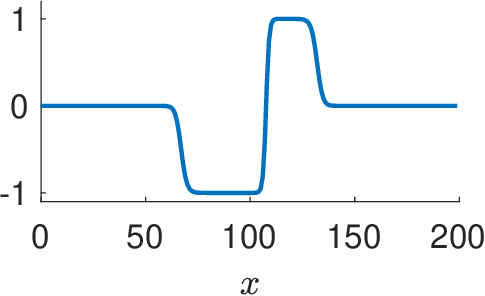}\hfill
\includegraphics[height=.8in]{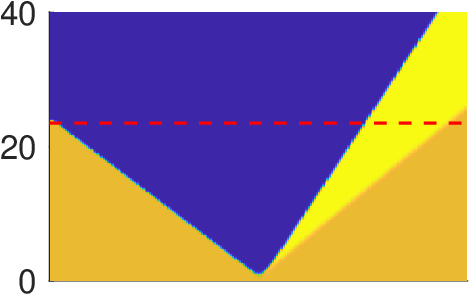}\hfill
\includegraphics[height=.8in]{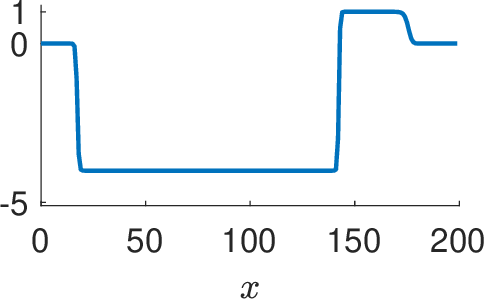}
\caption{Positive and negative fronts in the Nagumo equation~\eqref{e:nagumo} originating from a sign-changing initial condition at the center of the domain:  balanced nonlinearity $a=1$ (left two panels) and imbalanced nonlinearity $a=.2$ (right two panels). Space-time plots and snapshot at $t=20$ in both case. Note the slightly different speeds of propagation in the imbalanced case due to a pushed front propagating to the left in the right two panels (illustrated by the red dashed line). Invasion fronts leave behind a stationary (left) or traveling (right) kink.}\label{f:ppselect}
\end{figure}

\textbf{Pulled fronts, pushed fronts, and the selection of states in the wake.}
The limitation to ``small perturbations'' introduced in the absence of a comparison principle is to some extent necessary. Indeed, the invasion can be mediated by different selected fronts, depending on the initial condition.
Consider therefore the example of the Nagumo equation with $f(u)=u(1-u)(u+a)$, $0<a<1/2$. In addition to the front with $u_*(x)\to 1$ for $x\to-\infty$, there exists a front with $u_*(x;c)\to -a$ for $x\to -\infty$. Positive, step-like initial conditions lead to invasion at the pushed speed $c_\mathrm{push}=\frac{1+2a}{\sqrt{2}}$, while negative step-like initial conditions lead to invasion at the pulled speed $c_\mathrm{lin}=2\sqrt{a}$. In fact, a positive perturbation in the leading edge of a negative, step-like initial condition may lead to the selection of the positive front and the introduction of a kink in the wake; see for instance  Fig.~\ref{f:ppselect}. Global stability results without sign conditions can then clearly not hold. Rather, the leading edge of the invasion process selects not only a speed but also a state in the wake of the invasion. In the example of the Nagumo equation, the selected state is either $1$ with speed $c_\mathrm{push}$ or $-a$ with speed $c_\mathrm{lin}$.
\subsection{Pattern-forming fronts}
The selection of states in the wake is of particular interest in systems that admit many possible stable and unstable states, particularly pattern-forming equations with typically continua of equilibria. A simple example is a complex version of the Nagumo equation, known as the Ginzburg-Landau equation (CGL),
\begin{equation}\label{e:cgl}
A_t=A_{xx}+A-A|A|^2, \qquad A\in\C,
\end{equation}
with a family of stationary solutions $A_k(x;k)=\sqrt{1-k^2}\rme^{\rmi k x}$, $|k|<1$, that are stable for $|k|<1/3$. Restricting to real initial conditions, one finds the scalar equation discussed above, which leads to predicting propagation at speed $c=2$ and selection of the state $A_k(x;0)=1$ in the wake. It is conjectured that this selection happens for most, possibly complex, step-like initial conditions, up to possibly a complex rotation $\rme^{\rmi\varphi}$. Numerous results are concerned with existence and stability of front solutions in this equation with the strongest and most recent results establishing sharp linear decay estimates of the front propagating at the linear speed, but not convergence for step-like initial conditions; see~\cite{as5} and references therein.

\begin{figure}
\begin{minipage}{0.245\textwidth}
\centering
\includegraphics[height=1.5in]{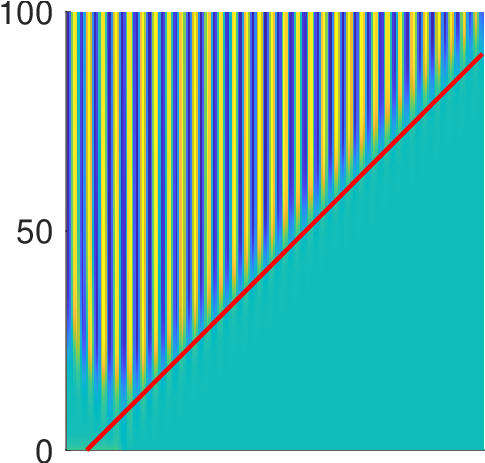}\\[0.1in]
\hspace*{0.07in}\includegraphics[height=0.90in]{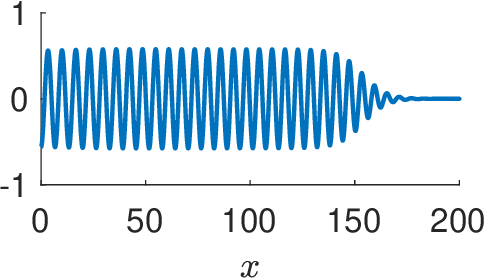}
\end{minipage}
\begin{minipage}{0.245\textwidth}
\centering
\includegraphics[height=1.5in]{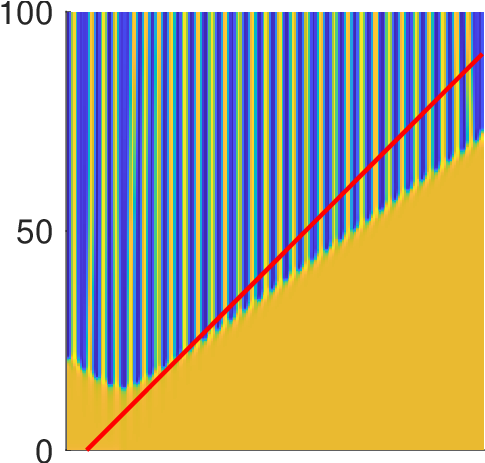}\\[0.1in]
\hspace*{0.07in}\includegraphics[height=0.90in]{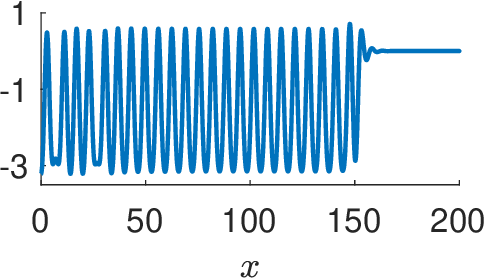}
\end{minipage}\hfill
\begin{minipage}{0.245\textwidth}
\centering
\includegraphics[height=1.5in]{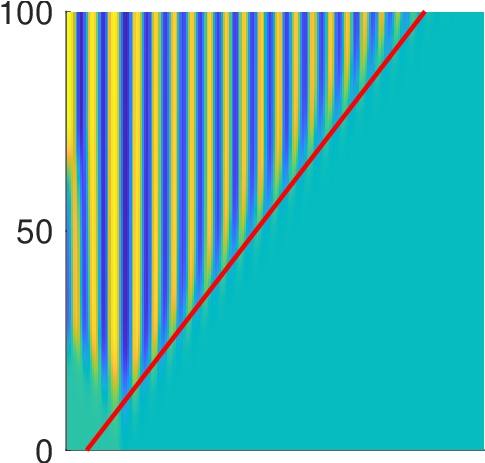}\\[0.1in]
\hspace*{0.07in}\includegraphics[trim={0.in 0in 0.in 0in},clip,height=0.90in]{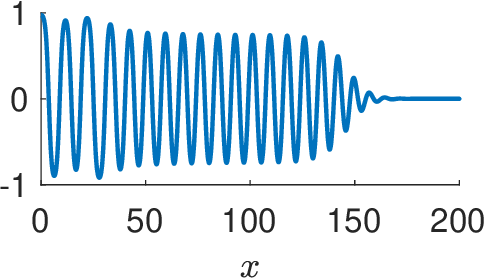}
\end{minipage}\hfill
\begin{minipage}{0.245\textwidth}
\centering
\includegraphics[height=1.5in]{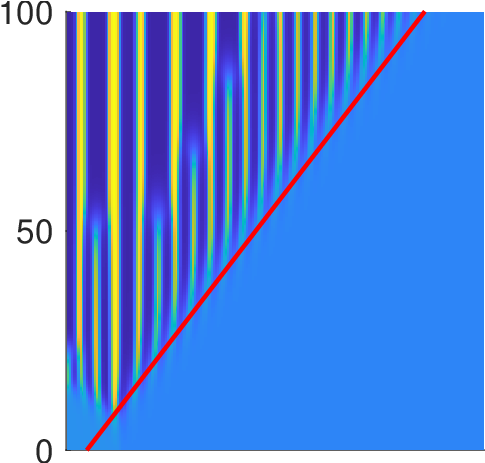}\\[0.1in]
\hspace*{0.0in}\includegraphics[height=0.87in]{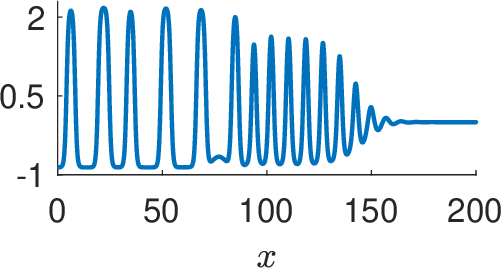}
\end{minipage}
\caption{Space-time plots (top row, time vertical) of pattern-forming fronts in SH~\eqref{e:sh}: pulled front (left, $\gamma=0,\eps=0.4$) and pushed front (center left, $\gamma=3,\eps=0-.4$); perturbations of $u\equiv 0$ in CH~\eqref{e:ch} with $\gamma=0$ (center right) and $\gamma=2$ (right). Red lines indicates linear spreading speeds from~\eqref{e:4cth4}. Note the faster, pushed, speed in the center left and the secondary coarsening in the right panel. Profiles of solution $u(t,x)$ (bottom row) at times $t=75,55,90,90$ (left to right).
%
}
\label{f:shch}\end{figure}

The Ginzburg-Landau equation arises in many situations as an amplitude and modulation equation near the onset of instability, in many contexts ranging from fluid dynamics to developmental biology; see for instance~\cite{schneideruecker,mielke,hoyle,crosshohenberg}. The explicit form of solutions $A(x;k)$ and the gauge invariance $A\mapsto \rme^{\rmi\varphi}A$ are due to an averaging symmetry that is generally absent in pattern-forming systems. The simplest example without such a gauge invariance is the Swift-Hohenberg equation,
\begin{equation}\label{e:sh}
u_t=-(\partial_{xx}+1)^2 u + \eps^2 u  +\gamma u^2 -u^3,\qquad u\in\R,
\end{equation}
with $\eps\gtrsim 0$ which we consider first for $\gamma=0$. Solutions are then well approximated by solutions to a scaled version of~\eqref{e:cgl} through $u(t,x)=\eps (A(\eps^2 t,\eps x)\rme^{\rmi x}+\bar{A}(\eps^2 t,\eps x)\rme^{-\rmi x})$. The trivial solution $u\equiv 0$ is unstable and there exists a family of stable periodic patterns which are  solutions of the form  $u_\mathrm{p}(x+\varphi;k)$, $u_\mathrm{p}(x+\frac{2\pi}{k};k)=u_\mathrm{p}(x;k)$, $k\sim 1, \varphi\in\R$. Invasion is conjectured at a linear spreading speed $c_\mathrm{lin}=4\eps+\rmO(\eps^2)$. Existence and stability of fronts is known for $c>c_\mathrm{lin}$ bur neither stability nor selection of the marginally stable front are known; see~\cite{es1,ce5} for some stability results and upper bounds on the speed of propagation, and references therein. In addition to the invasion speed, the linear analysis also predicts an invasion frequency $\omega_*$ and an associated wavenumber $k_*=\omega_*/c_\mathrm{lin}$ of patterns $u_\mathrm{p}(x;k)$ formed in the wake of the invasion process; see Fig.~\ref{f:shch} for an illustration of the pattern-forming process in the wake of the front.
We emphasize that the selection here appears to be largely independent of the choice of initial condition, in contrast to the example of the Nagumo equation where roughly a sign of the perturbation in the leading edge determines selected state and speed.

More complex invasion dynamics arise in the Cahn-Hilliard equation,
\begin{equation}\label{e:ch}
u_t=-(u_{xx}+u+\gamma u^2-u^3)_{xx}, \qquad u\in\R,
\end{equation}
where $u(x)\equiv\bar{u}\in(-\frac{1}{\sqrt{3}},\frac{1}{\sqrt{3}})$ is unstable. Linear predictions for the spreading speed match well with numerical simulations and existence of fronts is known; see~\cite{S2017}. Fronts select periodic patterns $u_\mathrm{p}(x;k)$ in the wake which however all are unstable, leading to subsequent coarsening. Linear predictions for the wavenumber appear to fail for values $|\bar{u}|\lesssim \frac{1}{\sqrt{3}}$ when $\gamma=0$; see Fig.~\ref{f:shch} for an illustration.

\begin{figure}
\begin{minipage}{0.245\textwidth}
\centering
\includegraphics[height=1.5in]{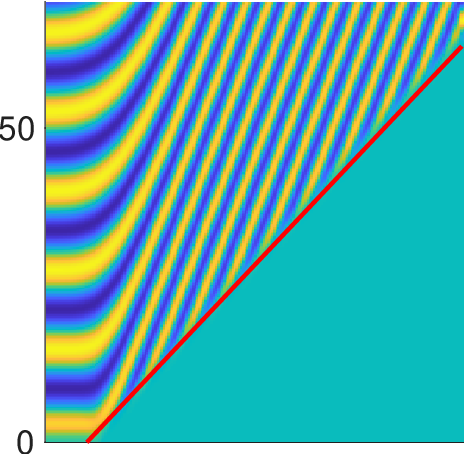}\\[0.1in]
\hspace*{0.07in}\includegraphics[height=0.90in]{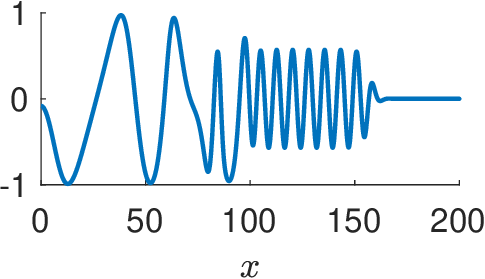}
\end{minipage}\hfill
\begin{minipage}{0.245\textwidth}
\centering
\includegraphics[height=1.5in]{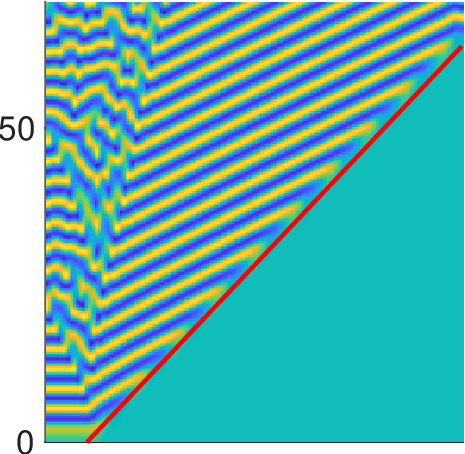}\\[0.1in]
\hspace*{0.07in}\includegraphics[height=0.90in]{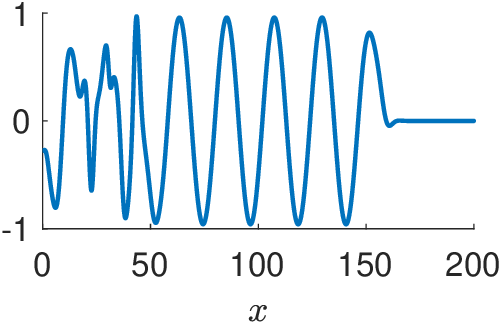}
\end{minipage}
\begin{minipage}{0.245\textwidth}
\centering
\includegraphics[height=1.5in]{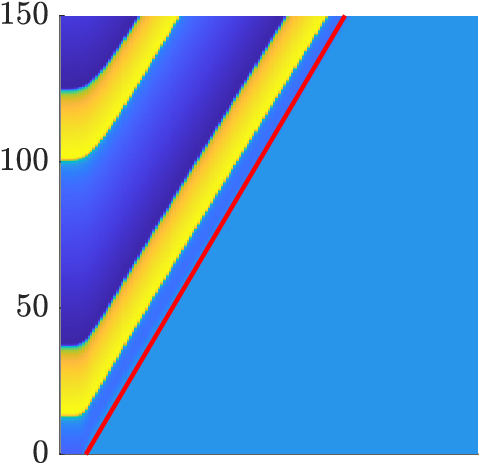}\\[0.1in]
\hspace*{0.07in}\includegraphics[height=0.90in]{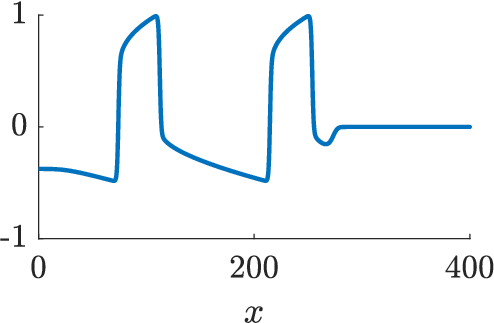}
\end{minipage}\hfill
\begin{minipage}{0.245\textwidth}
\centering
\includegraphics[height=1.5in]{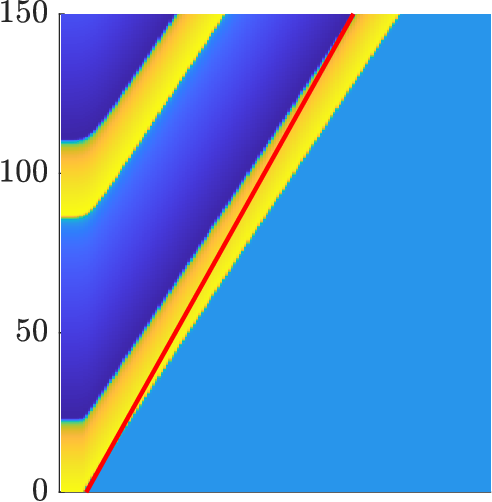}\\[0.1in]
\hspace*{0.07in}\includegraphics[height=0.90in]{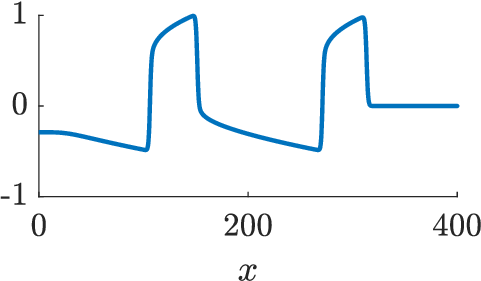}
\end{minipage}
\caption{Space-time plots (top row, time vertical) of pattern-forming fronts in CGL~\eqref{e:ccgl} and in FHN~\eqref{e:fhn}.  CGL with $\alpha=1$~\eqref{e:ccgl} showing $\Re A$: coherent invasion (left, $\beta=0.5$),  and secondary chaotic invasion  (center left, $\beta=-2$); FHN with $a=-0.2$, $b=\gamma=0$ showing a pulled front from positive perturbations (center right) and a pushed front from negative perturbations (right).
Red lines indicates linear spreading speeds from~\eqref{e:ccgllin} and~\eqref{e:fhnspsp2}.
Note the faster, pushed, speed in the right panel.
Profiles of solution $u(t,x)$ (bottom row) at times $t=60,60,150,150$ (left to right).
%
}
\label{f:cglfhn}
\end{figure}

Rather than spatially periodic, time-independent patterns, oscillatory media support temporally periodic solutions of the form $u_k(kx - \omega(k) t;k)$. Examples are the complex Ginzburg-Landau equation~\cite{aransonkramer}
\begin{equation}\label{e:ccgl}
A_t=(1+\rmi\alpha)A_{xx}+A-(1+\rmi\beta)A|A|^2,\qquad A\in\C,
\end{equation}
and the FitzHugh-Nagumo system (FHN)~\cite{tysonkeener}
\begin{equation}\label{e:fhn}
\begin{aligned}
u_t&=u_{xx}+u(1-u)(u-a) -v,\nonumber\\
v_t&=\eps(u-\gamma v + b),\nonumber
\end{aligned}
\end{equation}
with $0<\eps\ll 1$, $0<a<1$, and $b\geq 0$. In both cases, predictions for linear spreading speeds and spreading frequencies can be derived, and existence of invasion fronts can be established, at least in some parameter regimes. The FitzHugh-Nagumo system in fact exhibits two invasion fronts, one of which may be pushed~\cite{cartersch}. For pulled fronts, one can  predict wavenumbers in the wake of the invasion mostly correctly from the linearization; see Fig.~\ref{f:cglfhn} for an illustration.

We return to pattern-forming fronts in \S\ref{sec: modulated}, mimicking the strategy presented in \S\ref{s : FKPPintro}, from existence of fronts through linear marginal stability to selection.

\subsection{Outline}
{\color{black} In \S\ref{s:lms} through \S\ref{sec: modulated}, we focus on the role of marginal stability as a front selection mechanism.} In \S\ref{s:lms}, we explain how to predict linearized invasion speeds, and how this relates to a linearized marginal stability criterion, concluding with a discussion on practically computing linear spreading speeds, frequencies, and wavenumbers. We discuss the nonlinear invasion problem in \S\ref{s:nlms}, explaining strategies for finding front solutions and developing nonlinear marginal stability criteria which determine selected invasion speeds. In \S\ref{s:robust}, we discuss robustness of invasion processes under changes in parameters; in particular, this motivates the development of efficient numerical methods for continuing front speeds and detecting pushed-to-pulled transitions.
In \S\ref{s:nlmsp}, we formulate precise spectral conditions which rigorously guarantee pushed or pulled front invasion, and sketch the proofs of these selection results. We return to pattern-forming fronts in \S\ref{sec: modulated}, describing existing results on sharp linear stability in connection with open problems and strategies for understanding selection in this context. This section also includes practical strategies for finding selected states in the wake as well as problems related to secondary invasion processes that are initiated when the primary invasion selects an unstable pattern, as seen for instance in the Cahn-Hilliard equation. We conclude with a discussion, \S\ref{s:dis} that points to many open problems, revisiting in particular the connections with marginal positivity and steepest fronts.

\section{Linear marginal stability}\label{s:lms}

In \S\ref{s : FKPPintro}, we noted that for the Fisher-KPP (FKPP) equation, the spreading speed $c = 2$ for step-like initial conditions corresponds to the speed at which the associated front solution is marginally stable. In fact, more is true: the spreading speed $c = 2$ is precisely the speed at which we observe marginal pointwise stability in the linearized equation $u_t = u_{xx} + c u_x + u$. This further corresponds to the propagation speed of a fixed level set of the Green's function for the linearized equation, and we therefore refer to this speed as the {\em linear spreading speed}. In this section, we explain this correspondence and outline how to compute linear spreading speeds in general systems of constant-coefficient parabolic PDEs. 

\textcolor{black}{We start with some background on pointwise growth in constant-coefficient problems in \S\ref{s:2.1} and define pointwise spreading speeds through marginal pointwise growth in a fastest frame of reference. We relate the somewhat cumbersome pointwise growth criteria to algebraic double root computations in \S\ref{s:2.2} and to optimal choices of exponentially weighted spaces in \S\ref{s:2.3}.
The methods described in these section originate in the study of absolute and convective instabilities in problems from plasma physics; see~\cite{bers1983handbook,briggs} and~\cite{brevdo1,huerre90} for subsequent extensions to fluid dynamics problems.  We indicate how these linear criteria can predict patterns in the wake of invasion processes, with many caveats, in \S\ref{s:linpatt}. We collect alternate view points and consequences of the pointwise analysis in  
\S\ref{s:2.5} and conclude with a practical guide to the computation of pointwise growth modes and spreading speeds in \S\ref{s:practical_linearspreading} and many examples, some not documented in the literature previously, in \S\ref{s:2ex}. }


\subsection{Linear spreading speeds via pointwise growth modes}\label{s:2.1}
\textbf{General set-up.}
In this section, we will be concerned with systems of the form
\begin{equation}\label{e:lin}
u_t=\mathcal{P}(\partial_x)u, \qquad x\in\R,\quad u\in\R^N,\quad u(0,x)=u_0(x),
\end{equation}
with $\mathcal{P}(\nu)$ a matrix-valued polynomial of degree $2m$ such that the equation is parabolic and well-posed, that is, $\mathcal{P}(\nu)=d_{2m}\nu^{2m}+\rmO(\nu^{2m-1})$, $d_{2m}(-1)^{m}<0$. In order to track the solution $u(t,x)$ in a moving frame $x=\xi+ct$, we consider more generally,
\begin{equation}\label{e:linc}
u_t=\mathcal{P}(\partial_\xi)u+c\partial_\xi u, \qquad \xi\in\R,\quad u\in\R^N,\quad u(0,\xi)=u_0(\xi).
\end{equation}
Denote by  $\mathcal{L}=\mathcal{P}(\partial_\xi)+c\partial_\xi$ the operator appearing on the right hand side of (\ref{e:linc}). When considered as an unbounded, closed operator on spaces with translation invariant norms such as $X=L^p(\R,\R^N)$ or $X=BC^0_\mathrm{unif}(\R,\R^N)$, its spectrum is continuous and can be expressed as
\begin{equation}
\mathrm{spec}\,(\mathcal{L})= \{\lambda\in \mathbb{C} \ | \ d_c(\lambda,ik)=0 \ \text{for some $k\in\mathbb{R}$} \}
,\label{e:specdef}
\end{equation}
where $d_c(\lambda,\nu)$ is the dispersion relation in the co-moving frame from an ansatz $u(t,x)=\rme^{\lambda t +\nu \xi}\bar{u}$,
\begin{equation} d_c(\lambda,\nu)=\mathrm{det}\left( \mathcal{P}(\nu)+c\nu\mathrm{I}-\lambda I \right)\textcolor{black}{=d_0(\lambda-c\nu,\nu)}. \label{eq:disp} \end{equation}
Our focus is on unstable states $d_0(\lambda,\rmi k)=0$ for some $\Re\lambda>0$ and $k\in\R$. The norm $\|\rme^{\mathcal{L}t}\|_{X\to X}$ then grows exponentially for any wave speed $c\geq 0$. In fact, $\max\{\Re\,\mathrm{spec}\,(\mathcal{L})\}$ is independent of $c$ as is easily seen by shifting $\lambda\mapsto\lambda+c\rmi k$ when passing to a comoving frame.

\begin{figure}
    \centering\includegraphics[width=0.8\textwidth]{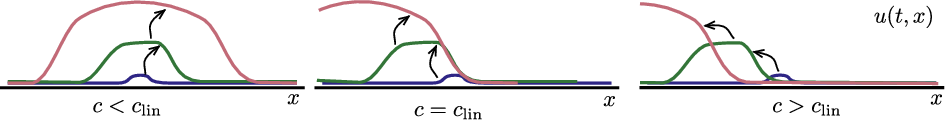}
    \caption{Small initial conditions grow in time due to the instability: in the stationary frame (left), the solution grows exponentially at each fixed point, while also spreading spatially; in a fast comoving frame, the solution decays exponentially at each fixed point, since the growth is advected to the left (right); for a critical intermediate speed, the growth is marginal since the frame tracks the leading edge of the temporal growth and spatial spreading (center). }\label{f:convabs}
\end{figure}
\begin{figure}[b]
    \centering\includegraphics[width=0.6\textwidth]{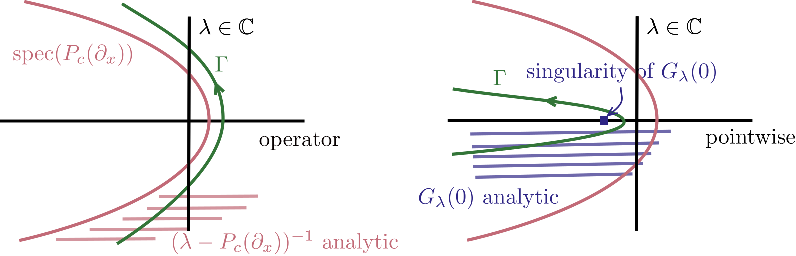}
    \caption{Construction of the temporal solution using the integration along $\Gamma$ over the resolvent: in function space (left) the contour can be deformed up to the spectrum of $\mathcal{L}$ which yields singularities of the resolvent as an operator; pointwise in $x$ with compactly supported initial conditions, obstructions are singularities in $\lambda$ of the resolvent kernel $G_\lambda(\xi)$, in our situation strictly to the left of the spectrum.  }\label{f:laplace}
\end{figure}

\textbf{Pointwise stability.}
Suppose that initial data $u_0(x)$ for (\ref{e:linc}) is compactly supported.  We are interested in quantifying how quickly the resulting solution will spread spatially.  Informally, this can be viewed from the perspective of a moving observer with the linear spreading speed being the slowest speed at which the observer outruns in the instability; see Fig.~\ref{f:convabs}.
A mathematical determination of this speed requires the calculation of pointwise exponential growth/decay rates as a function of the speed of the moving frame.  The primary aim of the remainder of this section is to explain how pointwise growth/decay rates can be computed from singularities of the pointwise Green's function which, in turn, may be computed by locating repeated roots of $d_c(\lambda, \nu)$ satisfying a ``pinching condition''.
Before proceeding, recall that the operator $\mathcal{L}$ generates an analytic semigroup on $L^2(\mathbb{R},\mathbb{R}^N)$ which may be defined using the inverse Laplace transform of the resolvent operator:
\[ \rme^{\mathcal{L}t}=\frac{-1}{2\pi i}\int_\Gamma  \rme^{\lambda t}(\mathcal{L}-\lambda)^{-1}\mathrm{d}\lambda.\]
Here, $\Gamma$ is a sectorial contour, oriented counterclockwise, and lying to the right of $\mathrm{spec}\,(\mathcal{L})$; see Fig.~\ref{f:laplace}. {\color{black} As mentioned above, for unstable background states the operator norm of $e^{\mathcal{L} t}$ induced by a fixed translation-invariant norm will grow exponentially.}

To obtain pointwise growth rates one must adopt a somewhat different perspective.  The resolvent operator appearing in the definition of $\rme^{\mathcal{L}t}$ can be expressed as an integral operator:
\begin{equation} (\mathcal{L}-\lambda)^{-1}u_0=\int_{-\infty}^\infty G_\lambda(\xi-y)u_0(y)\mathrm{d}y.\label{e:resolvent} \end{equation}
The function $G_\lambda(\xi-y)$ giving the resolvent kernel is called the pointwise Green's function; see~\cite{ZumbrunHoward}.  The advantage of this perspective is that the resolvent kernel may have an analytic continuation into the spectrum $\mathrm{spec}\,(\mathcal{L})$ whereas such an extension does not exist for the resolvent operator directly (for instance, if $G_\lambda(\xi)$ is pointwise analytic in $\lambda$ but is exponentially growing as $|\xi|\to\infty$); see again Fig.~\ref{f:laplace}.
Utilizing the pointwise Green's function and the fact that $u_0(x)$ is compactly supported, the solution to (\ref{e:linc}) can be expressed via the inverse Laplace transform as
\begin{equation} u(t,\xi)=\frac{-1}{2\pi i}\int_\Gamma \int_{-\infty}^\infty \rme^{\lambda t}G_\lambda(\xi-y)u_0(y)\mathrm{d}y\mathrm{d}\lambda, \label{e:ubyconvoluation} \end{equation}
where $\Gamma$ is a contour initially lying to the right of $\mathrm{spec}\,(\mathcal{L})$.  However, {\color{black} after swapping the order of integration,} Cauchy's integral theorem allows for the deformation of the inversion contour $\Gamma$ into the set $\mathrm{spec}\,(\mathcal{L})$ so long as singularities of $G_\lambda$ are avoided.  Optimal exponential growth/decay bounds are obtained by determining how far the Laplace inversion contour can be deformed while avoiding these singularities.

We therefore turn our attention to the location of singularities of the pointwise Green's function.   The pointwise Green's function can be constructed as a solution of the system of ODEs
\begin{equation} \label{eq:GlODE} \delta(\xi-y)\mathrm{I}=\mathcal{P}(\partial_\xi)G_\lambda+c\partial_\xi G_\lambda-\lambda G_\lambda ,\end{equation}
with $\lambda$ initially lying to the right of $\mathrm{spec}\,(\mathcal{L})$.  To obtain a general expression for $G_\lambda$ we find it convenient to express (\ref{eq:GlODE}) as a system of first order equations, setting $U=(u,\partial_\xi u,\ldots \partial_\xi^{2m-1}u)\in \C^{2mN}$, and solve for the Green's function for this first order equation via
\begin{equation} 
U_\xi=M_\lambda U+\delta(\xi-y) \mathrm{I}, \label{e:Ufirstorder} 
\end{equation}
where the $\lambda$-dependence of the matrix $M$ is linear, appearing in the bottom left entry, only.

In order to find bounded solutions to~\eqref{e:Ufirstorder}, we need to consider eigenvalues $\nu$ of $M_\lambda$, which we refer to in the sequel as \emph{spatial eigenvalues}, as opposed to spectral values $\lambda$ that are \emph{temporal eigenvalues}.  For $\lambda$ to the right of the essential spectrum ellipticity of the operator $\mathcal{P}(\partial_\xi)$ implies that the stable and unstable eigenspaces of $M_\lambda$ will have equal dimension.
For such $\lambda$ values, denote by $P^s_\lambda$ and $P^u_\lambda$  the complementary stable and unstable spectral projections for the matrix $M_\lambda$, that is, the projections on generalized eigenspaces $E^\mathrm{s}_\lambda$ and $E^\mathrm{u}_\lambda$  to spatial eigenvalues with negative and positive real parts, respectively.    Imposing decay as $\xi\to\pm\infty$, the equation (\ref{e:Ufirstorder}) can be directly solved as
\[ U(\xi)=\int_{-\infty}^\xi \rme^{M_\lambda (\xi-z)}P_\lambda^s \delta(z-y)\mathrm{d}z -\int_\xi^\infty \rme^{M_\lambda (\xi-z)}P_\lambda^u \delta(z-y)\mathrm{d}z \]
from which we obtain that the Green's function for the first order equation is given through
\begin{equation}\label{e:vecgreens} \mathcal{T}_\lambda(\xi-y) = \left\{ \begin{array}{cc} \rme^{M_\lambda (\xi-y)}P_\lambda^s & \xi\geq y \\ -\rme^{M_\lambda (\xi-y)}P_\lambda^u & \xi < y\end{array}\right. . \end{equation}
An expression for $G_\lambda(\xi-y)$ can be extracted from the entries of $\mathcal{T}_\lambda(\xi-y)$ and, importantly,  the domains of analyticity of both functions are identical; we refer to~\cite[\S2.3, Lemma 2.1]{HolzerScheelPointwiseGrowth} for details.

From~\eqref{e:vecgreens}, we see immediately that analyticity of $\mathcal{T}_\lambda$, and equivalently of $G_\lambda$, is determined by the analyticity of the spectral projections $P^{s/u}_\lambda$.  To the right of the essential spectrum these projections are analytic in $\lambda$:  stable and unstable spectral sets can be separated by a contour in the complex plane so that projections can be obtained using a Dunford integral over $\nu$ of the resolvent $(\nu-M_\lambda)^{-1}$, which depends analytically on $\lambda$.

\begin{figure}
    \centering\includegraphics[width=0.45\textwidth]{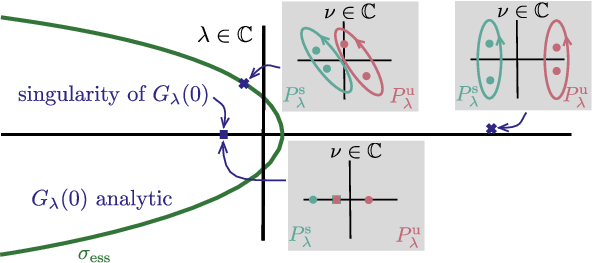}\hspace*{.3in}\includegraphics[width=0.45\textwidth]{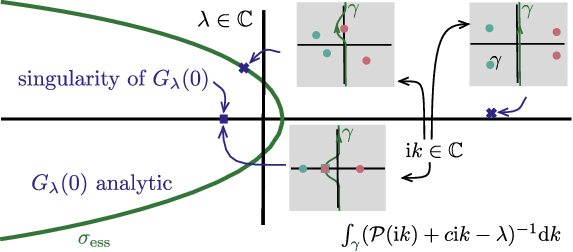}
    \caption{Singularities in the temporal eigenvalue parameter $\lambda$ of the pointwise Green's function induced by collision of spatial eigenvalues $\nu$ (left). Away from the singularities, Dunford integral contours can be used to continue the projections $P^\mathrm{s/u}(\lambda)$ in an analytic fashion starting at $\lambda\gg 1$. Also shown, the construction of the resolvent kernel via Fourier transform, possibly transforming the Fourier integration contour $\gamma$ into the complex $k$-plane (right).}\label{f:spatialev}
\end{figure}

The analysis thus far has concentrated on $\lambda$ values to the right of the spectrum $\mathrm{spec}\,(\mathcal{L})$.  By the construction above, analytic continuation of $G_\lambda$ into $\mathrm{spec}\,(\mathcal{L})$ depends only on our ability to analytically continue the spectral projections $P^{s/u}_\lambda$ into the spectrum.
This is possible so long as the spectral sets, that is, collections of eigenvalues defining these two projections, do not intersect. Conversely, intersections of these spectral sets typically lead to {\color{black} singularities} of the pointwise Green's function $G_\lambda$ and prevent further continuation of the Laplace inversion contour in (\ref{e:ubyconvoluation}); see Fig.~\ref{f:spatialev}. We emphasize that for us $P^{s/u}_\lambda$ always refer to the analytic continuations of the projections from $\mathrm{Re} \, \lambda \gg 1$; for a fixed $\lambda$ inside the essential spectrum, these may no longer correspond to the stable/unstable subspaces for $M_\lambda$; see again Fig.~\ref{f:spatialev}.

In this way, optimal pointwise growth/decay rates are obtained by consideration of the singularities of the analytically continued pointwise Green's function, which leads to the following definition.
\begin{definition}[Pointwise growth modes]
    The complex number $\lambda_*$ is a pointwise growth mode if $P_\lambda^s$ is not analytic as a function of $\lambda$ in a neighborhood of $\lambda_*$.  The pointwise growth rate of (\ref{e:linc}) is
    \[ \inf \left\{ \rho\in \mathbb{R} \ | \ P_\lambda^s \ \text{is analytic for } \ \mathrm{Re}(\lambda)>\rho \right\}.  \]
\end{definition}
One can show that pointwise growth modes give sharp bounds on pointwise exponential growth rates~\cite[Cor. 2.3]{HolzerScheelPointwiseGrowth}.
Pointwise growth rates can be used to differentiate between absolute and convective instabilities; see~\cite{bers1983handbook,briggs,huerre90,AbsoluteSpecArndBjorn}.  If the pointwise growth rate is positive, then the instability is {\em absolute} in the sense that the perturbation of the unstable state will grow both in norm and pointwise.  If the pointwise growth rate is negative then the instability is referred to as {\em convective} since pointwise decay is observed while the solution grows exponentially when measured in translation invariant norms.  When the pointwise growth rate is zero then the system is {\em marginally stable} and neither pointwise exponential growth nor decay is observed.

\begin{figure}
    \centering    \includegraphics[width=\textwidth]{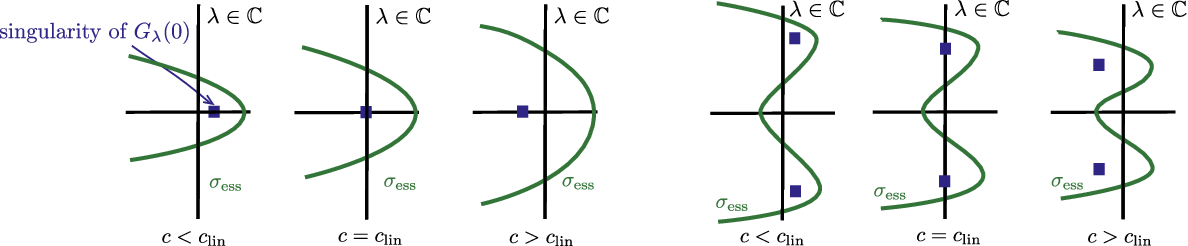}
    \caption{Pointwise growth modes (usually pinched double roots as explained in \S\ref{s:2.2}), cross the imaginary axis from right to left as $c$ is increased. The figure shows spectra and pointwise growth modes in the marginally stationary (left) and marginally oscillatory (right) case. At the linear spreading speed, the system is marginally pointwise stable.}\label{f:pgmclin}
\end{figure}

\textbf{Pointwise spreading speeds.}
We can now return to the original motivating problem for this section: the determination of the linear spreading speed.  The linear spreading speed characterizes how quickly compactly supported initial data spreads for the PDE (\ref{e:lin}).  To identify this speed, one can view this evolution in moving frames for a variety of speeds $c$. One then computes pointwise growth rates of (\ref{e:linc}) as a function of $c$.  If for some fixed $c$ the pointwise growth rate in this frame is positive then, by definition, $c$ provides a lower bound on the linear spreading speed.  Taking the supremum over all such lower bounds yields the linear spreading speed.
\begin{definition}[Spreading speeds]\label{def:lssfromPGR}
    The linear spreading speed associated to (\ref{e:linc}), $c_\mathrm{lin}$, is the supremum over all  speeds $c$ for which the pointwise growth rate of (\ref{e:linc}) is positive. 
\end{definition} 
\textcolor{black}{Using continuity properties of pointwise growth modes, one finds that in the frame with speed $c_\mathrm{lin}$, the system then is marginally stable, that is, its most unstable pointwise growth mode lies on the imaginary axis.} The linear spreading speed $c_\mathrm{lin}$ typically marks the transition between absolute and convective instabilities and occurs exactly when the pointwise Green's function has a singularity on the imaginary axis.
This crossing of pointwise growth modes is illustrated in Fig.~\ref{f:pgmclin}.
The supremum is finite for the parabolic equations that we are considering here: for large speeds, spatial eigenvalues $\nu=c^{1/{2m-1}}\hat{\nu}$ solve at leading order $-(-1)^m\hat{\nu}^{2m} +c\hat{\nu} =0$, and one finds $m$ roots with nonnegative real parts and $m$ roots with strictly negative real part. In other words, the spectrum of  the operator $\mathcal{L}$ in a space with exponential weight $\rme^{\eta x}$, $\eta=\delta c^{1/(2m-1)}$ is strictly negative for some small $\delta>0$ and we conclude that the spreading speed from Def.~\ref{def:lssfromPGR} is finite \cite[Lem.~6.3]{HolzerScheelPointwiseGrowth}. One can then find the linear spreading speed by decreasing the speed as a homotopy parameter until a pointwise growth mode appears on the imaginary axis; see also \S\ref{s:nlms} for an analogous procedure and more details. On the other hand, the maximum of $\Re\lambda$ in the spectrum gives rise to a pointwise growth mode in a frame moving with the group velocity $-\lambda'(\nu)$; see~\cite[Lem.~6.3]{HolzerScheelPointwiseGrowth}. As a consequence, the system is absolutely unstable for some speed and therefore the spreading speed is well-defined as a supremum and strictly positive in systems with reflection symmetry $x\to -x$.

We also note that without reflection symmetry, that is, if $\mathcal{P}$ is not even,  spreading speeds to the right may differ from those to the left and one should differentiate between left and right linear spreading speeds. In general, we may also have finitely many intervals of speeds where the system is absolutely unstable. The boundaries of these intervals are then spreading speeds in an appropriately refined sense. An example is the modulation equation near an oscillatory Turing instability~\cite{crosshohenberg},
\begin{equation}\label{e:coupcgl}
A_t=(1+\rmi\alpha)A_{xx}+c_\mathrm{g}A_x+A,\qquad B_t=(1+\rmi\alpha)B_{xx}-c_\mathrm{g}B_x+B,
\end{equation}
with spreading intervals $[c_\mathrm{g}-c_*,c_\mathrm{g}+c_*]\cup [-c_\mathrm{g}-c_*,-c_\mathrm{g}+c_*]$, where $c_*=2\sqrt{1+\alpha^2}$; see~\eqref{e:cglclin}.

\subsection{Pinched double roots and linear spreading speeds}\label{s:2.2}

The direct computation of pointwise growth modes is cumbersome since it requires an explicit calculation of the spectral projection matrices $P_\lambda^{s/u}$.  However, any loss of analyticity of these spectral projections can be traced to the coalescing of two eigenvalues: one corresponding to the stable projection and one corresponding to the unstable projection. At the collision, these eigenvalues are double roots of the dispersion relation (\ref{eq:disp}) and we are led to the following definition.
\begin{definition}[Pinched double roots]\label{def:PDR} A pair $(\lambda_*,\nu_*)$ is a  pinched double root if
\begin{itemize}
    \item[(i)] $d(\lambda_*,\nu_*)=0$ and $\partial_\nu d(\lambda_*,\nu_*)=0$;
    \item[(ii)] there is a continuous curve $\lambda(\tau)$, $\tau\geq 0$, $\Re\,\lambda(\tau)\nearrow$, $\lambda(0)=\lambda_*$ and $\Re\,\lambda(\tau)\to \infty$ as $\tau\to \infty$ and two roots $\nu_\pm(\lambda(\tau))$ with $\nu_\pm(0)=\nu_*$ and $\Re\,\nu_\pm(\lambda(\tau))\to \pm \infty$ as $\tau\to \infty$.
\end{itemize}
\end{definition}
The first of these conditions is a double root condition while the second is referred to as a pinching condition. The pinching condition is required so that the double root leads to a singularity of $G_\lambda$.
It enforces that these double roots form through the collision of one eigenvalue that originated in the stable eigenspace (associated with the spectral projection $P_\lambda^s$) for $\Re\,\lambda\to+\infty$ and one from the unstable eigenspace (associated to the projection $P_\lambda^u$).

Def.~\ref{def:PDR} leads to an alternate characterization of the linear spreading speed that is somewhat more tractable in practice than the one presented in Def.~\ref{def:lssfromPGR}.
\begin{lemma}[Pinched double roots and pointwise growth modes~\cite{HolzerScheelPointwiseGrowth}]\label{lem:lssfromPDR}
    For a given $\lambda\in\C$,
    \begin{itemize}
        \item $\lambda$ pointwise growth mode $\Longrightarrow$ $(\lambda,\nu)$ pinched double root for some $\nu$;
        \item $(\lambda,\nu)$ pinched double root $\Longrightarrow$ generically, $\lambda$ pointwise growth mode.
    \end{itemize}
    As a consequence, we can define a spreading speed for (\ref{e:lin})  associated with pinched double roots through
    \[ c_\mathrm{lin,dr}=\sup \left\{ c \ | \ d_c(\lambda,\nu) \ \text{has a pinched double root with} \ \mathrm{Re}(\lambda)>0  \right\}. \]
    We then have
    \[
    c_\mathrm{lin,dr}\geq c_\mathrm{lin},\qquad \text{and, generically, }
    \qquad c_\mathrm{lin,dr}= c_\mathrm{lin}.
    \]
\end{lemma}
For parabolic PDEs the dispersion relation will be a polynomial in the variables $\lambda$ and $\nu$ and so computing double roots, at least numerically, is straightforward.  The pinching condition is more difficult to verify; see \S\ref{s:practical_linearspreading} below.

We emphasize that marginally stable pinched double roots may overestimate the spreading speed since, as indicated in Lem.~\ref{lem:lssfromPDR}, they may not induce pointwise growth. A simple example occurs in a system of two uncoupled equations and is discussed in \eqref{e:cgl12lin};
see~\cite{HolzerScheelPointwiseGrowth} for more details and~\cite{HS2012}  for an example in a Lotka-Volterra system.


\subsection{Exponential weights and essential spectrum}\label{s:2.3}

Exponentially weighted function spaces play an important role in the analysis of invasion fronts in nonlinear systems. They enter the analysis in particular when one wishes to encode pointwise decay in function spaces in order to control nonlinearities, whose effects on the time evolution is not as explicitly controllable. We therefore pause to clarify the relationship between these spaces and the method for determining the linear spreading speed just explained. Exponentially weighted $L^p$ spaces were introduced in (\ref{e:Lpweighted}). The analysis of this section only requires a single weight so we consider $u\in L^2_{\mathrm{exp},\eta}(\mathbb{R})$ for some $\eta\geq 0$.  A function $u$ belongs to $L^2_{\mathrm{exp},\eta}(\mathbb{R}) $ if and only if $w=\rme^{\eta \xi}u$ is an element of $L^2(\mathbb{R})$. Statements throughout hold in $L^p$-based spaces, as well.

The resolvent operator,  given through (\ref{e:resolvent}) is bounded on $L^2(\mathbb{R})$ precisely when $G_\lambda \in L^1(\mathbb{R})$.  However, when $\lambda$ is continued into the spectrum $\mathrm{spec}\,(\mathcal{L})$ it is typically the case that the pointwise Green's function $G_\lambda$ loses this spatial localization, although $G_\lambda(\xi)$ is still analytic in $\lambda$. Indeed, $\mathrm{spec}\,(\mathcal{L})$ is characterized by those $\lambda$ values for which $d(\lambda, ik)=0$ for some $k\in\mathbb{R}$ and therefore $\lambda$ crossing through $\mathrm{spec}\,(\mathcal{L})$ generally corresponds to the crossing of a spatial eigenvalue from one side of the complex plane to the other leading to spatial exponential growth of the analytically extended Green's function $G_\lambda(\xi-y)$ on one half-line; compare also Fig.~\ref{f:spatialev}.

By working in an exponentially weighted space this localization can often be recovered.  For example, in order to bound  the resolvent operator acting on the space $L^2_{\mathrm{exp},\eta}(\mathbb{R})$ for some $\eta>0$, we need to bound
\begin{equation}
\rme^{\eta \xi}\left[(\mathcal{L}-\lambda)^{-1}u_0\right](\xi)=\int_{-\infty}^\infty \rme^{\eta (\xi-y)} G_\lambda(\xi-y)\left( \rme^{\eta y} u_0(y)\right)\mathrm{d}y \label{e:resolventeta}
\end{equation}
in $L^2$ in terms of $\|\rme^{\eta \cdot} u_0(\cdot)\|_{L^2}$. Inspecting the right-hand side, we then observe that $(\mathcal{L}-\lambda)$ is bounded invertible on $L^2_{\mathrm{exp},\eta}(\mathbb{R})$ if and only if $\rme^{\eta \cdot }G_\lambda(\cdot)\in L^1(\mathbb{R})$.   This integrability requirement implies that the spatial roots leading to the spectral projection $P_\lambda^s$ must all lie to the left of the line $\mathrm{Re}(z)=-\eta$ in the complex plane while the spatial roots lending to the definition of $P_\lambda^u$ must all lie to the right of this line. In other words, the  spectrum in $L^2_{\mathrm{exp},\eta}$ is given through
\begin{equation} \sigma_\eta(\mathcal{L})= \{\lambda\in \mathbb{C} \ | \ d_c(\lambda,ik-\eta)=0 \ \text{for some $k\in\mathbb{R}$} \}.\label{e:specdef2}
\end{equation}
Clearly, when interested in pointwise stability, one would like to optimize the chosen weight $\eta$. One simple approach would be to first for a given $\eta$, find the most unstable $\lambda$, that is, look for values of $\lambda$ so that $\Re\lambda$ is maximal in $\sigma_\eta$. To do this, one solves $d_c(\lambda,\rmi k -\eta)=0$ for $\lambda$ as a function of $\nu=\rmi k-\eta$ so that maximality implies
$
\frac{\rmd\,\Re\lambda}{\rmd\,\Im\nu}=0.
$ At the same time, one would look for an optimal choice of $\eta$ that minimizes that maximum, so that
$
\frac{\rmd\,\Re\lambda}{\rmd\,\Re\nu}=0.
$
By the Cauchy-Riemann equations for the complex analytic function function $\lambda(\nu)$ this then also implies $\frac{\rmd\,\Im\lambda}{\rmd\,\Re\nu}=0$ and  $\frac{\rmd\,\Im\lambda}{\rmd\Im\nu}=0$. Together, we see that the min-max condition {\color{black}
$
\min_\eta\max_k \Re \lambda(k-\rmi\eta)
$} leads to a solution of $\frac{\rmd\lambda}{\rmd\nu}=0$, satisfied at double roots $d(\lambda,\nu)=\partial_\nu(\lambda,\nu)=0$.

One comes to a similar conclusion when inspecting the stable and unstable spectral sets and attempting to separate the two by \emph{some} vertical line $\Re\nu=-\eta$ in the complex plane. One can show that instabilities always lead to $\Re\nu<0$ so that the weight in fact enforces decay in the direction of propagation~\cite[Remark 6.6]{HolzerScheelPointwiseGrowth}. Obstructions to this separation are the collision points between stable and unstable spatial roots $\nu$, precisely the pinched double roots described above, and more generally points where stable and unstable roots possess equal real part, a set in the $\lambda$-plane referred to as the absolute spectrum \textcolor{black}{when equipped with a Morse index condition}; see~\cite{AbsoluteSpecArndBjorn}. Those points $\lambda_*$ all form obstructions to finding a weight $\rme^{\eta \xi}$ so that $\rme^{\eta \cdot }G_\lambda(\cdot)\in L^1(\mathbb{R})$ for all $\Re\lambda\geq\Re\lambda_*$. \textcolor{black}{Commonly, pinched double roots indeed mark endpoints of curves of absolute spectrum but, unfortunately, not always: the pinching condition does not in general coincide with the Morse index condition, which requires an equal number of roots $\nu$ with real parts larger than and less than the double root. A simple counter example is the uncoupled system
\[
u_t=u_{xx}+2u_x+u,\qquad v_t= v_{xx} + \rmi v,\qquad \bar{v}_t= \bar{v}_{xx} - \rmi \bar{v}.
\]
In this example, one immediately finds spatial roots $\nu_u^\pm=-1\pm\sqrt{\lambda}$, $\nu_v^\pm= \pm\sqrt{\lambda-\rmi}$,  $\nu_{\bar{v}}^\pm= \pm\sqrt{\lambda+\rmi}$, with double roots at $\lambda=0$, $\lambda=\pm\rmi$ from collisions of roots $\nu_u^\pm$, $\nu_v^\pm$, and $\nu_{\bar{v}}^\pm$, respectively, and $\lambda=\pm\rmi/2$ from collisions of roots $\nu_v$ or $\nu_{\bar{v}}$ with roots $\nu_u$. The double root $\lambda=0$ is clearly pinched and a pointwise growth mode. At $\lambda=0$, the double root has $\nu=-1$, all other roots have $\Re\nu=\pm 1/\sqrt{2}>-1$, that is, the double root is not part of the absolute spectrum; see \cite{rss07} for background and strategies for computing absolute spectra in examples such as this one. }

Summarizing, suppose that the most unstable pinched double root is given by $(\lambda_*,\nu_*)$. Then, one has
\begin{equation}\label{e:pdrexp}
\max\Re\sigma_\eta(\mathcal{L})\geq \Re\lambda_*\text{ for all } \eta,\qquad  \text{ and, commonly,\quad  }
\Re\sigma_\eta(\mathcal{L})\leq \Re\lambda_*\text{ for }\eta=-\Re\nu_*.
\end{equation}
Examples where the second inequality does not hold have been analyzed in~\cite{FayeHolzerScheelSiemer}. In most examples, stability of spectrum in a carefully chosen exponentially weighted space does however give sharp conditions on pointwise stability.



\subsection{Linear pattern selection}\label{s:linpatt}
An imaginary pinched  double root $(\rmi\omega_\mathrm{lin},\nu_\mathrm{lin})$ at $c=c_\mathrm{lin}$ induces pointwise dynamics $\rme^{\rmi\omega_\mathrm{lin}t+\nu_\mathrm{lin}\xi}$ in the leading edge of the instability. For $\omega_\mathrm{lin}>0$, usually $\Im\nu_\mathrm{lin}=k_\mathrm{lin}\neq 0$, and  one finds spatio-temporal oscillations in the leading edge of the form $\rme^{\rmi(\omega_\mathrm{lin}t + k_\mathrm{lin} \xi)}\rme^{-\eta_\mathrm{lin} \xi}$,  with a super-imposed exponential decay. This suggests the creation of spatio-temporally periodic patterns derived from this linear prediction in the form
$ u(\omega_\mathrm{lin} t + k_\mathrm{lin} (x-c_\mathrm{lin} t))$, $u(\zeta)=u(\zeta+2\pi)$, in the wake of the leading edge. 
That is, a selected frequency $\omega_\mathrm{lin} - k_\mathrm{lin}c_\mathrm{lin}$ and spatial wavenumber $k_\mathrm{lin}$, with the caveat of exponentially decaying spatial modulations.

Notably a different wavenumber can be obtained by
%
%
%
%
tracking sign changes (or maxima) in the leading edge, expecting that each leads to a sign change in a nonlinear solution, a process referred to as node conservation~\cite{deelanger}. The oscillation frequency $\omega_\mathrm{lin}$ in the comoving frame leads to a spatial separation $L_\mathrm{node}=2\pi/k_\mathrm{lin}$, $k_\mathrm{node}=\omega_\mathrm{lin} /c_\mathrm{lin}$. 
Unfortunately, this prediction, while better than $k_\mathrm{lin}$ above, still needs refinement and may fail; see \S\ref{sec: modulated}, particularly \S\ref{s:modpract}.

\subsection{More about double roots and linear spreading speeds}\label{s:2.5}

\textbf{Variational criteria in order-preserving systems.}
For scalar equations obeying a comparison principle another approach can be taken to obtain the linear spreading speed.  For the system linearized near the unstable state, one seeks exponential solutions of the form $u(x,t)=\rme^{\nu x+\lambda t}$.  For each $\nu<0$, one can use the dispersion relation to obtain $\lambda=\lambda(\nu)$, so that these simple exponential solutions spread with speed $c(\nu)=-\frac{\lambda(\nu)}{\nu}$.   Due to the comparison principle these exponential solutions are automatically super-solutions which constrain compactly supported initial data and so $c(\nu)$ is an upper bound on the linear spreading speed.  One then defines
\[ c_\mathrm{lin}=\inf_{\nu<0} -\frac{\lambda(\nu)}{\nu}.\]
In more general contexts pertaining to inhomogeneous media this is referred to as the Freidlin-G\"{a}rtner formula; see~\cite{freidlingartner}.
Since $\lambda(\nu)$ is smooth we can obtain minimizers by simultaneously solving
\[ \lambda(\nu)+c\nu=0, \ -\frac{\lambda'(\nu)}{\nu}+\frac{\lambda(\nu)}{\nu^2}=0.\]
Multiplying the second equation by $-\nu$ and recalling that $c=-\frac{\lambda(\nu)}{\nu}$ this re-writes as
\[ \lambda(\nu)+c\nu=0, \qquad \textcolor{black}{\lambda'(\nu)+c=0,}\]
which implies that minimizers of the expression $c(\nu)$ are also double roots of the dispersion relation in a comoving frame.

\textbf{Double roots, group velocities, and phase velocities in scalar systems.} In a scalar system, $d(\lambda,\nu)=0$ if $\lambda=\lambda_\mathrm{st}(\nu)$, which we fix for now in the steady frame, that is, $\lambda_\mathrm{st}(\nu)=\mathcal{P}(\nu)$. When $\lambda,\nu\in\rmi\R$, one refers to $c_\mathrm{ph}=-\lambda/\nu$ as the phase velocity and to $c_\mathrm{g}=-\rmd \lambda/\rmd\nu$ as the group velocity.
In the frame moving with speed $c$, we then find $\lambda=\lambda_\mathrm{st}(\nu)+c\nu$, and the Galilean transformation $c_\mathrm{g}\to c_\mathrm{g}-c$, $c_\mathrm{ph}\to c_\mathrm{ph}-c$. Solving for the double-root criterion gives
\[
\textcolor{black}{\lambda_\mathrm{st}(\nu)+c\nu=\rmi\omega,\qquad \lambda_\mathrm{st}'(\nu)+c=0,
}\]
which can be mostly solved explicitly as
\[
\textcolor{black}{c=-\Re\lambda_\mathrm{st}'(\nu), \quad 0=\Im \lambda_\mathrm{st}'(\nu),\quad \omega=\Im \lambda_\mathrm{st}(\nu)+c\,\Im\nu,\quad 0=\Re \lambda_\mathrm{st}(\nu) +c\,\Re\nu.}
\]
The first and third equation define $c$ and $\omega$, the second and fourth equation translate into
\begin{equation}\label{e:spledr}
\frac{\rmd \mathcal{P}}{\rmd\nu}=\frac{\Re \mathcal{P}}{\Re\nu},
\end{equation}
which consists of two real equations, for real and imaginary part, in the two variables $\Re\nu$ and $\Im\nu$.
We emphasize however that this equation does not incorporate the pinching condition.

\textbf{Fourier-Laplace and the origin of pinching.}
The origin and meaning of the term ``pinching condition" in Def.~\ref{def:PDR} is not evident based upon our treatment here.  This phrase stems from the original analysis in the plasma physics literature where the solution to (\ref{e:linc}) is obtained using a combined Laplace (in time) and Fourier (in space) transform; see again~\cite{bers1983handbook,briggs,huerre90,vanSaarloosReview}.  In this setting, the goal remains the same: to deform the Laplace contour as far as possible towards the stable half of the complex plane so as to obtain optimal temporal estimates on the solution.  An immediate issue is encountered when the Laplace contour is deformed into the essential spectrum since one then encounters a singularity stemming from a root of the dispersion relation.  The Fourier contour must then also be deformed into the complex plane to avoid this singularity.  This simultaneous deformation continues until two curves of singularities collide (or pinch) and prevent further deformation of the contours.  The collision of these contours occurs exactly when there is a double root of the dispersion relation whose roots can be traced to opposite sides of the original Fourier contour.

The formulation in terms of Fourier-Laplace transform is more general than the one we have presented here in that it applies to more general evolution equations, including for instance nonlocal operators.  The Green's function approach that we present here is -- in our opinion -- a more straightforward route to motivating the double root criterion for computing linear spreading speeds. It also generalizes readily to variable-coefficient problems that we encounter in the discussion of nonlinear marginal stability in \S\ref{s:nlms}. Finally, it also forms the basis of a direct formulation for finding pinched double roots as spectral points of a generalized eigenvalue problem, thus enabling efficient iterative algorithms reminiscent of the inverse power method for regular eigenvalue problems that we briefly discuss below~\cite{S23}.


\textbf{Simple double roots and continuity.}
We say a double root is {\em simple} when $\partial_\lambda d(\lambda_*,\nu_*)\neq 0$ and $\partial^2_\nu d(\lambda_*,\nu_*)\neq 0$ and we say it is a {\em multiple} double root,  otherwise. At simple double roots, the linearization of the double-root system  with respect to $(\lambda,\nu)$ is invertible:
\[
f(\lambda,\nu)=\left(\begin{array}{c} d(\lambda,\nu)\\ \partial_\nu d(\lambda,\nu)\end{array}\right),\qquad \partial_{(\lambda,\nu)}f(\lambda,\nu)=\left(\begin{array}{cc} \partial_\lambda d& \partial_{\lambda\nu}d \\ \partial_{\lambda\nu} d& \partial_{\nu\nu} d\end{array}\right)=\left(\begin{array}{cc} \partial_\lambda d& 0 \\ *& \partial_{\nu\nu} d\end{array}\right),
\]
with invertible diagonal entries. Adding a parameter $\mu$ in $f$, we may then continue simple pinched double roots in the parameter with the implicit function theorem, analytically in the sense of perturbation theory and global continuation, and numerically using for instance arclength continuation.  We note that it is often advantageous to solve the equation for $\lambda$ and $\nu$ without taking determinants, solving at the same time for eigenvectors and generalized eigenvectors; see~\cite{rss07}. This approach  is particularly relevant when the system size $N$ is large stemming for instance from the discretization of a system in a cylinder or a periodically forced system; see \S\ref{s:infdim}.

Near a simple pinched double root, the roots of the dispersion relation have the expansion
\begin{equation}\label{e:diff}
\lambda-\lambda_*=d_\mathrm{eff}(\nu-\nu_*)^2+\rmO((\nu-\nu_*)^3).
\end{equation}
One clearly notices that~\eqref{e:diff}  mimics a diffusion equation for a slowly modulated eigenmode,   that is, $u\sim A(\eps^2 t,\eps x)\rme^{\lambda_* t + \nu_* x}$ is a solution at order $\eps^2$ if $A(T,X)$ solves
\begin{equation}\label{e:diffmodulation}
    A_T=d_\mathrm{eff} A_{XX}, \qquad \text{with }\quad d_\mathrm{eff}=-\frac{\partial_{\nu\nu} d(\lambda_*,\nu_*)}{2\partial_\lambda d(\lambda_*,\nu_*)}.
\end{equation}
We comment on the relevance of the diffusive refined asymptotics, next, \textcolor{black}{and refer to Fig.~\ref{f:instability_pdr_dns} for an illustration in direct simulations.}

\textbf{Simple double roots and effective diffusive spreading.} Solutions to the diffusion equation \eqref{e:diffmodulation} typically decay as $|A|\sim t^{-1/2}$ for localized initial conditions, and level sets of $|u|$ at the critical speed where $\Re\lambda_*=0$  then solve $\log t^{-1/2} +\Re\nu_* x=0$, that is, in a steady frame the propagation actually occurs with speed
\begin{equation}
    c(t)=c_\mathrm{lin}-\frac{1}{2\Re\nu_* t}+\rmO(1/t^2), \qquad \qquad x(t)=c_\mathrm{lin}t-\frac{1}{2\Re\nu_*}\log t +\rmO(1).
\end{equation}
We will see in \S\ref{s:nlmsp} that such a log-shift with slow convergence to the actual speed occurs in nonlinear equations as well, albeit with a prefactor of $\frac{-3}{2\Re\nu_*}$ corresponding to the decay in the heat equation with {\color{black} ``lossy'', for instance Dirichlet or positive outward-flux Robin,} boundary conditions.

In general, the effective diffusivity can be complex when $\lambda_*\not\in\R$, that is, for oscillatory invasion.  It need in fact  not be positive, in which case the modulation~\eqref{e:diffmodulation} is ill-posed; see~\cite{AbsoluteSpecArndBjorn,FayeHolzerScheelSiemer} for examples. In the examples that we present in this section, we however always find $\Re d_\mathrm{eff}>0.$ We will encounter this effective diffusive behavior in the leading edge of invasion fronts again  when discussing the reasons for nonlinear selection in \S\ref{s:nlmsp}.

\textbf{Multiple double roots.}
Both $\partial_\lambda d=0$ or $\partial_{\nu\nu}d=0$ lead to higher multiplicities. In particular, small perturbations lead to the emergence of multiple solutions to the double-root equation. One interesting example that we will encounter again later is the case when $\partial_\lambda d=0$, and the dispersion relation is quadratic in both $\nu-\nu_*$ and $\lambda-\lambda_*$, at leading order. Scaling $\nu=\nu_1\lambda$ as suggested by the Newton polygon, we typically find two roots $\nu_\pm(\lambda)$, analytic in $\lambda$, so that, locally and  up to a nonzero analytic function,
\[
d(\lambda,\nu)=(\nu -\nu_+(\lambda))(\nu-\nu_-(\lambda)),\qquad \qquad\nu_\pm(\lambda)=\nu_{1,\pm}\cdot\lambda+\rmO(\lambda^2).
\]
In the simplest example, this leading order corresponds to wave equation dynamics rather than the diffusive dynamics encountered above.

The interest in this example is that it naturally arises in uni-directionally coupled reaction-diffusion systems of the form, say
\begin{equation}\label{e:linuni}
u_t=u_{xx}+u+\beta v,\qquad v_t=dv_{xx}+\alpha v;
\end{equation}
see Example 5 in \S\ref{s:2ex}. The analysis of double double roots presents many subtleties that we shall discuss in \S\ref{s:res}. They are best understood in a broader picture of nontrivial resonances between modes $(\lambda,\nu)$ in the dispersion relation, generalizing the double roots we are interested in here.

\subsection{Practical considerations: linear marginal stability}\label{s:practical_linearspreading}

Understanding double roots and their pinching is of interest beyond finding spreading speeds. Locating double roots amounts to solving a system of two complex equations in two variables, which even in the polynomial case considered thus far is often difficult, mostly untractable analytically as we will see in \S\ref{s:2ex}. Determining if they satisfy the pinching condition is usually cumbersome. Adding the condition $\Re\lambda=0$ complicates the algebra significantly by adding a real variable, $c$, and a real condition, $\Re\lambda=0$, so that we end up with 4 real rather than 2 complex equations.

\textbf{Finding pointwise growth modes.}
Pointwise growth modes are a subset of the set of pinched double roots, which themselves are a subset of the set of all double roots. In order to find pointwise growth modes, one therefore usually starts with the dispersion relation $d_c(\lambda,\nu)$ and finds double roots, looking for instance for roots of the discriminant. It is not difficult to see that the discriminant of $d_c$ with respect to $\nu$ is a polynomial in $\lambda$ which has in general $2mN^2-N$ roots, when counted with
multiplicity; see also~\cite[Lem. 4.5]{rss07}. This excludes degenerate cases, which we shall encounter in the next section, when there are continua of  double roots. For specific parameter values, one then usually finds all double roots as roots of the discriminant and determines which of those are pinched on a case-by-case basis, tracking roots $\nu$ from the collision to  $\lambda=+\infty$. The task for a specific parameter value is not all that cumbersome since one usually proceeds by decreasing real part, starting with the double root with maximal $\Re\lambda$ that is yet to the left of the essential spectrum.

The pinching condition for the double root with maximal real part is continuous in parameters except when it collides with another double root or when another double root takes over as root with maximal real part, so that the pinching condition only needs to be checked in isolated instances when parameters are changed~\cite[Lem.~5.8]{HolzerScheelPointwiseGrowth}.
Simple pinched double roots are always pointwise growth modes; they lead to branch points of the pointwise Green's function~\cite[Lem.~4.4]{HolzerScheelPointwiseGrowth}. Pinched double roots of higher multiplicity need to be investigated on a case-by-case basis.

Several homotopies are possible to identify the most unstable pointwise growth mode. Short of a somewhat radical homotopy to a trivial equation, say with $\mathcal{P}(\partial_x)=\mathrm{diag}\,(d_j\partial_x^{2m}+a_j)$, with all $d_j,a_j$ different, one can start homotopies at a fixed $c_0$, for instance $c_0=0$, and compute the essential spectrum, $\lambda_j(k)$, $j=1,\ldots,N$; see~\cite{rss07} for strategies. One then identifies the most unstable point in the spectrum $\max_{j,k}\Re\lambda_j(k)=\Re\lambda_{j_*}(k_*)$ and its associated group velocity, $c_\mathrm{g}=-\Im \lambda'_{j_*}(k_*)$. Passing to a coordinate frame propagating at the group velocity, that is, with speed $c_0+c_\mathrm{g}$, the pair $(\lambda_j(k_*),\rmi k_*)$ becomes the most unstable pinched double root in this particular coordinate frame.
A similar strategy continues the essential spectrum in the exponential weight $\eta=\Re\nu=\Im k$. Tracking the maximum of the real part of the spectrum to a minimal value then yields a pinched double root.

We reemphasize that in all continuation procedures, one must be cautious that the pinched double root tracked in the homotopy remains the most unstable pinched double root, which usually requires tracking all double roots. We also note that knowing information about the absolute spectrum includes in particular information about double roots; we refer to~\cite{AbsoluteSpecArndBjorn,rss07} for both background information and computational algorithms.

\textbf{Finding spreading speeds.} With an algorithm to identify pointwise growth modes, or rather pinched double roots, one now simply continues these in the parameter $c$ associated with the comoving frame until the real part vanishes. If relying on numerics, one could also identify a good initial guess for the linear spreading speed from direct simulations and use this numerical data as an initial guess in the root finding procedures.

The arguably most direct method for finding spreading speeds is based on solving~\eqref{e:spledr}, using a simple root finding algorithm. We emphasize however that solutions to this equation do not in any known way guarantee pinching conditions or extremal properties of the spreading speed obtained in this fashion.

\textbf{Nonlinear eigenvalue problems.}  We mentioned, and will emphasize throughout, that pointwise growth modes share many properties of spectral values. One might therefore wonder if there are algorithms similar to the (inverse) power method, or more generally QR-iteration, that converge to pointwise growth modes. Such algorithms were recently proposed in~\cite{S23}, exploiting extensively the first-order formulation~\eqref{e:Ufirstorder}. The algorithm proposed there identifies the closest pointwise growth mode to a fixed initial guess $\lambda_0$ in the complex plane, or alternatively the closest $k$ pointwise growth modes, alias singularities of the pointwise Green's function. We emphasize that the algorithm does not compute double roots and then verify pinching conditions, but it rather identifies singularities of subspaces directly via an iterative nonlinear eigenvalue iteration. Examples in~\cite{S23} demonstrate linear convergence as usual for inverse-power type methods. Roughly speaking, the method first computes automated power series expansions of stable and unstable subspaces $E^\mathrm{s/u}(\lambda)$ as graphs over a fixed reference subspace. It then forms a square matrix $A(\lambda)$ with power series column vectors from $E^\mathrm{s/u}(\lambda)$ and treats $A(\lambda)u_0=0$ as a nonlinear polynomial eigenvalue problem. The inverse power method to this nonlinear eigenvalue problem consists of a straightforward multi-term recursion that provably converges to the nearest spectral value. We refer to~\cite{S23} for more details but note here that this method also generalizes to other situations.
In fact, treating the wave speed $c$ as a nonlinear eigenvalue parameter and setting $\lambda=0$ allows one to directly identify the spreading speed associated with a pointwise growth mode closest to an initial guess $c_0$. Unfortunately, no such algorithm is known for spreading speeds associated with a nonzero frequency $\rmi\omega$.

More generally, the algorithm also allows one to treat eigenvalue problems that we shall encounter in \S\ref{s:nlms}, with $x$-dependent coefficients, and thus investigate nonlinear marginal stability.

Compared to direct simulations, the approach here
\begin{itemize}
    \item introduces correct, or ``asymptotic'' boundary conditions;
    \item converges linearly, that is the error decays exponentially in the number of iterates, rather than algebraically when relying on direct simulations.
\end{itemize}

\subsection{Linear spreading speeds: examples}
\label{s:2ex}
We illustrate the concepts and tools presented in this section with several examples. \textcolor{black}{All examples have a reflection symmetry $x\mapsto -x$ so that it will be sufficient to restrict to positive spreading speeds $c_\mathrm{lin}$ with double root $(\lambda_\mathrm{lin},\nu_\mathrm{lin})$, understanding that $-c_\mathrm{lin}$ is also a spreading speed with double root $(\lambda_\mathrm{lin},-\nu_\mathrm{lin})$.}

\begin{figure}
    \centering\includegraphics[height=1.5in]{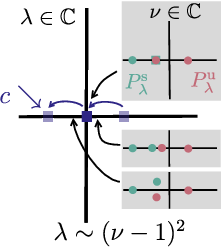}\qquad\qquad
    \includegraphics[height=1.5in]{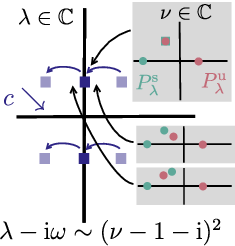}\qquad\qquad
    \includegraphics[height=1.5in]{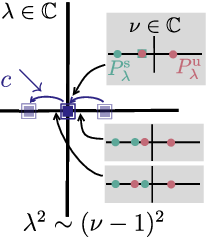}
    \caption{\textcolor{black}{Marginal linear stability induced by pinched double roots at $c=c_\mathrm{lin}$ as $c$ decreases in case of stationary leading edge (left), oscillatory leading edge (center), and a non-branched resonance (right). The large panels show movement of pinched double roots as $c$ decreases. The insets show the location of spatial roots $\nu$ for locations at (top insets) and near (bottom insets) the double root. For the branched double roots (left,center), the spatial roots rotate by $\pi/2$ as $\lambda$ passes through the double root; for the unbranched double root, the roots pass smoothly through each other. See examples, particularly Example 5, for a discussion of unbranched double roots.} }
    \label{f:instability_pdr}
\end{figure}

\textcolor{black}{\textbf{Overview.} We illustrate spreading mechanisms in 7 examples. Examples show that spreading speeds can be calculated explicitly in simple systems (Examples 1, 2, 3, 5, 6), generally not explicitly in more complex systems (Example 7), and perturbatively (Example 4). In addition to identifying the spreading speed, we find that the mode of invasion can be oscillatory or stationary, with simple pinched double roots crossing at $\lambda=0$ or $\lambda=\im\rmi\omega_*\neq 0$ as $c$ decreases through $c_\mathrm{lin}$; see Fig.~\ref{f:instability_pdr}. The transition occurs when double roots merge, leading to multiple double roots, which also give rise to more subtle spreading behavior illustrated in Example 5 and in the right panel of Fig.~\ref{f:instability_pdr}. We start with scalar equations with second-order and fourth-order diffusion in Examples 1 and 2, respectively. Example 1 shows the correspondence between pinched double roots, pointwise growth modes, singularities of the pointwise resolvent, and growth rates of the heat kernel explicitly. These explicit calculations become rather cumbersome in Example 2, where we resort to pinched double root calculations which predict either oscillatory or stationary spreading. Example 3 is the simplest example of a system, complex linear, with oscillatory spreading, and Examples 4 and 5 are perturbations of this simple example. The final examples demonstrate how in systems of 2 reaction-diffusion equations it is sometimes possible to find closed-form expressions for spreading speeds, Example 6, but generally not, Example 7. 
}

\begin{figure}[ht!]
    \begin{minipage}{0.2\textwidth}\centering $u$\\[0.1in]
    \includegraphics[width=0.9\textwidth]{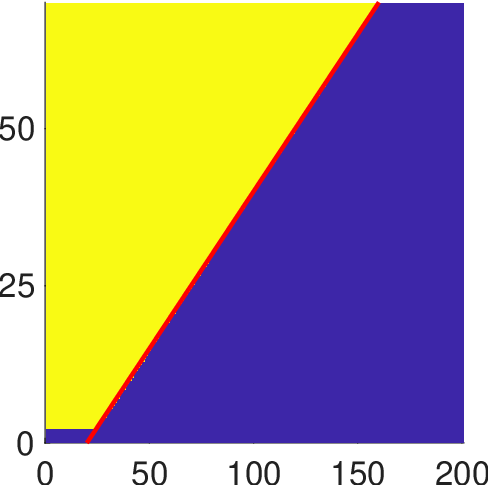}\\[0.1in]
    \includegraphics[width=0.9\textwidth]{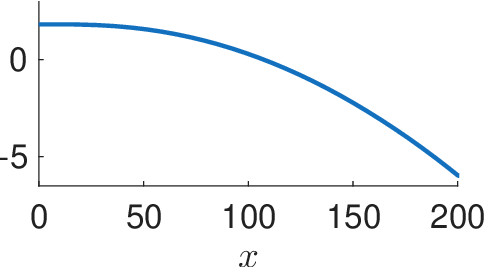}\\ [0.1in]   \includegraphics[width=0.9\textwidth]{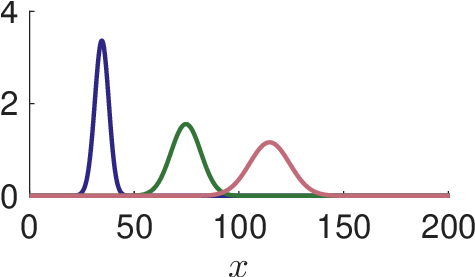}
    \end{minipage}\hfill
    \begin{minipage}{0.2\textwidth}\centering $\Re A$\\[0.1in]
    \includegraphics[width=0.9\textwidth]{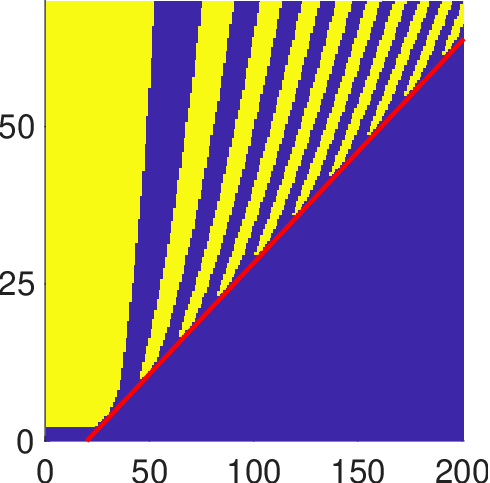}\\[0.1in]
    \includegraphics[width=0.9\textwidth]{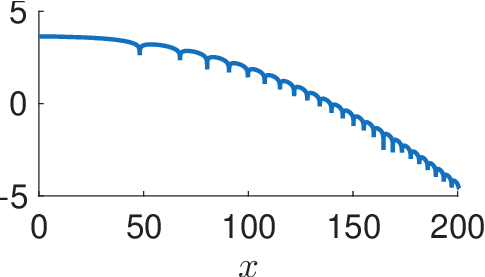}\\[0.1in]  \includegraphics[width=0.9\textwidth]{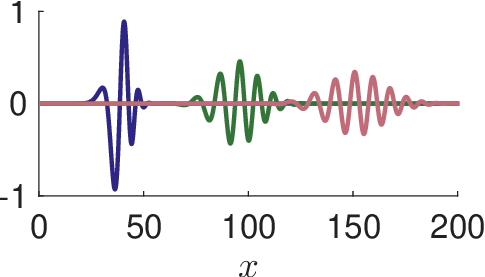}
    \end{minipage}\hfill
    \begin{minipage}{0.2\textwidth}\centering $u$\\[0.1in]
    \includegraphics[width=0.9\textwidth]{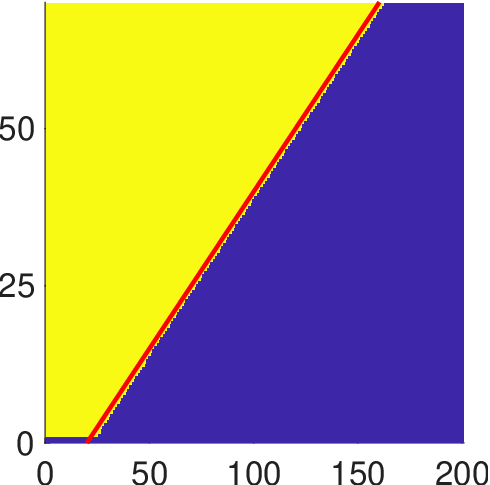}\\[0.1in]
    \includegraphics[width=0.9\textwidth]{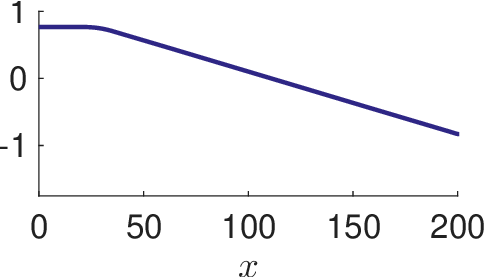}\\[0.1in]  \includegraphics[width=0.9\textwidth]{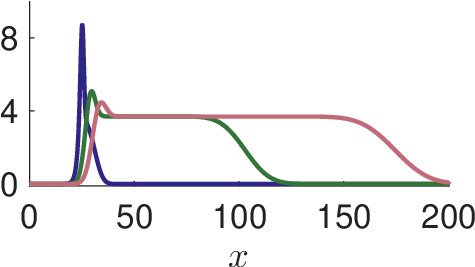}
    \end{minipage}\hfill
    \begin{minipage}{0.2\textwidth}\centering $v$\\[0.1in]
    \includegraphics[width=0.9\textwidth]{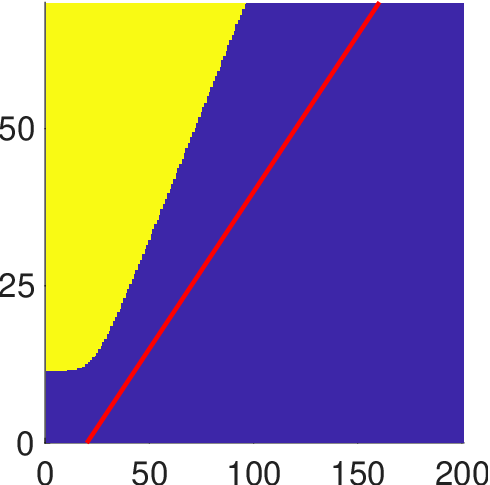}\\[0.1in]
    \includegraphics[width=0.9\textwidth]{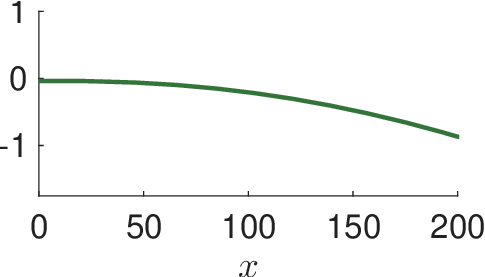}\\[0.1in]  \includegraphics[width=0.9\textwidth]{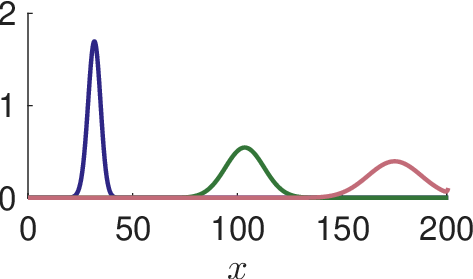}
    \end{minipage}\qquad\\
\begin{minipage}{0.2\textwidth}
    \[u_t=u_{xx}+u\]
    \end{minipage}\hfill
    \begin{minipage}{0.2\textwidth}
    \[A_t=(1+\rmi\alpha) A_{xx}+A\]
    \end{minipage}\hfill
    \begin{minipage}{0.4\textwidth}
    \[\begin{aligned}u_t=&(1+\delta)u_{xx} +(1+\gamma)u + v,\\
    v_t=&(1-\delta)v_{xx}+(1-\gamma)v
    \end{aligned}\]
    \end{minipage}\hfill
    \caption{\textcolor{black}{Direct simulations of equations exhibiting spreading according to the three mechanisms depicted in Fig.~\ref{f:instability_pdr} shown through space-time plots (top row), plots of the log of the amplitude (second row), and plots of the solution when multiplied with the exponential weight $\rme^{-\Re\nu_\mathrm{lin}\xi}$ (bottom row). Asymptotics are diffusive (left), oscillating in addition to diffusive-dispersive (center-left), and wave-equation like transport (center-right and right; $\delta=-0.9,\gamma=0.8$); compare~\eqref{e:diffmodulation} for diffusive and diffusive-dispersive asymptotics and Example 5, particularly Rem.~\ref{r:wave}, for wave-equation asymptotics. }  }
    \label{f:instability_pdr_dns}
\end{figure}

\textbf{Example 1: linear spreading speeds in a scalar reaction-diffusion equation.}
We first compute the linear spreading speeds and pointwise growth modes in the simplest example consisting of a scalar reaction-diffusion equation posed in a co-moving frame \textcolor{black}{with speed $c\geq 0$},
\begin{equation} u_t=u_{\xi\xi}+cu_\xi+u. \label{e:FKPPlin} \end{equation}
This equation is obtained for instance by linearizing  the FKPP equation $u_t=u_{\xi\xi}+u-u^2$ \textcolor{black}{near the equilibrium solution $u\equiv 0$} in a co-moving frame of speed $c$.   The invasion speed for this scaling of FKPP is well known to be $2$. \textcolor{black}{This can be seen very explicitly by computing the associated heat kernel $H(t,\xi)=\frac{1}{\sqrt{4\pi t}}\rme^{-\frac{(\xi-ct)^2}{4t}+t}$ and finding  exponential decay for fixed $\xi$ as $t\to\infty$ precisely when $|c|<2$ by completing the square in the exponent. We show how this is, somewhat equivalently, recovered from pinched double root computations or pointwise resolvents. The results of the computations here are summarized graphically in  Fig.~\ref{f:instability_pdr_dns} (left column), which shows the diffusive leading-edge dynamics,  and the dispersion relation is described in the left column of Fig.~\ref{f:instability_pdr}. }

From  the dispersion relation 
\[ d(\lambda,\nu)=\nu^2+c\nu+1-\lambda,\]
we find the equation for  double roots 
\[ 0=\left(\begin{array}{c} d(\lambda,\nu) \\
\partial_\nu d(\lambda,\nu)\end{array}\right)=\left(\begin{array}{c} \nu^2+c\nu+1-\lambda \\
 2\nu+c\end{array}\right),\]
from which we obtain the existence of a pointwise growth mode at
\[ \lambda_*=1-\frac{c^2}{4}.\]
Requiring that $\Re\lambda_*=0$ leads to the linear spreading speed $c_{\mathrm{lin}}=2$, \textcolor{black}{recalling that we restrict to $c_\mathrm{lin}>0$ by reflection symmetry.}

As described previously, these pointwise growth modes arise from singularities of the pointwise Green's function.  It is instructive to work through the derivation of $G_\lambda$ in this example.  The Green's function satisfies the following differential equation:
\begin{equation} \delta(\xi-y) = G_\lambda ''+cG_\lambda ' +(1-\lambda)G_\lambda. \label{e:GKPPdef} \end{equation}
Construction of $G_\lambda$ requires the identification of exponential solutions which decay as $\xi\to \pm \infty$.  These decay rates correspond to roots of the dispersion relation,
\[ \nu_\pm(\lambda) = -\frac{c}{2}\pm \frac{1}{2}\sqrt{c^2-4+4\lambda}. \]
Note that $\mathrm{Re}(\nu_-(\lambda))<0$ for all $\lambda$ while $\nu_+(\lambda)$ crosses the imaginary axis when $\lambda=1-k^2+c\rmi k$ lies in the spectral set $\mathrm{spec}_{L^2}\,(\mathcal{L})$. \textcolor{black}{The solution to \eqref{e:GKPPdef} can be found for instance rewriting it as a vector equation as in \eqref{e:vecgreens} and computing projections explicitly using only spatial roots $\nu_\pm(\lambda)$}, finding
\begin{equation}
G_\lambda(\xi-y)= \left\{ \begin{array}{cc} \frac{1}{\nu_-(\lambda)-\nu_+(\lambda)} \rme^{\nu_-(\lambda)(\xi-y)}, & \xi\geq y, \\ \frac{1}{\nu_-(\lambda)-\nu_+(\lambda)} \rme^{\nu_+(\lambda)(\xi-y)}, & \xi<y \end{array}\right. . \end{equation}
Note that $G_\lambda\notin L^1(\mathbb{R})$ for $\lambda$ with $\mathrm{Re}(\nu_+(\lambda))\leq 0$ but that for fixed $\xi$, $G_\lambda(\xi)$ can be continued analytically in $\lambda$ into this set regardless.  Substituting the explicit formulas for the roots $\nu_\pm(\lambda)$, the Green's function takes the form
\begin{equation}
G_\lambda(\xi-y)= \left\{ \begin{array}{cc} \frac{-1}{\sqrt{c^2-4+4\lambda}} \rme^{-\frac{c}{2}(\xi-y)} \rme^{-\frac{1}{2}\sqrt{c^2-4+4\lambda}(\xi-y)}, & \xi\geq y \\ \frac{-1}{\sqrt{c^2-4+4\lambda}}\rme^{-\frac{c}{2}(\xi-y)} \rme^{\frac{1}{2}\sqrt{c^2-4+4\lambda}(\xi-y)}, & \xi<y \end{array}\right. . \label{e:GKPPexplicit} \end{equation}
Note that singularities of the Green's function, i.e. pointwise growth modes,  occur at $\lambda=1-c^2/4$, where the square root vanishes --- precisely as computed from double root criterion.

We conclude this example with some comments on the use of exponential weights.  Starting again from the original PDE (\ref{e:FKPPlin}) we could take $w(\xi,t)=\rme^{-\eta \xi}u(\xi,t)$ after which the PDE is transformed to
\[ w_t= w_{\xi\xi}+(c+2\eta)w_\xi+\left(\eta^2+c\eta+1\right)w. \]
Observe that the coefficient $\eta^2+c\eta+1$ contributes the pointwise growth or decay of the solution in this weighted space while the advective coefficient $(c+2\eta)$ affects the pointwise decay rate via transport.  Setting both coefficients to zero, and thus obtaining the simple heat equation predicted more generally in~\eqref{e:diffmodulation}, amounts to solving the following system of equations
\[  0=\left(\begin{array}{c} \eta^2+c\eta+1 \\
 2\eta+c\end{array}\right).\]
This is precisely the double root condition and recovers the linear spreading speed $c=2$.  In a frame moving with this linear speed and in a weighted space with $\eta=-1$ we find that the dynamics reduce to the heat equation thereby recovering the diffusive dynamics predicted more generally from~\eqref{e:diff}.

\begin{figure}
    \centering\includegraphics[width=0.65\textwidth]{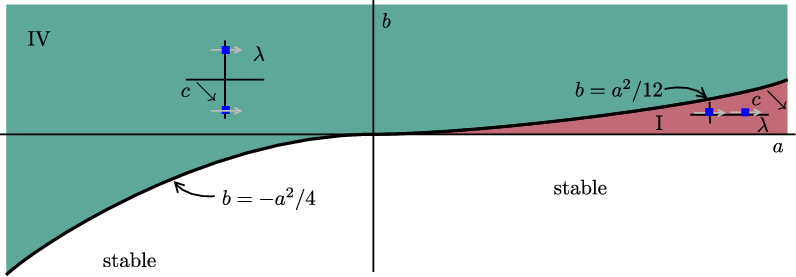}
    \caption{Parameter regions (green and red) where the origin in~\eqref{e:4thorder} is unstable, divided into a region (IV) (green) where the spreading speed is induced by a complex conjugated pair of pinched double roots, with ensuing oscillations in the leading edge, and (I) (red), where the leading edge behavior is stationary. \textcolor{black}{Gray arrows in insets indicate movement of double roots in the $\lambda$-plane as $c$ decreases.}  In case (I), two real double roots cross the origin as $c$ decreases but only the more stable one is pinched; see~\eqref{e:regions IIIIV} for a description of regions.}\label{f:spsp4}
\end{figure}

\textbf{Example 2: Spreading speeds in fourth order scalar equations, multiple and complex double roots.}
More interesting yet still explicit examples are scalar fourth-order equations, such as the Swift-Hohenberg~\eqref{e:sh} and the Cahn-Hilliard equation~\eqref{e:ch}. Spreading speeds were listed in~\cite{vanSaarloosReview} and we give some details here. \textcolor{black}{While we constructed explicit heat kernels and alternatively Green's functions  in the previous example, both of which point to pointwise growth and spreading speeds, both become rather unwieldy in the example here, so that pinched double-root computations give the preferred computational method toward determining spreading behavior. Even the determination of pinched double roots is, however, not immediate.} There are in fact multiple double roots and it requires care to identify which ones are pinched.  The results are summarized in Fig.~\ref{f:spsp4} and formulas for roots $\nu$  in~\eqref{e:4th2},~\eqref{e:4th4}, and speeds and frequencies in~\eqref{e:4cth2},~\eqref{e:4cth4}, for the regions (II) and (IV), respectively. \textcolor{black}{We find both stationary leading edge dynamics as shown in the left columns of Figs.~\ref{f:instability_pdr} and~\ref{f:instability_pdr_dns}, but also oscillatory dynamics as shown in the corresponding center-left columns. }

Specifically, consider
\begin{equation}\label{e:4thorder}
    u_t=-u_{xxxx} + a u_{xx}+ b u.
\end{equation}
The essential spectrum $\{\lambda=-k^4-ak^2+b,\, k\in\R\}$ possesses a maximum at $k=0$, $\lambda=b$, when $a\geq 0$, and at $\lambda=\frac{a^2}{4}+b$, $k=\sqrt{-a/2}$, when $a<0$. We therefore restrict to the regions
\begin{equation}\label{e:unstableregion}
A=\{ a\geq0,\ b>0\}; \qquad \qquad \text{and }\quad
B=\{ a<0,\ \frac{a^2}{4}+b>0\}.
\end{equation}
In~\eqref{e:spledr}, we then use $\mathcal{P}(\nu)=b+a\nu^2-\nu^4$ with $\nu=\nu_\mathrm{r}+\rmi\nu_\mathrm{i}$, to find
\begin{equation}
    \begin{aligned}
        -2 a \nu_\mathrm{r}-12\nu_\mathrm{i}^2\nu_\mathrm{r} + 4\nu_\mathrm{r}^3&=\frac{-b+a\nu_\mathrm{i}^2+\nu_\mathrm{i}^4}{\nu_\mathrm{r}}-a\nu_\mathrm{r}-6\nu_\mathrm{i}^2\nu_\mathrm{r}+\nu_\mathrm{r}^3,\\
          -2 a \nu_\mathrm{i}+12\nu_\mathrm{r}^2\nu_\mathrm{i} -4 \nu_\mathrm{i}^3&=0.
    \end{aligned}
\end{equation}
The value of $\nu_\mathrm{i}$  from the second equation can be substituted into the first equation to yield solutions (listing only solutions with $\nu_\mathrm{r}<0$ corresponding to $c>0$~\cite[Rem.6.6]{HolzerScheelPointwiseGrowth}, and with  $\nu_\mathrm{i}\geq 0$, thus
omitting complex conjugates $\pm\nu_\mathrm{i}$),
\begin{align}
       \text{(I):}&& \nu_\mathrm{r,I}&=-\frac{\sqrt{a-\sqrt{a^2-12 b}}}{\sqrt{6}},&\nu_\mathrm{i,I}&=0,\label{e:4th1}\\
          \text{(II):}&& \nu_\mathrm{r,II}&=-\frac{\sqrt{a+\sqrt{a^2-12 b}}}{\sqrt{6}},&\nu_\mathrm{i,II}&=0,\label{e:4th2}\\
         \text{(III):}&&  \nu_\mathrm{r,III}&=-\frac{\sqrt{a-\sqrt{7a^2+24 b}}}{2\sqrt{6}},&\nu_\mathrm{i,III}&=\frac{\sqrt{-3a-\sqrt{7a^2+24 b}}}{2\sqrt{2}},\label{e:4th3}\\
         \text{(IV):}&&  \nu_\mathrm{r,IV}&=-\frac{\sqrt{a+\sqrt{7a^2+24 b}}}{2\sqrt{6}},&\nu_\mathrm{i,IV}&=\frac{\sqrt{-3a+\sqrt{7a^2+24 b}}}{2\sqrt{2}}.\label{e:4th4}
\end{align}
We find roots in (III) are not real in region $A\cup B$ from~\eqref{e:unstableregion}, and roots are real in (I), (II), (IV) when
\begin{equation}\label{e:regions IIIIV}
    \text{(I): } \{a>0,0<b<\frac{a^2}{12}\};\  \text{(II): } \{a>0\}\cup \{a<0,b<0\};  \   \text{(IV): } \{a>0,b>\frac{a^2}{12}\}\cup \{a<0,b>-\frac{a^2}{4}\};
\end{equation}
Clearly, complex values of $\Re\nu$ and $\Im\nu$ are irrelevant.
The speeds and frequencies associated with double roots are
\begin{align}
       \text{(I):}&& c_\mathrm{lin,I}&=\frac{2}{3\sqrt{6}} (2a+\sqrt{a^2-12 b})\sqrt{a-\sqrt{a^2-12 b}},&&\nonumber\\
       &&\omega_\mathrm{I}&=0 ,\label{e:4cth1}\\
       \text{(II):}&& c_\mathrm{lin,II}&=\frac{2}{3\sqrt{6}}(2a-\sqrt{a^2-12 b})\sqrt{a-+\sqrt{a^2-12 b}},&&\nonumber\\
       &&\omega_\mathrm{II}&=0 ,\label{e:4cth2}\\
       \text{(IV):}&& c_\mathrm{lin,IV}&=\frac{2}{3\sqrt{6}} (-2a+\sqrt{7a^2+24 b})\sqrt{a+\sqrt{7a^2+24 b}},&&\nonumber\\
       &&\omega_\mathrm{IV}&=\frac{1}{8\sqrt{3}} (-3a+\sqrt{7a^2+24 b})^{3/2}\sqrt{a+\sqrt{7a^2+24 b}} ,\label{e:4cth4}
\end{align}
In case (II), $c_\mathrm{lin,II}<0$ when $a<0$ so that only the region $a>0$, $0<b<a^2/12$ is admissible. In this region, one easily checks that $c_\mathrm{lin,II}>c_\mathrm{lin,I}$.  We claim that the double root associated with $c_\mathrm{lin,II}$ is not pinched. For this, we consider the essential spectrum in a weighted space, setting $\nu=-\sqrt{a/6} + \rmi k$, with the weight being an average of the weights from $\nu_\mathrm{r,I}$ and $\nu_\mathrm{r,II}$, to find, in a frame with speed $c_\mathrm{lin,II}$,
\[
\Re\lambda=\frac{1}{36}\left(a^2- 12 b -a\sqrt{a^2-12b}\right)-k^4=\frac{a^2}{36}
\left((1-\tau)-4\sqrt{1-\tau} \right)-k^4, \quad \text{ with } b=\tau a^2/12,\tau\in(0,1).
\]
Since the right-hand side is strictly negative, we conclude that at speed $c_\mathrm{lin,II}$, the system is exponentially stable in this weighted space so that the double root at the origin cannot be pinched.

Altogether, we found spreading speeds $c_\mathrm{I/IV}$ with double roots $(\lambda,\nu)=(\rmi\omega_\mathrm{I/IV},\nu_\mathrm{r,I/IV})$ in
\[
\text{(I): } \{a>0,0<b<a^2/12\},\qquad\qquad \text{(IV): } \{a>0,b>a^2/12\}\cup \{a<0,b>-a^2/4\},
\]
respectively.
When crossing  the curve $b=a^2/12, a>0$ with increasing $b$, two real double roots meet at the origin and become complex conjugates. The spreading speeds are obtained by setting the real parts of those double roots equal to zero.
It remains to check that the double roots are in fact pointwise growth modes. In both cases, there necessarily is a crossing of a pinched double root for some speed $c$, and, having eliminated all other options, the speeds $c_\mathrm{lin,I/IV}$ are necessarily the correct ones.

Continuing in the $(a,b)$-plane from $b>a^2/12$ into $b<a^2/12$ while $a>0$, the two complex pinched double roots merge into two real double roots, where only the more stable one is pinched. We shall later see an example~\eqref{e:cgl22real} where both double roots are pinched.

Finally, we also list effective diffusivities~\eqref{e:diffmodulation}, obtained by Taylor expanding at the double root:
\begin{align}
      \text{(I):}&& d_\mathrm{eff,I}&=\sqrt{a^2-12b};\\
      \text{(IV):}&& d_\mathrm{eff,IV}&=\frac{3a-\sqrt{7a^2+24b}}{2} - \frac{\rmi}{2} \sqrt{3} \sqrt{\sqrt{7 a^2+24 b}-3 a}
   \sqrt{\sqrt{7 a^2+24 b}+a}+\frac{3 a}{2}.
\end{align}

\textbf{Example 3: Complex Ginzburg-Landau equation, complex double roots.}
The onset of oscillations in spatially extended systems is described at leading order by the universal complex Ginzburg-Landau modulation equation
\begin{equation}\label{e:cgln}
A_t=(1+\rmi\alpha)A_{xx}+(1+\rmi\omega)A-(1+\rmi\beta) A|A|^2,
\end{equation}
where $\omega$ measures the detuning of the frequency with the parameter, $\alpha$ the detuning with spatial modulation (linear dispersion), and $\beta$ the detuning with amplitude (nonlinear dispersion). With the gauge invariance, one could eliminate $\omega$ but we shall keep this parameter.
The linearization at the (unstable) origin is
\begin{equation}\label{e:ccgllin}
A_t=(1+\rmi\alpha) A_{xx} + (1+\rmi\omega) A, \qquad A\in\C.
\end{equation}
\textcolor{black}{The spreading dynamics of this equation are illustrated in the center-left columns of Figs.~\ref{f:instability_pdr} and~\ref{f:instability_pdr_dns} for $\alpha\neq 0,\, \omega=0$.}

Quickly setting $A=\rme^{\lambda t+\nu (x- ct)}$, we find $d_c(\lambda,\nu)=\lambda-c\nu-(1+\rmi\alpha)\nu^2-(1+\rmi\omega)$. Then solving  $\partial_\nu d_c=0$ gives $\nu=-c/(2(1+\rmi\alpha))$, which, substituting into $d=0$ and setting $\lambda=\rmi\Omega$ gives
\begin{equation}\label{e:cglclin}
c=2\sqrt{1+\alpha^2},\qquad \Omega=\alpha+\omega, \qquad \nu=\frac{-1+\rmi\alpha}{\sqrt{1+\alpha^2}}.
\end{equation}
Clearly, in this scalar equation, $d_\mathrm{eff}=(1+\rmi\alpha)$.
Alternatively, one writes $A=u+\rmi v$, finds
\[
u_t=u_{xx}-\alpha v_{xx} + u-\omega v,\qquad v_t=v_{xx}+\alpha u_{xx} + v+\omega u,
\]
with dispersion relation at $c=0$ \textcolor{black}{computed as in \eqref{eq:disp} as $d_0(\lambda,\nu)=(\nu^2(1+\rmi\alpha)+1+\rmi\omega-\lambda)(\nu^2(1-\rmi\alpha)+1-\rmi\omega-\lambda)$, and double roots of $d_c(\lambda, \nu)=d_0(\lambda-c\nu,\nu)$ }at $\lambda=\pm \rmi(\alpha+\omega)$,  $\nu=\frac{-1\pm\rmi\alpha}{\sqrt{1+\alpha^2}}$ for $c = 2 \sqrt{1+\alpha^2}$.
We will return to this example later when studying pattern forming fronts.

\textbf{Example 4: Parametrically forced Ginzburg-Landau equation, splitting of multiple double roots.}
We now turn to a variation  of~\eqref{e:cgln} with parametric forcing: the original oscillations with approximate frequency $\omega_\mathrm{Hopf}$ are subjected to a weak, spatially homogeneous  external forcing with frequency  $\omega_\mathrm{forcing}=2\omega_\mathrm{Hopf}$. In the modulation approximation, this forcing breaks the gauge symmetry and introduces a term $\gamma\bar{A}$
\begin{equation}\label{e:cgl22}
 A_t=(1+\rmi\alpha)A_{xx}+(1+\rmi\omega)A-(1+\rmi\beta) A|A|^2+\gamma \bar{A};
\end{equation}
see \cite{doelmanschielen,cgl_force_valid,kamphuis2025patternformationswifthohenbergequation} for background,  derivations, and validity. Higher resonances in the forcing $\omega_\mathrm{forcing}=\ell\omega_\mathrm{Hopf}$, $\ell=3,4,5,\ldots$ lead to terms $\bar{A}^{\ell-1}$ and do not change the linear part at leading order.
\textcolor{black}{The goal in this example is to both illustrate the basic techniques and potential pitfalls of a perturbative analysis, as well as to document interesting physical effects corresponding to frequency locking and unlocking in the leading edge, marked by the change from stationary to oscillatory leading edge behavior, switching between left and center-left leading-edge dynamics in Figs.~\ref{f:instability_pdr} and~\ref{f:instability_pdr_dns}.}

Possibly rotating $A\mapsto \rme^{\rmi\varphi}A$, we can assume that $\gamma>0$.
The linearization
\begin{equation}\label{e:cgl2lin}
A_t=(1+\rmi\alpha)A_{xx}+(1+\rmi\omega)A+\gamma \bar{A},
\end{equation}
of~\eqref{e:cgl22}  leads to a  dispersion relation \textcolor{black}{in the steady frame}
\[
d_0(\lambda,\nu)=\left(\alpha ^2+1\right) \nu
   ^4+(2 \alpha  \omega
   -2 \lambda +2)\nu ^2 -\gamma
   ^2+\lambda ^2-2 \lambda
   +\omega ^2+1.
\]
Finding pinched double roots, one can solve $\partial_\nu d_c=0$ for $\lambda$ as a function of $\nu$ and substitute into $d_c=0$, to find a 6th order polynomial in $\nu$ with no obvious way to  analyze the solutions analytically.
We therefore focus on a  perturbative analysis near $\omega=\alpha=\gamma=0$ where the dispersion relation is
\[
d_c(\lambda,\nu)=(\lambda-c\nu-1-\nu^2)^2,
\]
simply an uncoupled product of two copies of the linearized Fisher equation. We readily conclude that the spreading speed is $c=2$ associated with a double root $\lambda=0,\nu=-1$. Near the double root, setting $\nu=-1+\hat{\nu}$, we have
\[
d_c({\lambda}-2(-1+\hat{\nu}),-1+\hat{\nu})=(\lambda-\hat{\nu}^2)^2, \qquad c=2.
\]
The double root in this equation is degenerate as there is a family of double root solutions for all $\lambda$, albeit not pinched. It is therefore not possible to assign a multiplicity to this double root. In fact, depending on perturbations, there may be two or four double roots in a vicinity of $\lambda=0,\nu=-1$.  For generic $\alpha,\gamma,\omega$, however, the system has 6 double roots; see~\cite[\S4]{rss07} for a homotopy argument.

We next set $\alpha=\eps\alpha_1$, $\gamma=\eps\gamma_1$, $\omega=\eps\omega_1$.
In an exponential weight $\eta=1$ and a frame with speed $2+\hat{c}$, we find that the essential spectrum is given through
\[
\sigma_\mathrm{ess}=\{\rmi  \hat{c} k-\hat{c}-k^2-\sqrt{\gamma ^2-\alpha
   ^2-2 \alpha  \omega -\alpha ^2
   k^4-4 \rmi \alpha ^2 k^3+6 \alpha ^2 k^2+2
   \alpha  k^2 \omega +4 \rmi \alpha ^2 k+4 \rmi
   \alpha  k \omega -\omega ^2}|k\in\R\}.
\]
In particular, $\mathrm{max}\,\Re\lambda=\rmO(\eps)$ when $\alpha,\omega,\gamma,\hat{c}=\rmO(\eps)$. We may therefore scale
$
\alpha=\alpha_1\eps,\ \gamma=\gamma_1\eps,\ \omega=\omega_1\eps,\ \lambda=\lambda_1\eps,\nu=\nu_1\eps,
$
and track only double roots with $\lambda_1,\nu_1=\rmO(1)$ as $\eps\to 0$. We then set $\lambda_1=\rmi \Omega_1$ and find at $\rmO(\eps^2)$,
\begin{equation}\label{e:expo2}
    \begin{aligned}\alpha_1 ^2+2 \alpha_1  \omega_1
       +c_1^2-\gamma_1 ^2+\omega_1 ^2-\Omega_1 ^2,-2
       \left(c_1^2+2 c_1 \nu_\mathrm{1,r}+2 \left(\alpha_1
       ^2+\alpha  \omega_1 -\nu_\mathrm{1,i} \Omega_1
       \right)\right)&=0,\\
       2 c \Omega_1 ,-2 (c (2 \nu_\mathrm{i1,}+\Omega_1
       )+2 \nu_\mathrm{1,r} \Omega_1 )&=0.
    \end{aligned}
\end{equation}
This polynomial system can be solved  in $(\Omega_1,c_1,\nu_\mathrm{1,r},\nu_\mathrm{1,i})$ with solutions, setting
$
\mathcal{D}=(\alpha_1+\omega_1)^2-\gamma_1^2,
$
\begin{itemize}
    \item Case 1: $\mathcal{D}>0$:
\begin{equation}\label{e:cgl22imag}
    c_1=0,\quad \Omega_1=\mathcal{D},\quad \nu_\mathrm{1,r}=0,\quad \nu_\mathrm{1,i}=\frac{\alpha_1(\alpha_1+\omega_1)}{\sqrt{\mathcal{D}}};
\end{equation}
    \item Case 2: $\mathcal{D}<0$:
\begin{equation}\label{e:cgl22real}
    c_1=\sqrt{-\mathcal{D}},\quad \Omega_1=0,\quad \nu_\mathrm{1,r}=\frac{\alpha_1^2+\gamma_1^2-\omega^2}{2\sqrt{\mathcal{-D}}},\quad \nu_\mathrm{1,i}=0;
\end{equation}
\end{itemize}
In fact, there is also the complex conjugate solution with $\Omega_1=-\mathcal{D}$ in Case 1, and a solution with $c_1=-\sqrt{-\mathcal{D}}$ in Case 2. The former is simply the complex conjugate mode, the latter is irrelevant since the speed is less than the other spreading speed.


Computing the essential spectrum in a weighted space with weight $\eta=1-\eps\nu_{1,\mathrm{r}}$, one finds at leading order that the spectrum is pinched at $\lambda=\eps\rmi\Omega_1$, so that the spreading speeds correspond to pinched double roots and give accurate predictions.

The transition between Cases 1 and 2 is equivalent to the merging of two real eigenvalues that then split into a complex pair as the parameter is varied. The parameter here is $\mathcal{D}$, and the eigenvalues are rather double roots of the dispersion relation, while the parameter $c$ is adjusted so that the rightmost eigenvalue is located on the imaginary axis. A similar transition has been observed in multiple other instances; see for instance~\cite{GMS1} for an example in a phase-field system and Example 2.

The sign of $\mathcal{D}$ compares the size of  the  linear frequency $\alpha+\omega$ to  the size of the forcing $\gamma$. Heuristically, large linear frequency favors oscillations, while large $\gamma$ favors locking, zero frequency. Case 2 then reflects locking of the leading edge into the periodic forcing, while  Case 1 implies oscillatory detuning by $\alpha$ and $\omega$ and suggests the creation of non-$(1:2)$-resonant oscillations.

\textbf{Example 5: Mode competition, resonant and anomalous spreading.}
In the previous example, $\gamma$ effectively favored the real part of the amplitude, breaking the gauge invariance and attempting to lock to this real, in-phase oscillation. More complicated competition arises when including higher-order terms and simultaneous effects of other resonances. \textcolor{black}{This example will illustrate the robust presence of a double double root that results in wave equation type dynamics in the leading edge illustrated in the center-right and right panel of Figs.~\ref{f:instability_pdr} and~\ref{f:instability_pdr_dns}.}

\begin{figure}
    \raisebox{0.05in}{\includegraphics[width=0.24\textwidth]{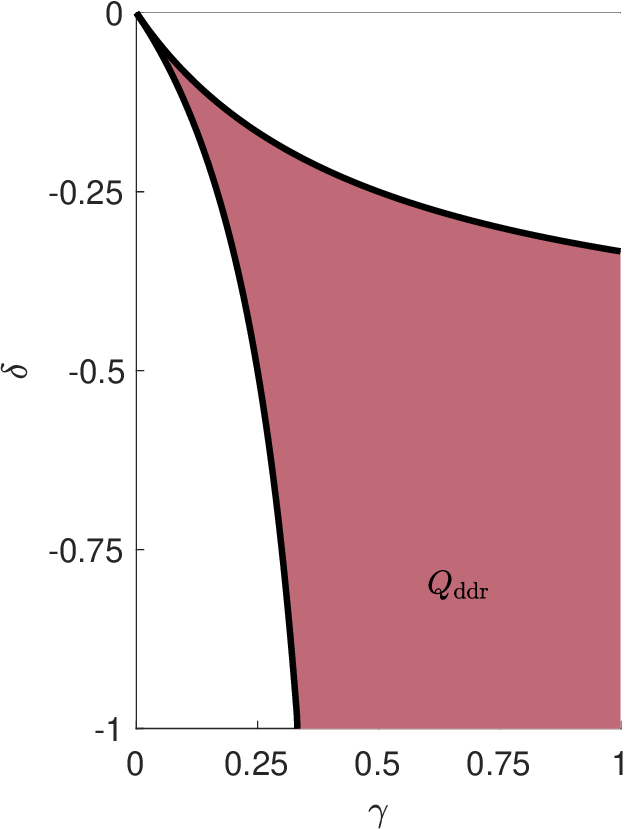}}\hfill%
    \includegraphics[width=0.2\textwidth]{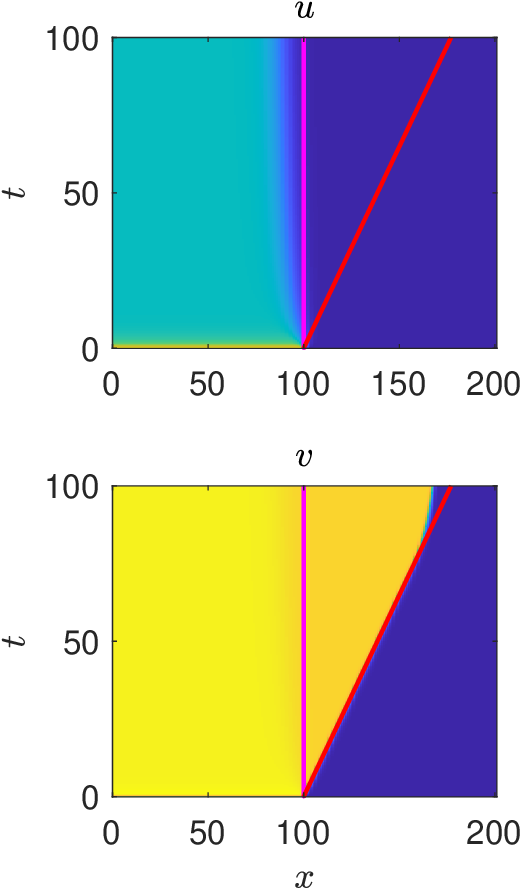}\hfill%
    \includegraphics[width=0.2\textwidth]{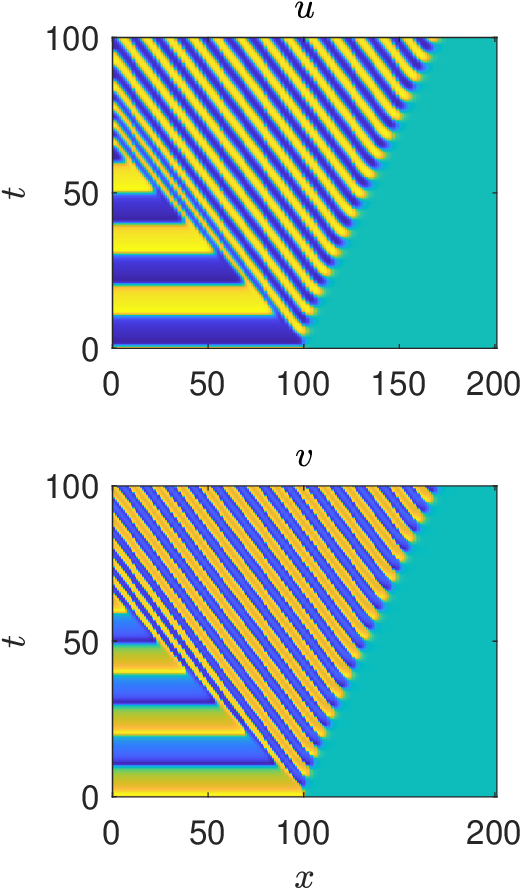}\hfill%
    \includegraphics[width=0.2\textwidth]{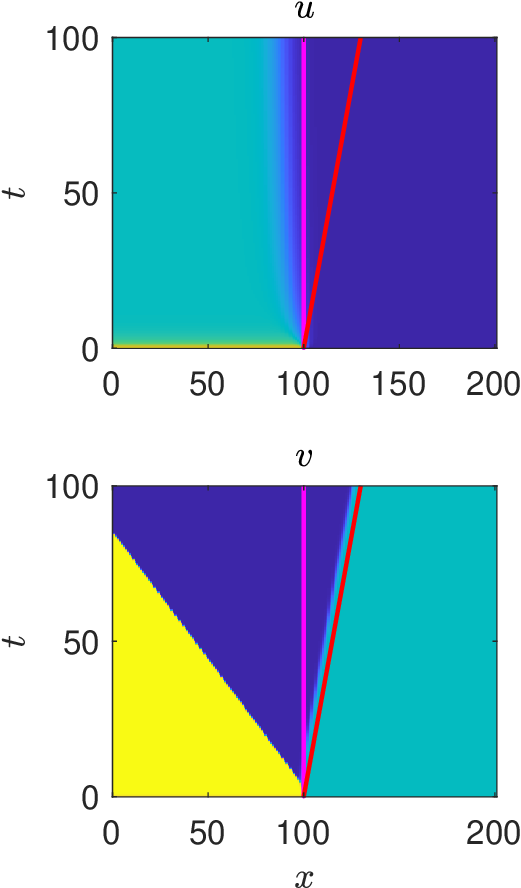}%
    \caption{Region $Q_\mathrm{ddr}$ in parameter space where the most unstable double root of the dispersion relation is a double double root coupling $u$ and $v$ equations in~\eqref{e:cgl12lin}. Space-time plots of numerical solutions to~\eqref{e:cmc} with $\gamma=0.8,\,\delta=-0.9$, in the region $Q$,  in a frame moving with the linear spreading speed $c=2\sqrt{1-\delta}\sqrt{1-\gamma}$ of~\eqref{e:cgl12lin}. The left panel with $\eps_1=1,\eps_2=\iota=0$ illustrates ``relevant'' coupling; the center panel with $\eps_1=1,\eps_2=-1,\iota=0$ breaks the double double root into a pair of complex double root leading to oscillations; the right panel with $\eps_1=\eps_2=0$, $\iota=1$ illustrates $1:2$-resonant back-coupling. Also shown in left and right panel are the linear spreading speed (vertical, magenta) and the predicted spreading speed associated with the double double root (left) and the resonant speed discussed in \S\ref{s:res}.}\label{f:ressliver}
\end{figure}

For simplicity, we focus here on the real-coefficient setting, introducing however differences in diffusivities between the in-phase and anti-phase modes,
\begin{equation}\label{e:cgl12}
    A_t=A_{xx}+\delta\bar{A}_{xx} +A -A|A|^2 +\gamma \bar{A} + \beta \bar{A}^2,
\end{equation}
with linearization for $A=u+\rmi v$,
\begin{equation}\label{e:cgl12lin}
        u_t=(1+\delta)u_{xx}+(1+\gamma)u,\qquad\qquad
        v_t=(1-\delta) v_{xx}+(1-\gamma)v,
\end{equation}
and dispersion relation.
\[
\begin{aligned}
    d(\lambda,\nu)&=d_u(\lambda,\nu)\cdot d_v(\lambda,\nu),\
    d_u(\lambda,\nu)&=\lambda-(1+\delta)\nu^2-(1+\gamma),\
    d_v(\lambda,\nu)&=\lambda-(1-\delta)\nu^2-(1-\gamma).
\end{aligned}
\]
Clearly, $\delta,\gamma>0$ favors the real mode spatio-temporally, leading to faster spreading. We therefore focus on $\delta\gamma<0$, allowing for a competition between the two contributions to spreading. Without loss of generality, after possibly switching $u$ and $v$, we may assume $\gamma>0>\delta>-1$.  As long as $|\delta|,|\gamma|<1$, we find spreading speeds
$
c_u=2\sqrt{(1+\delta)(1+\gamma)}$ and $c_v=2\sqrt{(1-\delta)(1-\gamma)}
.$
It is worth noting though that the dispersion relation possesses more double roots, arising through the collision of roots $\nu$ in the $u$ equation with roots $\nu_v$ in the $v$-equation, some of which may in fact induce a stronger instability and, nominally, faster spreading. To find those, we set both factors in $d$ to zero and find
\[
\lambda-(1+\delta)\nu^2-(1+\gamma)=0,\qquad \lambda-(1-\delta)\nu^2-(1-\gamma)=0.
\]
In the comoving frame of speed $c$, we find in this fashion a double root $  (\lambda_\mathrm{ddr},\nu_\mathrm{ddr})$
and a spreading speed $c_\mathrm{ddr}$ which leads to $\lambda_\mathrm{ddr}=0$,
\begin{equation}
    (\lambda_\mathrm{ddr},\nu_\mathrm{ddr})=\left(\frac{\gamma-\delta-c\sqrt{-\gamma\delta}}{\delta},\sqrt{-\frac{\gamma}{\delta}}\right),\qquad c_\mathrm{ddr}=\frac{\gamma-\delta}{\sqrt{-\gamma\delta}}.
\end{equation}
Of course, this linear spreading speed is not relevant in the sense that it does not lead to spreading at this speed. Technically, the second implication of Lem.~\ref{lem:lssfromPDR} fails. Near this double root, the dispersion relation has a quadratic expansion,
\begin{multline}\label{e:ddrexp}
d\left(\lambda_1-c_\mathrm{ddr}(\nu_\mathrm{ddr}+\nu_1\right),\nu_\mathrm{ddr}+\nu_1)=
-\frac{ \left(\lambda_1
   \sqrt{-\gamma  \delta }+\nu_1
   (-2 \gamma  \delta +\gamma
   +\delta )\right) \left(\lambda_1
   \sqrt{-\gamma  \delta }+\nu_1  (2
   \gamma  \delta +\gamma +\delta
   )\right)}{\gamma  \delta }\\+\rmO\left(|\lambda_1|^3+|\nu_1|^3\right).
\end{multline}
In particular, $\partial_\lambda d=0$ at the double root, which violates the nondegeneracy condition. In fact, this double root is a double solution to the system $(d,\partial_\nu d)=0$, a double double root.
Checking if this double root corresponding to $c_\mathrm{ddr}$ is pinched, we find after some sign checking that the associated root $\nu$ is the larger of the two roots of $d_u(-c_\mathrm{ddr} \nu,\nu)$ if and only if $\gamma+\delta+2\gamma\delta>0$ and it is the smaller of the two roots of $d_v(-c_\mathrm{ddr} \nu,\nu)$ if and only if $\gamma+\delta-2\gamma\delta<0$. The double root is then pinched if
\[
(\gamma,\delta)\in Q_\mathrm{ddr}=\left\{ (\gamma,\delta)\mid\frac{-\gamma}{1-2\gamma} <\delta<\frac{-\gamma}{1+2\gamma}\right\}, \text{ so that }
(\gamma+\delta+2\gamma\delta)(\gamma+\delta-2\gamma\delta)<0,
\]
a region depicted in Fig.~\ref{f:ressliver}. One also verifies that precisely in this region, the linear spreading speed is less than the double double root speed,
\[
c_\mathrm{lin}=\max\left\{2\sqrt{(1+\delta)(1+\gamma)},2\sqrt{(1-\delta)(1-\gamma)}\right\}<c_\mathrm{ddr}=\frac{\gamma-\delta}{\sqrt{-\gamma\delta}}, \qquad \text{ for all } (\gamma,\delta)\in Q.
\]
Adding linear or nonlinear coupling terms, for instance through the parameters $\eps_{1/2}$ and $\iota$,
\begin{equation}\label{e:cmc}
    \begin{aligned}
        u_t&=(1+\delta)u_{xx}+(1+\gamma)u-u(u^2+v^2) +\eps_1 v + \iota(u^2-v^2)\\
        v_t&=(1-\delta)v_{xx}+(1-\gamma)v-v(u^2+v^2) + \eps_2 u - \iota u v
    \end{aligned}
\end{equation}
exhibits the effect of these double double roots through an acceleration to the speed $c_\mathrm{ddr}$ when $\eps_1\neq 0$ and a more subtle resonant speed when $\iota\neq 0$ %
; see Fig.~\ref{f:ressliver} for simulations and \S\ref{s:res} for a discussion.
\textcolor{black}{\begin{remark}[Wave equation asymptotics for anomalous spreading]\label{r:wave}
        When $\eps_1\neq 0$ and $\eps_2=0$, the double double root spreading in an exponential weight and frame moving with $c_\mathrm{ddr}$ leads to diffusive transport 
        \[u_t=(1+\delta)u_{xx} + c_u u_x+v, \qquad v_t=(1-\delta)v_{xx}-c_v v_x,
        \]
        which, neglecting diffusive corrections, gives the wave equation 
        \[
        (\partial_t - c_u \partial_x)(\partial_t + c_v \partial_x)u=0,
        \]
        with $c_u=(-\gamma\delta)^{-1/2}(-\gamma-\delta+2\gamma\delta)<0$ and $c_v=(-\gamma\delta)^{-1/2}(-\gamma-\delta-2\gamma\delta)>0$; see Fig.~\ref{f:instability_pdr_dns}, right panels, for an illustration.
\end{remark}}
\textbf{Example 6: FitzHugh-Nagumo, real double roots, oscillatory fronts.}
\textcolor{black}{Our last example with a closed-form expression for spreading speeds does not appear to be documented in the literature and explores the  limits of  algebraic approaches.}

The linearization at the unstable origin of~\eqref{e:fhn}, when $a<0$, is
\begin{equation}
    \begin{aligned}
u_t&=u_{xx}-a u -v,\\
v_t&=\eps(u-\gamma v),
    \end{aligned}\label{e:fhnl}
\end{equation}
We focus here on the case $\gamma=0$ which allows for an explicit solution of the polynomial equations for double roots and spreading speeds.  The dispersion relation in the comoving frame is
\[
d_c(\lambda,\nu)=c\nu^3+(c^2-\lambda)\nu^2-(2c\lambda+ac)\nu+\lambda^2+a\lambda +\eps,
\]
with derivative at $\lambda=0$, $
\partial_\nu d_c(0,\nu)=c(-a+2c\nu+3\nu^2).$
This latter equation is solved for $c>0$ as $c=(a-3\nu^2)/(2\nu)$ and the result substituted into $d_c(0,\nu)=0$ gives  a quartic equation in $\nu$,
$
4\eps-a^2+2 a\nu^2 + 3\nu^4=0,
$
with solutions and associated speeds (listing only potentially positive speeds)
\begin{align}
  \nu_\mathrm{lin}&=- \frac{1}{\sqrt{3}} \sqrt{-a-2\sqrt{a^2-3\eps}},\qquad &
  c_\mathrm{lin}&= \frac{\sqrt{3}\left(-a-\sqrt{a^2-3\eps} \right)}{\sqrt{-a-2\sqrt{a^2-3\eps}}},\label{e:fhnspsp}\\
  \nu_\mathrm{lin}&=-\frac{1}{\sqrt{3}} \sqrt{-a+2\sqrt{a^2-3\eps}},\qquad &
  c_\mathrm{lin}&= \frac{\sqrt{3}\left(-a+\sqrt{a^2-3\eps} \right)}{\sqrt{-a+2\sqrt{a^2-3\eps}}}.\label{e:fhnspsp2}
\end{align}
For $a^2>4\eps$ there is precisely one negative real root, for $4\eps>a^2>3\eps$ there are two, and for $a^2<3\eps$ there are none. For $a^2<3\eps$, we necessarily find spreading caused by a pair of complex conjugate pinched double roots not visible in this calculation. For $0<\eps\ll 1$, perturbation from $\eps=0$ immediately shows that the first root in~\eqref{e:fhnspsp} indeed gives the correct spreading speed, that is, it is associated with the crossing of a right-most pinched double root in $\C$  as $c$ is decreased.
In summary, we have shown that for all $\eps>0$, sufficiently small, and $\gamma=0$, we have a spreading speed $c$ with  real crossing double root $\lambda=0$ and associated $\nu$ given by~\eqref{e:fhnspsp2}. The stationary dynamics and monotone leading edge are in contrast to the observed pattern formation in the wake; see Fig.~\ref{f:cglfhn}.

\textcolor{black}{\textbf{Example 7: General 2-species reaction-diffusion systems.}
 More general reaction-diffusion systems with two species 
 \begin{equation}
    \begin{aligned}
    u_t&=u_{xx} + c u_x + \alpha  u + \beta v,\\  
    v_t&=D v_{xx} + c v_x + \gamma  u + \delta v,
    \end{aligned}\label{e:rd2}
\end{equation}
do not appear to allow for closed-form computation of linear spreading speeds. One finds a dispersion relation as a fourth-order polynomial in $\nu$. Double roots can then be obtained as roots of the resultant which is an unwieldy fifth order polynomial in $\lambda$ without explicit roots or ways to determine when roots cross the imaginary axis as $c$ varies. The examples presented here are thus exceptions to the rule that spreading speeds are typically not explicit, with the last example of the FitzHugh-Nagumo linearization with $\gamma=0$ a novel addition to the list of exceptions.
}




\section{Nonlinear marginal stability}\label{s:nlms}

Invasion is usually not a small-amplitude phenomenon, governed entirely by the linearization at the unstable state. It is rather mediated by nonlinear propagating structures. We nevertheless mimic the viewpoint   detailed in \S\ref{s:lms},  starting with an informal statement of the marginal stability conjecture that reflects the speed selection process from compactly supported initial conditions for linear equations in \S\ref{s:3.0}. We then provide context, identifying the selected front as member of a family of nonlinear fronts, \S\ref{s:3.1}, characterized by marginal pointwise stability, \S\ref{s:3.2}, and conclude with practical strategies in \S\ref{s:3.3} and examples in \S\ref{s:3.4}. We focus throughout on \emph{rigidly propagating} fronts leaving behind a stable equilibrium, and turn to time- and space-periodicity in \S\ref{sec: modulated}.

\subsection{Marginal stability and selection}\label{s:3.0}

We consider the nonlinear parabolic system
\begin{equation}\label{e:parnl}
u_t=\mathcal{P}(\partial_x)u + f(u), \qquad u\in\R^N, \quad \mathcal{P}(\nu)=\sum_{j=0}^{2m}\mathcal{P}_j\nu^j\in \R^{N\times N} \text{ polynomial of degree } 2m.
\end{equation} We assume that there is an  unstable equilibrium $u=0$, so $f(0)=0$, and assume $f'(0)=0$  after possibly absorbing $f'(0)$ into $\mathcal{P}(\partial_x)$. We also assume that there is a stable equilibrium $u_-\in\R^n$, that is,
$\mathcal{P}(0)u_-+f(u_-)=0$ and all eigenvalues of $\mathcal{P}(\rmi k)+f'(u_-)$ have negative real part for all $k\in\R$. In particular, setting $k=0$, we see that $u_-$ is a stable equilibrium in the ODE $u_t=\mathcal{P}(0)u+f(u)$. Clearly this includes the FKPP equation where $\mathcal{P}(\partial_x)=\partial_{xx}+1$,  $f(u)=-u^3$, with $u_-=1$ and negative eigenvalues $-1-k^2$ of  the scalar $\mathcal{P}(\rmi k)+f'(u_-)$. We could allow more general $f=f(u,u_x,\ldots)$ with straightforward modifications.

In the comoving frame of speed $c$, we look for front-like equilibria $u_*(x;c)$, which solve the ODE
\begin{equation}\label{e:twode}
\mathcal{P}(\partial_\xi)u+c\partial_\xi u + f(u)=0, \qquad u(\xi)\to 0 \text{ for } \xi\to+\infty, \quad  u(\xi)\to u_- \text{ for } \xi\to-\infty.
\end{equation}
As we shall see in \S\ref{s:3.1}, front solutions typically exist for a range of speeds $c$, and the question then is which of those fronts are actually selected, that is,  observed when starting from initial conditions with support in $x<0$, say. This ``selection'' of fronts was made precise in~\cite[Def.~1]{as1}, and states, mildly rephrased, the following.
\begin{definition}[Selected fronts]\label{def: selected front}
    A speed $c_*$ and associated front $u_*(x;c_*)$ are \emph{selected} if an open class of steep initial conditions propagates with asymptotic speed $c_*$ and stays close to translates of the front $u_*$. More precisely, we require that there exists a non-negative continuous weight $\rho : \R \to \R$ and,  for any $\eps>0$, a set of initial data $\mathcal{U}_\eps \subseteq L^\infty (\R)$ such that:
	\begin{itemize}
		\item[(i)] for any $u_0 \in \mathcal{U}_\eps$, there exists a function $h(t) = \mathrm{o}(t)$  such that the solution $u(x, t)$ to~\eqref{e:parnl} with initial data $u_0$ satisfies, for $t$ sufficiently large,
		\begin{align}
		\| u(\cdot + c_* t + h(t), t) - u_* (\cdot) \|_{L^\infty(\R)} < \eps; \label{e: basin propagation}
		\end{align}
		\item[(ii)] there exists $u_0 \in \mathcal{U}_\eps$ such that $u_0 (x) = 0$ for all $x>0$ sufficiently large;
		\item[(iii)]  $\mathcal{U}_\eps$ is open in the topology induced by the norm $\| g \|_{\rho} = \| \rho g \|_{L^\infty}$.
	\end{itemize}
\end{definition}
\textcolor{black}{We reemphasize that most of Def.~\ref{def: selected front} is simply a statement on orbital asymptotic stability, that is, asymptotic stability of the family of translates, in suitable norms --- with the exception of (ii), which insists that the basin of attraction contains steep initial conditions, more precisely initial conditions that vanish in the leading edge $x>0$.}

The key observation we wish to highlight here is that fronts selected in this sense can be found as marginally stable members of the larger family of fronts. We therefore shall discuss in \S\ref{s:3.2} concepts of pointwise and spectral stability of nonlinear fronts, in contrast to the   stability of the trivial state $u=0$ in \S\ref{s:lms}. 
%

%
\begin{theorem}[Marginal stability conjecture~\cite{as1,avery2}]\label{t:msc}
Suppose that the front $u_*(x;c_*)$ at speed $c_*$ is marginally stable due to
\begin{itemize}
    \item \emph{(pulled)} a simple branch point {\color{black}$(\lambda,\nu)=(0,-\eta_\mathrm{lin})$} at $\lambda=0$; \emph{or}
    \item \emph{(pushed)} a simple eigenvalue at $\lambda=0$.
\end{itemize}
Then this front is selected in the sense made precise above in Def.~\ref{def: selected front}.
\end{theorem}
\textcolor{black}{We believe a converse theorem should hold as well, stating that: the basin of attraction of fronts with $c<c_*$ is not open since they are strongly unstable, and the basin of attraction of fronts with speed $c>c_*$ does not contain steep initial data as in Def.~\ref{def: selected front}(ii).}

The marginal stability assumption in pushed and pulled cases is \emph{spectral} (in a fixed exponentially weighted space --- see \S\ref{s:2.3}), rather than pointwise.
In the pulled case, the positional shift is $h(t)=-\frac{3}{2\eta_\mathrm{lin}}\log(t)+\rmO(1)$. In the pushed case it is $h(t)=h_0+\rmO(\rme^{-\delta t})$ for some $\delta>0$. Higher order asymptotics of the positional shift in the pulled case have been derived with matched asymptotics~\cite{EbertvanSaarloos} and proven when comparison principles are available~\cite{NRRkppasy}. The marginal stability of pushed fronts has been established much earlier, going back to~\cite{Satt1977,HadelerRothe}.
%
%
%
%
We will sketch the proof of the marginal stability conjecture in \S\ref{s:nlmsp}, explaining in particular how the marginal stability is necessary and sufficient for a matching of a diffusive leading edge to the appropriately shifted front profile. The remainder of this section focuses on first the existence of a family of fronts and then the characterization of marginal stability through spectral properties.

\subsection{Families of fronts}\label{s:3.1}

We may write \eqref{e:twode}  as a first-order differential equation for $U=(u,\partial_\xi u,\ldots \partial_\xi^{2m-1}u)\in\R^{2mN}$,
\begin{equation}\label{e:twsys}
         U_{0,\xi}=U_1,\quad
         U_{1,\xi}=U_2,\qquad\ldots\qquad,
         U_{2m-1,\xi}=\mathcal{P}_{2m}^{-1}\left(-cU_1-f(U_0)  - \sum_{j=0}^{2m-1}\mathcal{P}_jU_j,  \right)
\end{equation}
where we  look for heteroclinic orbits in~\eqref{e:twsys}  connecting  $U_-=(u_-,0,\ldots,0)$ and $U_+=(0,0,\ldots,0)$.

\textbf{Counting arguments.}
Heteroclinic orbits lie in the intersection of the stable manifold of $U_+$ and the unstable manifold of $U_-$. A first step in understanding those intersections consists of determining the dimension of those manifolds. The dimension of the unstable manifold of $U_-$ is $mN$, that is, precisely half the dimension of the ambient space, independent of $c$. To see this, consider the matrices $M_\lambda$ that arise when writing the linearization
$\mathcal{P}(\partial_\xi)u+c\partial_\xi u + f'(u_-)u=\lambda u$ as a first-order system $U'=M_\lambda U$:
\begin{itemize}
\item For $\lambda\gg 1$, the equation can be written as a small perturbation, after suitable rescaling, of the principal part $\mathcal{P}_{2m}\partial_\xi^{2m}u-\lambda u$, which does not possess solutions of the form $\rme^{\rmi k \xi}$ when
$(-1)^{m+1}\mathcal{P}_{2m}$ is negative, that is, when the equation is well-posed.
\item During a homotopy from $\lambda\gg1$ to $\lambda=0$, $M_\lambda$ does not have purely imaginary eigenvalues  $\nu=\rmi k$ since those correspond to a solution $u\sim\rme^{\lambda t+ \rmi k \xi}$ and hence instability of $u_-$.
\end{itemize}
As a consequence, the unstable dimension is $mN$ during the entire homotopy.
The dimension of the unstable manifold of $U_+=0$ does in general depend on the equation and on $c$. Mimicking the homotopy that we outlined for $U_-$, crossings of eigenvalues $\nu$ of $M_\lambda$ now may occur due to unstable spectrum, that is, we cannot exclude solutions of the form  $u\sim\rme^{\lambda t+ \rmi k \xi}$ for some $\lambda>0$. A homotopy in $c$ at $\lambda=0$ gives some general insight:
\begin{itemize}
    \item For $c\to\infty$, only the highest derivatives  contribute, distributed in each component according to $(-1)^{m-1}\nu^{2m}+c\nu=0$, for all $N$ components of the system, resulting in $N$ zero roots, $mnN$ stable roots, and $(m-1)N$ unstable roots, all at size $c^{1/(2m-1)}$. At next order, one finds that the zero roots split  according to $c\nu+f'(0)=0$. If we denote by $i_f$ the number of unstable eigenvalues of $f'(0)$, that is the unstable dimension in the \emph{kinetics}, we find that the dimension of the stable manifold is $mN+i_f$ when $c$ is sufficiently large.
    \item Eigenvalues cross the imaginary axis precisely when there are solutions of the form $u\sim \rme^{\rmi k (\xi-ct)}$ for some $c$ and some $k\neq 0$. We shall exclude this type of more complex instability in the sequel.
\end{itemize}
Altogether,  we find that the dimension of the stable manifold of $U_-$ and the dimension of the unstable manifold of $U_+$ add up to $2mN+i_f$, predicting generically  transverse  intersections of dimension $i_f$. For an instability with only one real positive eigenvalue, as say in the FKPP equation, this predicts generically the unique intersection that we find in the phase plane for any speed $c$.

\textbf{Existence at large speeds.}
For large speeds, we rescale $c\partial_\xi=\partial_y$ to find  formally the pure backward-in-time kinetics $u_y=-f(u)+\rmO(\eps)$, with $\eps=c^{-1}$, where the higher-order terms contain higher derivatives in $y$. Rigorously, we invoke Fenichel's geometric singular perturbation theory~\cite{fenichel}, setting for instance $U=(u,u_1,u_{2m-1})$, defined through
$u_{j}=\eps^{\frac{j}{2m-1}}\partial_\xi^ju$, $j=1,\ldots,2m-1$, and desingularizing time via $\partial_\xi=\eps^{\frac{-1}{2m-1}}\partial_y$, to obtain
\begin{equation}\label{e:gstp}
u_y=u_1,\quad u_{1,y}=u_2,\quad \ldots\quad  u_{2m-1,y}=(-1)^m\left(u_1+\eps^{\frac{2m}{2m-1}}f(u)\right)+\rmO(\eps u_1,\eps^{\frac{2m-2}{2m-1}} u_2,\eps^{\frac{2m-3}{2m-1}} u_3,\ldots).
\end{equation}
At $\eps=0$, there is an $N$-dimensional slow manifold $u_1=u_2=\ldots u_{2m-1}=0$, with linearization having roots $0$ of multiplicity $N$ and roots of unity $(-1)^{\frac{m}{2m-1}}$ off the imaginary axis. Expanding the manifold in $\eps$ gives $u_1=-\eps^{\frac{2m}{2m-1}}f(u)$ at leading order and reduced flow $u_y=-\eps^{\frac{2m}{2m-1}}f(u)$, or
$cu_\xi=-f(u)$ as the formal calculation, dropping higher $\xi$-derivatives, suggests.

As a consequence, we find fast traveling waves at leading order as heteroclinic orbits in the kinetics $u_t=f(u)$, in slow reversed time $\xi=-\eps t$.  Such heteroclinics  exist quite generally, following for instance trajectories in the unstable manifold of $u=0$. In summary,  we typically do find traveling waves for $c\gg 1$.

\textbf{Homotopies.} We wish to continue from $c\gg 1 $ to $c=\rmO(1)$ and therefore  focus on the case $i_f=1$. Of course, there are numerous potential obstacles to guaranteeing the existence of heteroclinics for all positive $c$ and we list some possibilities here.
\begin{itemize}
    \item \emph{Front splitting:}
        A simple example arises in the Nagumo equation, where the front connecting $u=1$ to $u=0$ splits into a front connecting $u=1$ and $u=-a$ and a front connecting $u=-a$ and $u=0$; see Fig.~\ref{f:nagumophaseportraitsplit}. More generally, the derivative of the front could split its support as one of the options in the concentration-compactness trichotomy~\cite{lions_1984}.
        \begin{figure}
        \centering\includegraphics[width=\textwidth]{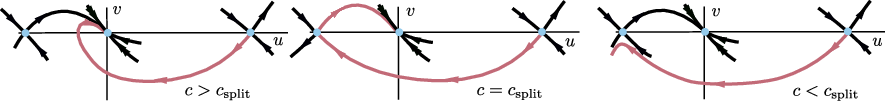}
        \caption{Phase portrait sketches in the $u-u_x$-plane of~\eqref{e:kpptw}, $f(u)=u(1-u)(u+a)$ for $a<\frac{1}{2}$, with $c$ decreasing from left to right. The front connecting $u=1$ to $u=0$ ceases for exist for $c<c_\mathrm{split}=(1-a)/\sqrt{2}$.} \label{f:nagumophaseportraitsplit}
        \end{figure}
\item \emph{Loss of transversality or saddle-nodes of fronts} seems rare; we are not aware of bifurcations of invasion fronts due to a loss of transversality in the literature other than \cite[\S2,Ex.4]{pnp}.
\item \emph{Bifurcations of asymptotic states} can lead in particular to staged invasion fronts; see \cite[\S2]{pnp}. Somewhat subtle examples, motivated for instance by multi-dimensional spreading in coupled amplitude equations, were studied for instance in~\cite{HSaccelerated,FSlocked}.
A simple example is
\begin{equation}
    \label{e:kpppf}
    u_t=u_{xx}-u(u+1)(u-1)((u-1)^2-\mu), \qquad \mu\sim 0,
\end{equation}
with unstable equilibrium $u=0$ and linear spreading speed $c_\mathrm{lin}=2\sqrt{1-\mu}$. One can readily show that for $\mu>0$, the front selects the bifurcating equilibrium  $1-\mu$; see Fig.~\ref{f:scalpf}.
\begin{figure}
    \raisebox{0.2in}{\includegraphics[width=0.25\textwidth]{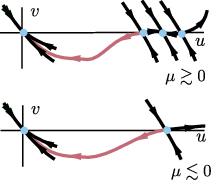}}
    \hfill\raisebox{0.in}{\includegraphics[width=0.25\textwidth]{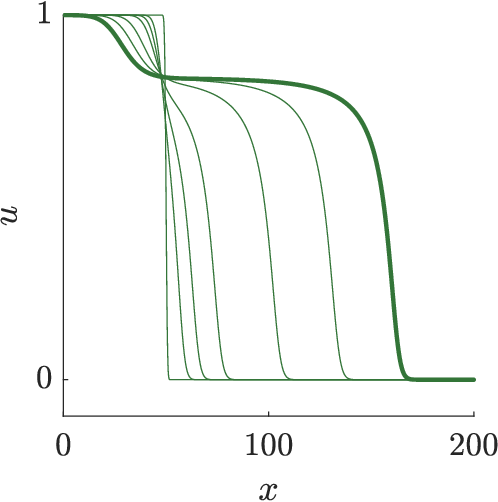}}
    \hfill\raisebox{0.2in}{\includegraphics[width=0.35\textwidth]
{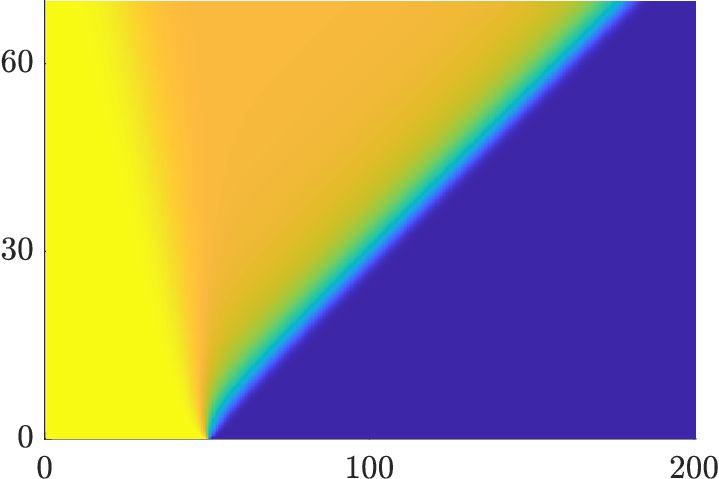}}
    \caption{The invasion front of~\eqref{e:kpppf} with $\mu=0.03$ selects the lowest of the three equilibria emerging near $u=1$: phase portrait of traveling-wave equation (right) showing how only the lower, closest of the equilibria connects to the origin with a front; direct simulation with step function $u=1$ for $x<50$, snapshots in time (center) and space-time plot (right), show the selection through the formation of a secondary front propagating slowly to the left; see text for details.}\label{f:scalpf}
\end{figure}
A more subtle example, discussed also in \cite{pnp}, arises in coupled FKPP-type equations,
    \begin{equation}\label{e:frontmorsechange}
        u_t=u_{xx}+u-u^3,\qquad v_t=dv_{xx}+g(\mu,u)v - v^3.
    \end{equation}
    Bifurcations occur when $g(\mu,1)=0$ or $g(\mu,0)=0$, leading to instabilities in the wake or additional instabilities in the leading edge.
    More intricate bifurcations occur in situations with differential transport speeds. Those may arise naturally as amplitude equations near an oscillatory Turing bifurcation~\cite[\S2.2]{SRadial},
    \begin{equation}\label{e:cpw}
        u_t=u_{xx}+c_\mathrm{g} u_x+\mu u-u(u^2+\beta v^2),\qquad v_t=u_{xx}-c_\mathrm{g}v_x+\mu v-v(v^2+\beta u^2).
    \end{equation}
    The linearization at $u=v=0$ in a comoving frame of speed $c$  possesses a pair of imaginary eigenvalues that crosses the imaginary axis as $c$ decreases through $|c_\mathrm{g}|$.
\end{itemize}
\subsection{Stability of fronts}\label{s:3.2}
We now turn to the linearization at a family of fronts $u_*(\xi;c)$ parameterized by the wave speed $c$,
\begin{equation}\label{e:linfront}
    v_t=\mathcal{P}(\partial_\xi)v + c\partial_\xi v+ f'(u_*)v=:\mathcal{L}_c v.
\end{equation}
We can now in fact mirror the the strategy from \S\ref{s:lms} to identify marginal stability:
\begin{enumerate}
    \item determine pointwise stability of $\mathcal{L}_c$ in the leading edge;
    \item find the speed and associated front with $\mathcal{L}_c$ marginally stable in the leading edge.
\end{enumerate}
The idea again is that an instability in the leading edge will lead to an acceleration of the front, while stability leads to decay in the leading edge and a slow-down.
In order to determine stability, we would like to understand the ``heat kernel'' $H(t,\xi,y)$, solution to~\eqref{e:linfront} with $v(0,\xi)=\delta(\xi-y)$, which one in turn calculates from the pointwise resolvent $G_\lambda(\xi,y)$,
\begin{equation}\label{e:resx}
    (\mathcal{L}_c-\lambda)G_\lambda(\xi,y)=\delta(\xi-y),
\end{equation}
via inverse Laplace transform as in~\eqref{e:resolvent}.  Analyticity of $G_\lambda$ in $\lambda$ then allows one to deform the integral contours in the inverse Laplace transforms up to the rightmost singularity of $G_\lambda$ in the complex plane. Marginal stability as in (ii) then is equivalent to a singularity of $G_\lambda$ being located on the imaginary axis. We outline below some general strategies to analyze the pointwise resolvent $G_\lambda$ but emphasize that in general, the spatio-temporal growth of disturbances can be rather complex; see for instance~\cite{brevdo1,brevdo2}, and~\cite{HSaccelerated} for examples where such intricate spatio-temporal growth governs front invasion.


\textbf{Pointwise analysis.} The most refined results on the structure of $G_\lambda$ are often obtained using dynamical systems techniques, reformulating~\eqref{e:resx} as a first-order system as in~\eqref{e:Ufirstorder},
\begin{equation}\label{e:res1stx}
    U_\xi=M_\lambda(\xi)U+\delta(\xi-y)I.
\end{equation}
The subspaces $E^\mathrm{u}_{\lambda,y}$ and $E^\mathrm{s}_{\lambda,y}$ of initial conditions at $\xi=y$ that lead to bounded solutions as $\xi\to -\infty$ and $\xi\to+\infty$, respectively, are of dimension $mN$ for $\Re\lambda\gg 1$ and depend analytically on $\lambda$. There are then two main sources for singularities of $G_\lambda$:
\begin{itemize}
    \item singularities at $\infty$, related to branch points: the limiting subspaces lose analyticity as elements of the Grassmannian;
    \item singularities at $y$, related to intersections of $E^\mathrm{u}_{\lambda,y}$ and $E^\mathrm{s}_{\lambda,y}$.
\end{itemize}
The latter singularities are referred to, depending on context and specific situation, as extended point spectrum~\cite{AbsoluteSpecArndBjorn}, as resonances, as embedded eigenvalues, or as zeros of the analytic extension of the Evans function; see~\cite{KapitulaPromislow,FiedlerScheel,sandstab,ArndBjornEvansfcn} for background on the Evans function and its extensions. The former singularities are simply the branch points of the dispersion relation identified as sources of instability in \S\ref{s:lms}. The fact that one can largely decouple these two sources of singularities relies on extending stable and unstable subspaces defined for the limiting problems at $\pm\infty$ to finite $\xi$ using fixed point arguments similar to stable and strong stable manifold theorems; see for instance~\cite{ZumbrunGardner,KapitulaSandstede}.

We remark that it can be useful to analyze existence and stability jointly, that is, to append the stability ODE to the existence problem, finding an autonomous nonlinear problem
\begin{equation}\label{e:exstab}
    \mathcal{P}(\partial_\xi)u + c \partial_\xi u +f(u)=0,\qquad  \mathcal{P}(\partial_\xi)v + c \partial_\xi v +f'(u)v=\lambda v,
\end{equation}
in which one of course is interested in a specific solution $u$, only. Stable and unstable subspaces then are found simply within stable and unstable manifolds.  The $v$-equation can also be analyzed with a specific focus on stable or unstable subspaces by writing subspaces as graphs and studying the resulting matrix  Riccati equations; see for instance~\cite[\S4]{ArndBjornEvansfcn} and~\cite{rdr}.

Marginal stability  occurs at $c=c_*$ when $G_\lambda$ is analytic in $\Re\lambda>0$ for $c>c_*$ and  either
\begin{enumerate}
    \item $G_\lambda$ has a branch point singularity for some $\lambda\in\rmi\R$ for $c=c_*$; or
    \item $G_\lambda$ has a pole for some $\lambda\in\rmi\R$ for $c=c_*$.
\end{enumerate}
The first case (i) is commonly referred to as a pulled front, since the marginal stability is induced by marginal stability ahead of the front, at the trivial state, which can be thought of pulling the front into the unstable region. The second case is commonly referred to as a pushed front, since the marginal stability is related to the shape of the interface, which then pushes the front into the unstable region; see~\S\ref{s:nlmsp} for precise definitions.

\textbf{Spectral analysis.}
In most cases, singularities of $G_\lambda$ relate to spectrum of the operator $\mathcal{L}_c$ in appropriately weighted function spaces. Therefore, consider a weight  $\omega_\eta(\xi)=1,\xi<-1$, $\omega_\eta(\xi)=\rme^{\eta \xi}, \xi>1$, $\omega_\eta(\xi)>0$ and smooth, with associated function spaces $L^p_\eta$ with norm
$
\|u\|_{p,\eta}=\|\omega u\|_{L^p}.
$
It is often possible to choose $\eta>0$, so that in $L^p_\eta$,  the singularities of $G_\lambda$ that detect marginal stability correspond to spectrum, that is, $\Re\mathrm{spec}\,\mathcal{L}_c<0$ for $c>c_*$, and (compare Figs.~\ref{f:kppphaseportrait} and~\ref{f:nagumophaseportrait})
\begin{itemize}
    \item $\mathcal{L}_c$ has essential spectrum on $\rmi\R$ for $c=c_*$ (pulled fronts); or
    \item $\mathcal{L}_c$ has point spectrum on $\rmi\R$ for $c=c_*$ (pushed fronts).
\end{itemize}
In order to obtain the correct value of $c_*$, it is however necessary to choose an optimal $\eta$, pushing the essential spectrum as far as possible to the left. This minimization over $\eta$ of the maximum real part of $\lambda$ in the spectrum naturally yields the branch points of the dispersion relation as discussed in \S\ref{s:lms}.

We emphasize here that spectral stability is a stronger condition than pointwise stability. Clearly, spectral stability implies the existence of a bounded resolvent in $\Re\lambda\geq 0$, analytic in $\lambda$, from which one readily concludes the analyticity of the resolvent kernel in $\lambda$. On the other hand,  instability in a fixed norm may amount to instabilities that do not correspond to singularities of the pointwise resolvent, so that in general,  the speed detected by spectral marginal stability may be larger than the speed corresponding to pointwise marginal stability. Examples arise for instance when the right-most point of the absolute spectrum is not a branch point; see~\cite{FayeHolzerScheelSiemer} for an example.

\textbf{Stability for $c\gg1$.} As an illustration, we show stability in weighted spaces for large $c$, thus giving an a priori upper bound for \emph{nonlinear speeds} similar to the a priori bounds on linear speeds in \cite[Lem.~6.4]{HolzerScheelPointwiseGrowth}.
The existence and stability problem~\eqref{e:exstab} can be scaled jointly as in~\eqref{e:gstp} to find
\begin{equation}\label{e:gstpstab}
\begin{aligned}
u_y&=u_1,\quad u_{1,y}&=u_2,\quad\ldots\quad u_{2m-1,y}&=(-1)^m\left(u_1+\eps^{\frac{2m}{2m-1}}f(u)\right)+\rmO(\eps u_1,\eps^{\frac{2m-2}{2m-1}} u_2,\ldots),\\
v_y&=v_1,\quad v_{1,y}&=v_2,\quad\ldots\quad v_{2m-1,y}&=(-1)^m\left(v_1+\eps^{\frac{2m}{2m-1}}(f'(u)-\lambda)v\right)+\rmO(\eps v_1,\eps^{\frac{2m-2}{2m-1}} v_2,\ldots).
\end{aligned}
\end{equation}
For bounded $\lambda$ and $\eps\ll 1$, the $v$-equation is a small perturbation of a constant-coefficient equation with eigenvalues  0 and $(-1)^{\frac{m}{2m-1}}$, each of multiplicity $N$. In particular, there is $\delta>0$, small, so that there are $mN$ eigenvalues with real part less than $-\delta$ and $mN$ eigenvalues with real part larger than $-\delta$, so that the equation for weighted eigenfunctions $\rme^{\delta y} v$ is hyperbolic and does not possess bounded solutions. Fenichel's theory~\cite{fenichel}, guarantees that this hyperbolic structure persists for small $\eps$. Similarly, scaling  in $\lambda\gg 1$, gives hyperbolicity in the $v$-equation related to well-posedness of the parabolic equation and absence of eigenvalues for  $\lambda\gg 1$. Together, we find stability in a weight $\eta=\delta c^{\frac{1}{2m}} \xi$, $0<\delta\ll 1$.

\subsection{Practical considerations: nonlinear marginal stability}\label{s:3.3}
In order to find  selected invasion fronts numerically, both  speeds and states selected in the wake, one usually turns to direct simulations in a steady frame. We discuss here several strategies that can complement such a first investigation in light of the marginal stability criterion.

\begin{figure}
     \includegraphics[width=0.3\textwidth]{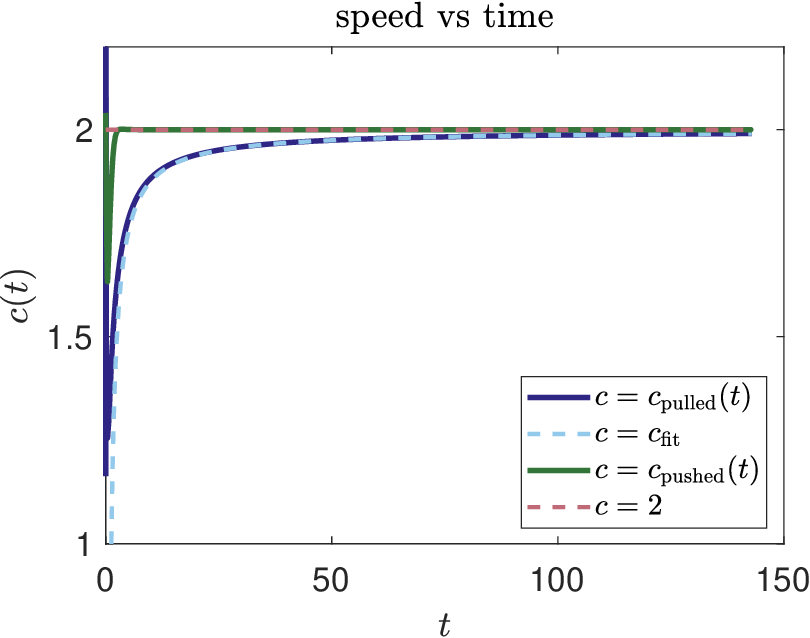}\hfill
     \includegraphics[width=0.3\textwidth]{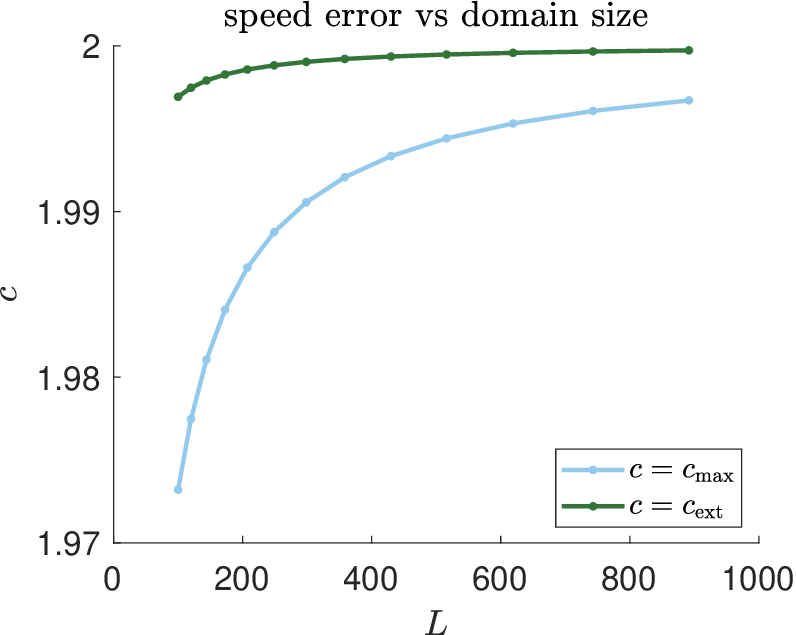}\hfill
     \includegraphics[width=0.3\textwidth]{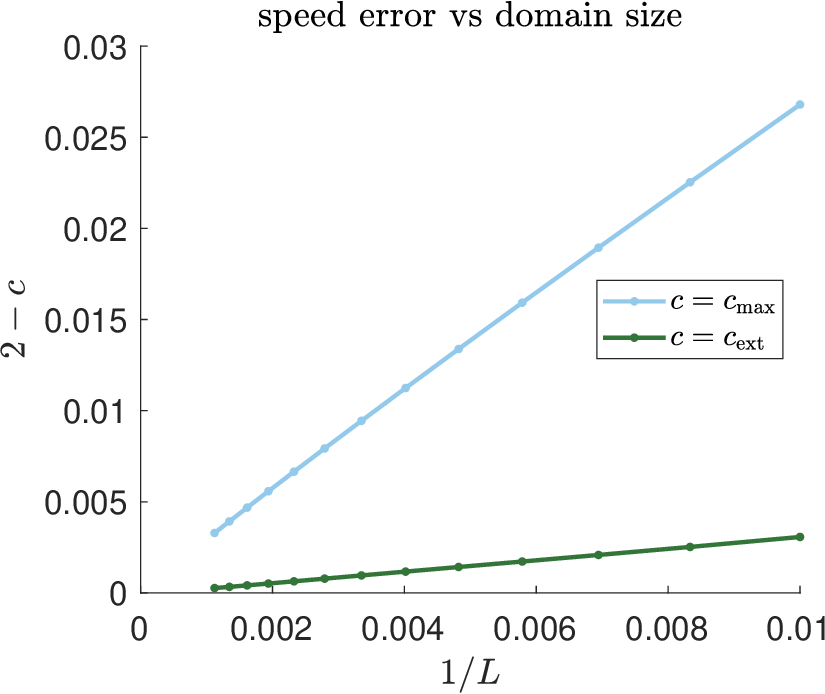}
    \caption{Speed of the front interface in the FKPP equation on a domain of size $L$ plotted against time for $L=300$, showing that the speed increases monotonically (left). Also shown for comparison is the rapid, exponential convergence in the case of a pushed front with appropriate cubic nonlinearity. The maximal speed $c_\mathrm{max}$ increases slowly with $L$ (center). Extrapolated speeds $c_\mathrm{ext}$, using that $c\sim c_0+c_1/t$, gives better approximations. The error of the speed decreases as $1/L$ (right).
    }\label{f:speed convergence}
\end{figure}

\textbf{Direct simulations.}  One first discretizes~\eqref{e:parnl} in a large bounded domain $x\in[0,L]$ with, say, Dirichlet boundary conditions, $\partial_{x}^{2j} u=0$, $j=0,\ldots m-1$ at $x=L$, Neumann boundary conditions $\partial_{x}^{2j+1} u=0$, $j=0,\ldots m-1$ at $x=0$, and initial conditions $u_0(x)$ supported in $x\in[0,1]$. One should keep in mind that the state in the leading edge is unstable and will thus amplify round-off errors. Simulations quickly become meaningless, determined by round-off error fluctuations. A first remedy is to adjust variables to guarantee that the leading-edge state is the trivial state $u\equiv 0$ rather than $u\equiv 1$, say, so that disturbances are represented with a resolution of $\sim10^{-308}$ rather than $\sim10^{-16}$ in double precision. The speed is then estimated from the front position  $x_*(t)$, defined through say thresholding,
\[
x_*(t)=\mathrm{argmax}_x\,\{|u(t,x)|\geq\delta\},\qquad \delta>0 \text{ not too small,}
\]
after interpolating on the grid.
The derivative $x_*'(t)$ converges to the speed $c_*$, albeit slowly, with rate $1/t$ for pulled fronts. Extrapolation, fitting $x_*'(t)=a_0+a_1/t$ for large $t$, can improve the estimates. Since speed estimates converge with rate $1/L$ and computational time increases with rate $L^2$, due to the increased domain size and length of the computation, the error $\Delta c$ decreases for pulled fronts slowly like
\[
\Delta c\sim \mathrm{(Cost)}^{-1/2}.
\]
Exploring regions of parameter space in this fashion can quickly become quite expensive. We show some of those convergence rates in Fig.~\ref{f:speed convergence} in the case of pulled fronts.

\textbf{Comoving frames.}
With a rough estimate $c_0$ of the speed, one could simulate in a comoving frame $u_t=\mathcal{P}(\partial_\xi)u +c_0\partial_\xi u + f(u)$, where fronts propagate with speed $c_*-c_0$, allowing for observations for times $T\sim L/|c_0-c_*|$. More sophisticated methods would improve $c_0$ continuously during the simulation using an appropriate phase condition~\cite{beyn2007phase}, or, one can simply shift the solution at discrete time steps, extending it by zeros  at the upstream edge of the domain. Yet more directly, since one is looking for a stationary solution in such a frame, one could then resort to a Newton method to compute the front solving for $u$ and $c$ in, for instance,
\begin{equation}\label{e:newtph}
    \mathcal{P}(\partial_\xi)u +c\partial_\xi u + f(u)=0,\qquad \int_{L_-}^{L_+}u = M(L_+-L_-),\qquad \xi\in [0,L],\qquad u'(0)=0,\ u(L)=0,
\end{equation}
where $[L_-,L_+]$ is an interval centered roughly at $\xi=L/2$ and $M$ a value of $u$ at the front interface.

Locally, this problem possesses a unique solution (roughly, the phase condition fixing $\int u$ selects a distinguished front out of the family of fronts with varying speed) and the speed $c$ will converge to the selected spreading speed, with finite-size corrections
\begin{equation}\label{e:finitesize}
    c(L)=c_*+\rmO(1/L^2)  \text{ for pulled fronts;}\qquad c(L)=c_*+\rmO(\rme^{-\eta L}) \text{  for pushed fronts.}
\end{equation}
Existence, local uniqueness, and error estimates were derived in~\cite{ADSS21}. Since effort for Newton iterations in a domain of size $L$ scales linearly, we find the improved error estimate for pulled fronts
\[
\Delta c\sim \mathrm{(Cost)}^{-2}.
\]
For pushed fronts, the error is again exponentially small in $L$ and hence in the computational effort.

An advantage of the methods proposed thus far is that they apply to both pulled and pushed fronts, although they are less effective for pulled fronts. We describe an algorithm that achieves exponential convergence for both pulled and pushed fronts in \S\ref{s:robust},
\[
\Delta c\sim \rme^{-\delta\cdot\mathrm{(Cost)}}, \quad \delta>0.
\]
\textbf{Linear predictions.}
Estimates for speeds obtained in this fashion can be compared to linear predictions from \S\ref{s:lms}. On the other hand, linear predictions can provide initial guesses for comoving frame speeds, in particular when a larger parameter space is explored.

\textbf{Global considerations.}
We expect (generically) selected fronts to be isolated, that is, there are discrete sets of fronts $u_{*,j}(\xi;c_j)$ with possibly different speeds. We discussed the example of the  Nagumo equation 
$
u_t=u_{xx}+u(1-u)(u+a),
$
with $a<1/2$ in the introduction: negative step-like initial conditions lead to spreading with the pulled speed $c_{*,1}=2\sqrt{a}$, while positive step-like initial conditions lead to pushed fronts with speed $c_{*,2}=(1+2a)\sqrt{2}>c_{*,1}$. In general, there can be many more fronts; see \cite{pnp}. There do not appear to be good strategies on how to find those fronts or determine their basin of attraction, although  $c\gg 1$ may be a good starting point for homotopies in many situations.



\textbf{Stability of fronts.} Establishing stability and detecting marginal stability in some generality is possible using either comparison principles where available or a perturbative scenario that we shall discuss in \S\ref{s:robust}. When computing spectra directly for a discretized problem, one clearly needs to be cautious about boundary conditions. Routines that phrase the eigenvalue  problem as an ODE can systematically compute and incorporate boundary conditions. Implementations then rely on computing determinants and find eigenvalues through winding number computations or root finders~\cite{stablab}, or, alternatively, use iterative methods that interpret the eigenvalue problem as the characteristic equation to a multi-term recursion~\cite{S23}. Rather than continuing in $c$ and detecting marginal stability, it is often advantageous to continue marginally stable fronts and speeds in system parameters starting in a somewhat well understood situation; see  \S\ref{s:robust} where we shall discuss effective algorithms that in particular track the pushed-to-pulled transition.



\subsection{Nonlinear spreading speeds and marginally stable fronts: examples}\label{s:3.4}

\textbf{Example: Scalar reaction-diffusion --- constructing fronts.}
A variety of methods can be used to construct fronts in  simple scalar FKPP-type equations. We present here a template that somewhat separates the phase-plane structure of the traveling-wave equation from the general strategy. We start with the construction of traveling waves, pursuing here the dynamical systems perspective of finding heteroclinic orbits. Consider therefore traveling waves to~\eqref{e:fkpp} that solve
\begin{align}
    u_\xi=v,\qquad \qquad
    v_\xi=-cv-f(u).\label{e:kpptw2}
\end{align}
where we assume $f(0)=f(1)=0$, $f(u)>0$ for $u\in(0,1)$, $f'(0)=a>0$, and $f'(1)<0$, and focus throughout on $c>0$. A traditional analysis finds fronts by constructing trapping regions, using the explicit form of the vector field. We instead explore a continuation from the large $c$ limit, as outlined in \S\ref{s:3.1}.
For $c=1/\eps$, $\eps\ll 1$, we can rescale $\xi =\eps  y$ and find
\begin{align}
    u_y=\eps v,\qquad \qquad    v_y=-v-\eps f(u).\label{e:kpptw3}
\end{align}
At $\eps=0$, the manifold of equilibria $\mathcal{M}_0=\{(u,v)|v=0\}$ is normally hyperbolic, that is, the linearization at each equilibrium possesses a simple zero eigenvalue associated with the tangent vector to $\mathcal{M}_0$, and a second eigenvalue off the imaginary axis, here -1. This manifold then continues smoothly as a smooth manifold $\mathcal{M}_\eps$ for small $\eps\neq 0$ given as a graph $\{v=\eps h(u;\eps)|\,u\in\R\}$ \cite{fenichel}. Substituting this expression into the equation for $v$ in~\eqref{e:kpptw3} and using the equation for $u$ we find
$
(\eps h_u)\cdot ( \eps v) = -\eps h - \eps f,
$ which gives $h(u;0)=-f(u)$, and, from the first equation, the vector field on $\mathcal{M}_\eps$ projected onto the $u$-direction as
\[
u_y=-\eps f(u) + \rmO(\eps^2),
\]
\textcolor{black}{Since $f$ has nondegenerate zeros at $0$ and $1$ they persist. They in fact do not depend on $\eps$, as is clear from the equation, and are connected by a heteroclinic orbit since $f(u)>0$ for $u\in(0,1)$, with explicit expansion in the case $f(u)=u(1-u)$,}
\[
u(\xi)=\frac{1}{1+\rme^{\eps \xi}}\left(1+\eps R(\eps \xi,\eps)\right), \qquad R \text{ smooth }.
\]
Note also that in this construction, both equilibria have one unstable eigenvalue, in addition to one unstable eigenvalue for the origin and one stable eigenvalue at $u=0$, reflecting the more general formula~\eqref{s:3.1}. We note that more direct methods are available to treat this singular limit; see for instance~\cite{RS07}.

{\color{black} The above construction gives existence of fronts for $\varepsilon$ sufficiently small, that is, $c = \frac{1}{\varepsilon}$ sufficiently large. One then naturally wants to continue these fronts as $c$ decreases. }Depending on $f$,  the heteroclinic orbit to~\eqref{e:kpptw2} may persist for all values of $c$, limit on a chain of heteroclinics, or become unbounded. To see this, one inspects the energy $E(u,v)=\frac{1}{2}v^2+F(u)$, $F'=f$, which is non-increasing for all non-negative $c$ and forces all bounded solutions to be heteroclinic for nonzero $c$. Compactness of a sequence of bounded solutions in the local topology then ensures convergence to a heteroclinic after appropriate shifts.

For a quadratic nonlinearity, one finds that trajectories in the branch of the unstable manifold of $(1,0)$ with $u<1$ are bounded, staying inside the sublevel set of the energy formed by the homoclinic trajectory to $(1,0)$ at $c=0$, so that the heteroclinic exists for all values of $c\neq 0$. In the case of a cubic nonlinearity, $f(u)=u(1-u)(u-a)$ with $a<0$, the global bifurcation picture depends on the magnitude of $a$. For $-a\leq -1$, the heteroclinic exists for all values of $c$, arguing as in the case of the quadratic nonlinearity. For $0>a>-1$, the heteroclinic connecting $u=1$ to $u=0$ limits on a heteroclinic connecting the two PDE stable equilibria $u=1$ to $u=-a$ and a ``subsequent'' heteroclinic connecting $u=-a$ to $u=0$ at $c_\mathrm{bist}=(1+a)/\sqrt{2}>0$. Before this heteroclinic splitting, the heteroclinic passes through the strong unstable manifold of $u=0$ at the pushed front speed $c_\mathrm{push}=(1-2a)/\sqrt{2}$; see Fig.~\ref{f:nagumophaseportrait} for the crossing of the strong stable manifold near $c_\mathrm{push}$ and Fig.~\ref{f:nagumophaseportraitsplit} for the phase portrait near $c=c_\mathrm{bist}$. Establishing details here relies on monotonicity properties of invariant manifolds in $c$.


\textbf{Example: Scalar reaction-diffusion --- stability of fronts.}
In the FKPP equation, fronts $u_*(\xi;c)$ are monotone in $\xi$  and have asymptotics $u_*(\xi;c)\sim \rme^{\nu_+\xi}$ for $c>2$, where $\nu_+=-c/2+\sqrt{c^2/4-1}$ is the larger of the two roots of the dispersion relation. In the exponentially weighted space $L^2_\eta$ with $\eta=-c/2$, the essential spectrum of the linearization is stable for all $c>2$ as one readily infers from the dispersion relation. For $c\gg1$, we saw fronts never have unstable point spectrum. Since the linearization is a Sturm-Liouville operator, potential point spectrum is real. We therefore only need to exclude point spectrum crossing at $\lambda=0$, in which case the only potential eigenfunction bounded at $\xi=-\infty$ is $u_*'(\xi;c)$. As a consequence, point spectrum crossing corresponds precisely to the presence of a front with steep decay, a pushed front, which one can exclude in the specific case of a quadratic nonlinearity; see \cite{AHR2} for more examples.

\textbf{Example: Local bifurcation.}
Fronts can be constructed near onset of instability using  reduction methods. Assuming for instance that $f(u;\mu)$ depends on a parameter $\mu\in\R$ and undergoes a generic bifurcation, such as saddle-node, transcritical, or cusp (with $\mu\in\R^2$), at $u=0,\mu=0$, one can try to find heteroclinic orbits nearby using a center-manifold reduction for the traveling-wave ODE. The ideas, in particular assumptions on $\mathcal{P}$, were generally described in~\cite[\S2]{SRadial}, without emphasis on spreading speeds. A more detailed analysis of stability for bifurcating fronts was carried out in the case of a transcritical bifurcation both for fronts faster than the linear spreading speed~\cite{RK1996} and for the critical, selected fronts~\cite{RK1998}. Nonlinear selection of the critical fronts was established in~\cite{avery2}.

\textbf{Example: Anomalous spreading.}
Supplementing FKPP with a simple diffusion equation
\begin{equation}\label{e:fkppdiff}
    u_t=u_{xx}+u(1-u)+\alpha v,\qquad v_t=dv_{xx},
\end{equation}
we still find the front solutions in $u_*(\xi;c)$ simply setting $v=0$. For $d>2$, however, the linearization at the front with speed 2 is
\[
u_t=(\partial_\xi+1)^2 u - 2u_* u + \alpha v,\qquad v_t=d(\partial_\xi+1)^2v + 2(1-d)(\partial_\xi+1)v+(d-2)v,
\]
so that in an exponential weight with growth $\eta=1$ the $u$-equation is marginally stable but the $v$-equation exhibits exponential growth with rate $\rme^{(d-2)t}$. Clearly, for $\alpha=0$ the front with speed 2 is the selected front, yet it is unstable in any exponentially weighted space\footnote{One can in this case choose different weights in the $u$- and $v$-equation, and obtain stabilization for $\alpha=0$, but more general examples would not allow for this strategy~\cite{FayeHolzerScheelSiemer}.}.  In fact,  $\alpha>0$ and $v\neq 0$ lead to an acceleration of the front in the $u$-equation. An analysis of the pointwise resolvent will of course not exhibit singularities in $\lambda>0$ since resolvent computations simply decouple. This phenomenon can also occur in examples with nonlinear coupling; see also \S\ref{s:res}.


\textbf{Example: Order preservation and marginal stability.}
Another setting where marginal stability can be rigorously verified are monotone systems such as the Lotka-Volterra competition system
\begin{equation}\label{eq:LV}
    \begin{split}
        u_t &= u_{xx}+u(1-u-av) \\
        v_t &= dv_{xx}+rv(1-bu-v).
    \end{split}
\end{equation}
When $0<a<1<b$ the homogeneous steady state $(u,v)=(1,0)$ is unstable.  Traveling front solutions invading this unstable state have been studied extensively and there exist parameter regimes for which these fronts are pulled~\cite{lewisliweinberger,weinbergerlewisli}, while for other parameter values the invasion fronts are pushed~\cite{HS2012,hosono03,huang10}.  In either case, the system (\ref{eq:LV}) preserves the ordering
\begin{equation} (u_1,v_1) \succeq (u_2,v_2) \ \iff u_1(x)\geq u_2(x) \quad \text{and} \quad  v_1(x)\leq v_2(x) \quad \text{for all} \ x\in\mathbb{R}. \label{eq:LVordering} \end{equation}
This ordering can be used to exclude unstable eigenvalues in a manner analogous to scalar equations obeying a comparison principle.  Intuitively this is achieved by comparison to the non-zero translational eigenfunction which precludes exponential growth of any other eigenfunction that can be scaled to satisfy the ordering (\ref{eq:LVordering}).  In practice, the details are slightly more complicated and we refer the interested reader to~\cite{bates06,FayeHolzerLotKaVolterra,leung11}.


%

\section{Robustness of invasion processes}\label{s:robust}
We study  structural stability of invasion, focusing on robustness of spreading speeds, pushed and pulled fronts, and the pushed-to-pulled transition under small perturbations.
We introduce our approach based on farfield-core decompositions and then extend to systems and singularly perturbed equations. We rely on Fredholm properties for linearized operators. Roughly speaking, a linear operator is Fredholm if it is invertible after factoring out a finite dimensional kernel and/or co-kernel. Fredholm properties of traveling waves can be computed by considering their end states, only; see \cite{SandstedeReview, FiedlerScheel, KapitulaPromislow} for a review.

We stress throughout that the perspective on wave speed selection emphasizing marginal stability is advantageous as it reduces the selection problem for reaction-diffusion PDEs to the problem of existence and stability of traveling fronts which can be treated using ODE methods.  This is particularly convenient for the study of singularly perturbed systems where there is a well established, generally applicable literature on reduction methods at our disposal.

\subsection{Robustness and bifurcations in scalar reaction-diffusion}\label{s:regularpert}
\textcolor{black}{We start this discussion with an analysis of bistable fronts, that is, fronts connecting two stable states, rather than the fronts relevant so far where the state ahead of the front is unstable. This situation is simpler yet similar to the case of pushed fronts and will be adapted to study pulled fronts and the transition between pulled and pushed fronts. }

\textbf{Motivation: bistable fronts.} We consider a bistable front in the scalar reaction diffusion equation
\begin{equation} \label{eq:mainrobust} u_t=u_{xx}+f(u;\delta), \quad f(0;\delta)=0=f(1;\delta), \quad f_u(0;\delta)<0, \quad f_u(1;\delta)<0, \end{equation}
where both asymptotic states are stable and rescaled to one and zero, respectively, for any $\delta$ sufficiently small. The simplest example would be the Nagumo equation with $f(u)=u(1-u)(u+a+\delta)$ where now $-1<a<0$! We suppose we know  existence of a monotone traveling front solution $q_0(x)$ propagating with speed $c_0$ when $\delta=0$ and wish to show that this front survives small changes in the parameter $\delta$. This problem is classical and can be attacked from several different approaches. We emphasize an approach using the implicit function theorem in preparation for the subsequent study of pushed and pulled fronts.

When (\ref{eq:mainrobust}) is recast as a system of first order equations,
\begin{eqnarray}
u'&=&v \nonumber \\
v'&=& -cv-f(u;\delta), \label{eq:robustTWeqns}
\end{eqnarray}
the traveling front is identified as a heteroclinic orbit in the intersection of the one-dimensional unstable manifold of the rest state $(1,0)$ and the one dimensional stable manifold of the rest state $(0,0)$.  As such this heteroclinic orbit lacks structural stability and is expected to be destroyed under small perturbations. To retain the front solution for $\delta\neq 0$, we must vary another parameter which in this case is simply the wave speed $c$. We argue using the implicit function theorem applied to the nonlinear equation
\[
F(q,c;\delta):=q''+cq'+f(q;\delta)\stackrel{!}{=}0, \qquad \text{with root }\quad  F(q_0,c_0;0)=0.
\]
We wish to work in spaces of localized functions and therefore set  $q=\chi_- + w$,  where $\chi_-$ is a smooth cut-off function, $\chi(\xi)=1$ for $\xi<-1$ and $\chi(\xi)=0$ for $\xi>1$, and where the ``core function'' $w$ is localized in the sense that $w\in L^2(\mathbb{R})$, in fact exponentially localized.    One now checks that $\tilde{F}(w;c):=F(\chi q+w,c;\delta)$ maps  $\tilde{F}:H^2\times\R^2\to L^2$ locally, in fact smoothly. To apply the implicit function theorem directly we would need the linearization $\partial_wF(\chi_-+w_0,c_0;\delta)=\partial_qF(q_0,c_0;\delta)$ to be bounded invertible.  However this fails due to the presence of a kernel induced by translational invariance which gives a family of solutions $q_0(\xi+\xi_0)$ and and associated translational eigenfunction $q_0'$ by differentiation with respect to $\xi_0$. We remove this translational invariance by  augmenting the operator $F$ with a phase condition and consider
\begin{equation} G(w,c;\delta) =\left(\begin{array}{cc} q''+cq'+f(q;\delta) \\ \langle q-q_0,q_0'\rangle \end{array}\right),\qquad  \partial_{(w;c)} G(w_0,c_0;0)=\left(\begin{array}{cc} B_0 & q_0' \\ \qquad \langle \cdot, q_0'\rangle & 0 \end{array}\right).
\label{eq:Gdefpersistence} \end{equation}
where we understand $q=u_-(\delta)\chi_-+w$, and wrote $B_0= \partial_\xi^2+c_0\partial_\xi+f_u(q_0;0)$.  Note that $B_0$ is a Fredholm operator with index zero: we have $\mathrm{ker}(B_0)=\mathrm{span}\{q_0'\}$ and $\mathrm{coker}(B_0)=\mathrm{span}\{ \rme^{c_0\cdot}q_0'\}$. Then $\partial_{(w,c)} G(w_0,c_0;0)$, the joint linearization of $G(w,c;0)$ with respect to $w$ and $c$, is also Fredholm of index 0, essentially by virtue of simultaneously adding one condition and one degree of freedom.  The kernel is  trivial since $q_0'$ does not lie in the range of $B_0$, that is, $\lambda=0$ is algebraically simple, so that $\partial_{(w;c)}G(w_0,c_0;0)$ is  bounded invertible.  This implies that the front solutions can be continued in $\delta$ as a family of front solutions $q(\delta)$ with speeds $c(\delta)$ for $|\delta|\ll 1$.
Monotonicity of the front is preserved and therefore Sturm-Liouville implies that the front will remain stable.  Furthermore, expansions of the wave speed in the parameter $\delta$ can be obtained using a Lyapunov-Schmidt reduction and projecting onto the co-range one obtains,
\[ c(\delta)=c_0-\frac{\langle f_\delta(q_0;0),\psi\rangle}{\langle q_0',\psi\rangle}\delta+\mathcal{O}(\delta^2), \quad \text{ where }\psi= \rme^{c_0\xi}q_0'(\xi) \text{ spans the cokernel. }\]

\textbf{Pushed fronts.}  The persistence result for bistable fronts can be adapted in a somewhat straightforward fashion to pushed fronts in scalar reaction-diffusion equations.  \textcolor{black}{With the set-up as in (\ref{eq:mainrobust}) we now assume that the zero state is unstable, $f_u(0;0)>0$, considering for instance again the Nagumo equation as the simplest example with $f(u)=u(1-u)(u+a+\delta)$ and $\frac{1}{2}>a>0$.}  The pushed  front that exists at $\delta=0$ is assumed to have steep exponential decay in the sense that
\begin{equation} q_0(\xi)\sim \rme^{\frac{1}{2}\left(-c_0-\sqrt{c_0^2-4f_u(0;0)}\right) \xi} \quad \text{as} \ \  \xi\to\infty. \label{eq:decayforpushedscalar}\end{equation}
Repeating the approach for bistable fronts, we set $q=\chi_-+w$ and seek solutions of  $G(w,c;\delta)=0$.  An important difference arises here, however, as the linear operator $B_0$ is no longer Fredholm index zero. Indeed, it retains a one-dimensional kernel spanned by the translational eigenfunction $q_0'$, but the cokernel is now trivial since the candidate adjoint eigenfunction $\psi= \rme^{c_0\xi}q_0'(\xi)$ grows exponentially as $\xi\to\infty$ and is thus not square integrable. Therefore, the Fredholm index of $B_0$ is one. We could now simply solve $G(w,c_0;0)=0$ without varying $c$, which would however give solutions without the steep decay that we require for pushed fronts. This is also clear from the planar phase portrait, where the front corresponds to a saddle-sink heteroclinic which is of course robust.

The resolution is to work in an exponentially weighted space to exclude solutions with weak exponential decay, for instance $L^2_{\mathrm{exp}, 0,\frac{c_0}{2}}(\mathbb{R})$; see~\eqref{e:Lpweighted}.  The $L^2$-adjoint then acts on spaces with the inverse weight,
\[
B_0:H^2_{\mathrm{exp}, 0,\frac{c_0}{2}}(\mathbb{R})\to L^2_{\mathrm{exp}, 0,\frac{c_0}{2}}(\mathbb{R}),\qquad\qquad B_0^*:H^2_{\mathrm{exp}, 0, -\frac{c_0}{2}}(\mathbb{R})\to L^2_{\mathrm{exp}, 0, -\frac{c_0}{2}}(\mathbb{R}).
\]
In particular, the domain of the adjoint now contains $\psi$ and $B_0$ is Fredholm of index zero, with
 $\mathrm{ker}(B_0) =\mathrm{span}\{q_0'\}, \quad \mathrm{coker}(B_0)=\mathrm{span} \left\{\rme^{c_0\cdot}q_0'\right\}. $

We can now proceed as above, noting that $q_0'$ does not belong to the range of $B_0$, and solve $G(w,c;\delta)=0$ with the implicit function theorem. Monotonicity combined with the steep exponential decay of the front imply that the perturbed fronts remain marginally stable due to an isolated zero eigenvalue with translational eigenfunction $q_\delta'(\xi)$.  Altogether -- in the scalar reaction-diffusion equation setting -- marginally stable pushed fronts are robust to small changes in system parameters!

\textbf{Pulled fronts.} We now consider pulled fronts, where the analysis deviates somewhat from the template followed thus far. Suppose that for $\delta=0$ there exists a traveling front $q_0(\xi)$ propagating with the linear spreading speed $c_0=2\sqrt{f_u(0;0)}$. Further assume that this front is monotone with asymptotics
\begin{equation} q_0(\xi)\sim \xi \rme^{-\frac{c_0}{2}\xi}, \quad \xi\to\infty. \label{eq:decaypulledpersistence} \end{equation}
We establish persistence of this marginally stable pulled front for small $\delta\neq0$.

The linearization at the origin has a Jordan block which leads to the algebraic prefactor in the asymptotics. In contrast to the bistable and pushed case, the wave speed is not a free parameter but needs to be fixed as the linear spreading speed  $c_\mathrm{lin}(\delta)=2\sqrt{f_u(0;\delta)}$, which also preserves the Jordan block structure in the traveling-wave equation. Nevertheless, as in the case of pushed fronts, the heteroclinic orbit in the traveling wave equation~\eqref{eq:robustTWeqns}  is structurally stable as a saddle-sink connection upon varying $\delta$ and $c=c_\delta$.

We can now consider $G(w,c_\mathrm{lin}(\delta);\delta)=0$ and solve for $w$, only, with the implicit function theorem.
In the following, we design however a functional-analytic approach to persistence that extracts more precise information on the asymptotics in the leading edge. This turns out to be essential when studying stability and identifying a transition to pushed fronts.

These refined asymptotics are obtained, following~\cite[Theorem 2]{as1}, with a modified ansatz for $q$, namely
\begin{equation} q=\chi_-(\xi)+w(\xi)+\chi_+(\xi)(a\xi+b)\rme^{-\frac{c_\mathrm{lin}(\delta)}{2}\xi}, \label{eq:qdef} \end{equation}
where $\chi_\pm$ are smooth cut-off functions on $\mathbb{R}^{\pm}$, and $a,b$ are real parameters. The ``core'' function $w\in L^2_{\mathrm{exp}, 0, \eta_0+\varepsilon}(\mathbb{R})$ is localized in the sense that it has stronger decay than the front itself.  We then wish to solve
\begin{equation} G(w,a,b;\delta) =\left(\begin{array}{cc} q''+c_\mathrm{lin}(\delta)q'+f(q;\delta) \\ \langle q-q_0,q_0'\rangle \end{array}\right), \label{eq:Gpulldefpersistence} \end{equation}
where we again use the shorthand $q$ for the sum in~\eqref{eq:qdef}.  We first verify that $G:H^2_{\mathrm{exp}, 0,\eta_0+\varepsilon}\times\mathbb{R}^2\times\R\to L^2_{\mathrm{exp}, 0, \eta_0+\varepsilon}\times \mathbb{R}$. Since the ansatz for $q$ involves far-field terms that decay with rate $\rme^{-\eta_0 \xi}$, the fact that the map $G(w,a,b;\delta)$ preserves the strong localization of the core function $w$ is not immediately apparent from the definition of the operator. It is rather implied by the fact that the far field terms are chosen to be exact solutions of the asymptotic system. This can be seen explicitly setting $\psi(\xi)=\chi_+(\xi)(a\xi+b)\rme^{-\frac{c_\mathrm{lin}(\delta)}{2}\xi}$ and focusing on $\xi>1$ where $\chi_+(\xi)=1$ and $\chi_-(\xi)=0$.  Then
\[  q''+c_\mathrm{lin}(\delta)q+f(q;\delta)=w''+c_\mathrm{lin}(\delta)w'+f(w+\psi;\delta)-f_u(0;\delta)\psi +\left(\psi''+c_\mathrm{lin}(\delta)\psi'+f_u(0;\delta)\psi\right). \]
Since $\psi$ solves the linear equation in  $\xi>1$, the terms in the parenthesis vanish. Since also  $f(w+\psi;\delta)-f_u(0;\delta)\psi\in L^2_{\mathrm{exp}, 0, \eta_0+\varepsilon}$, we find the desired localization of $G$. It is also not difficult to show smoothness of $G$, so that we may try to  apply the implicit function theorem. Since $G(w_0,a_0,b_0;0)=0$ for some $a_0>0$ and appropriate $w_0,b_0$,  we need to examine the linearization
\[
\partial_{(w,a,b)} G(w_0,a_0,b_0;0)=\left(\begin{array}{ccc} B_0 & B_0 \left[\xi \rme^{-\frac{c_0}{2}\xi}\chi_+(\xi)\right] &  B_0 \left[\rme^{-\frac{c_0}{2}\xi}\chi_+(\xi)\right]  \\ \langle \cdot, q_0'\rangle & \langle \chi_+(\xi)\xi\rme^{-\frac{c_0}{2}\xi} , q_0'\rangle &  \langle \chi_+(\xi)\rme^{-\frac{c_0}{2}\xi} , q_0'\rangle \end{array}\right).
\]
Strong localization eliminates one direction in the asymptotics, thus adds a condition and thereby reduces the Fredholm index of $B_0$ to $-1$. Adding two parameters and one condition implies by Fredholm bordering theory that $\partial_{(w,a,b)}G(w_0,a_0,b_0;0)$ is Fredholm with index zero.  It remains to verify that the kernel of this operator is trivial.  While this may appear more intricate, it follows  exactly as in the bistable and pushed cases: the only bounded solution of $B_0$ is given by $q_0'$, but this potential kernel element is not compatible with the linearization of the phase condition which requires  $\beta \langle q_0',q_0'\rangle=0$ forcing triviality of the kernel.   The implicit function theorem then implies persistence of the pulled front.  Marginal stability is obtained from monotonicity owing to the fact that the parameter $a(\delta)>0$ for all $\delta$ sufficiently small.

\textbf{The pushed-to-pulled transition.}
The most interesting situation at this point is the boundary between pushed and pulled propagation. Reinspecting persistence of pulled fronts, a monotone pulled front at $\delta=0$ persists for $\delta$ sufficiently small. Monotonicity in the leading edge was encoded in the fact that the parameter $a$ measuring leading edge decay $a \xi \rme^{-c_\mathrm{lin}(\delta)/2)\xi}$  was positive when $\delta=0$.  The approach outlined in the pulled section works without modification in the scenario where the front instead has pure exponential decay in the sense that $a=0$ when $\delta=0$, yielding a family of traveling fronts parameterized by $\delta$. Generically, $a'(0)\neq 0$ and one obtains that fronts propagating with the linear speed lose monotonicity whenever $a'(0)\delta<0$. Fix $a'(0)>0$.
Then for $\delta>0$ and sufficiently small the fronts are pulled, but for $\delta<0$ the front with speed $c_\mathrm{lin}(\delta)$ is non-monotone, hence unstable and not selected.  So what happens for $\delta<0$? It turns out that in this case a marginally stable pushed front bifurcates. We analyze this bifurcation with some small refinements to the method we used for robustness of pushed fronts undertaken previously.

Recall that $c_{\mathrm{lin}}(\delta)=2\sqrt{f_u(0;\delta})$ is the linear spreading speed in the unperturbed case.
The bifurcating pushed fronts we seek will have speed  $c=c_{\mathrm{lin}}(\delta)+\sigma^2$ where we now posit the following decomposition
\begin{equation} q(\xi)=\chi_-(\xi)+w(\xi)+\chi_+(\xi)\beta \rme^{\nu(\sigma,\delta)\xi}, \quad \text{with }  \nu(\sigma,\delta)=-\frac{1}{2}\left(c_{\mathrm{lin}}(\delta)-{\sigma^2} -\sqrt{(c_{\mathrm{lin}}(\delta)+\sigma^2)^2-4f_u(0;\delta)}\right).\label{eq:pushpullansatz} \end{equation}
In analogy to~\eqref{eq:Gpulldefpersistence} this leads to finding roots of  $G:H^2_{\mathrm{exp}, 0, \eta_0+\varepsilon}\times\mathbb{R}^2\times \mathbb{R}\to L^2_{\mathrm{exp}, 0, \eta_0+\varepsilon}\times \mathbb{R}$ for some $\varepsilon>0$,
\begin{equation} G(w,\beta,\sigma) =\left(\begin{array}{cc} q''+(c_{\mathrm{lin}}(\delta)+\sigma^2)q'+f(q;\delta) \\ \langle q-q_0,q_0'\rangle \end{array}\right), \label{eq:Gtransitiondefpersistence} \end{equation}
with linearization again Fredholm of index zero,
\[
\partial_{(w,\beta,\sigma)} G(w_0,\beta_0,0;0)=\left(\begin{array}{ccc} B_0 & B_0 \left[\rme^{-\frac{c_0}{2}\xi}\chi_+(\xi)\right] &  B_0 \left[\xi\frac{\partial \nu}{\partial \sigma}(0,0) \rme^{-\frac{c_0}{2}\xi}\chi_+(\xi)\right]  \\ \langle \cdot, q_0'\rangle & \langle \chi_+(\xi)\rme^{-\frac{c_0}{2}\xi} , q_0'\rangle &  \langle \xi\frac{\partial \nu}{\partial \sigma}(0,0)\chi_+(\xi)\rme^{-\frac{c_0}{2}\xi} , q_0'\rangle \end{array}\right).
\]
The kernel is  trivial as before since, (i) the only bounded solution of $B_0$ is $q_0'$, and (ii) the potential kernel element is not compatible with the linearization of the phase condition; see~\cite{pp} for details. The implicit function theorem then generates a family of monotone pushed fronts for $\delta \lesssim 0$ with speed
$c_{\mathrm{ps}}(\delta)=c_{\mathrm{lin}}(\delta)+\sigma'(0)^2\delta^2 +\mathcal{O}(\delta^3). $
In this scalar setting, the imposition of strong exponential decay of the front as encoded in (\ref{eq:pushpullansatz}) together with monotonicity imply marginal stability and hence that a pushed, selected front has bifurcated.

\textbf{Stability and extensions to higher order parabolic equations and systems.}
Extensions of these methods to systems of reaction-diffusion equations or, more generally, higher order parabolic equations are possible adapting the methods described here see~\cite{as1}.  A crucial difference is that stability is no longer inherently intertwined with monotonicity and so marginal stability of the constructed fronts must be established directly in those settings.  We illustrate the idea, first developed in \cite{PoganScheel}, once again restricting to the scalar setting and to robustness of marginal stability of pulled fronts.

Assume the existence of a front with asymptotics (\ref{eq:decaypulledpersistence}).  We study the linear eigenvalue problem and seek eigenfunctions of the operator $B_0-\lambda I$, bounded in the weighted space $L^2_{\mathrm{exp},0,\eta_\delta}$, where $\eta_\delta=c_{\mathrm{lin}}(\delta)/2$ and recall that
$B_0=\partial^2_\xi+c_{\mathrm{lin}}(\delta)\partial_\xi +f_u(q_\delta;\delta).$
The essential spectrum of $B_0$ in the weighted space $L^2_{\mathrm{exp},0,\eta_\delta}$ touches the imaginary axis and so our principal concern is that an eigenvalue emerges from the essential spectrum upon perturbation.  We first record the exact solution at $\xi=+\infty$
\[ u(\xi)=C\rme^{\left(-\frac{c_{\mathrm{lin}}(\delta)}{2}-\sqrt{\lambda}\right)\xi},\qquad \qquad u''+c_\mathrm{lin}(\delta)u' +f_u(0;\delta)u-\lambda u=0,\]
for any $\mathrm{Re}(\lambda)\geq 0$.  We again apply a far-field core ansatz.  Letting $\gamma^2=\lambda$ we then seek bounded solutions of the form
\[ u(\xi,\delta,\gamma)=w(\xi)+\beta \chi_+(\xi) \rme^{-(\eta_\delta+\gamma) \xi},\]
where $w\in L^2_{\mathrm{exp},0,\eta_\delta+\varepsilon}$.  Recall that $B_0:H^2_{\mathrm{exp},0,\eta_\delta+\varepsilon} \to L^2_{\mathrm{exp},0,\eta_\delta+\varepsilon}$ is Fredholm with index $-1$ with $\mathrm{coker}(B_0)=\mathrm{span}\{\psi\}$.  Writing $P$ for the projection of $B_0$ onto its range, one can show that eigenvalues correspond to values of $\gamma=\sqrt{\lambda}$ for which there is a non-trivial solution to
\begin{equation}
    \begin{split}
    P\left( (B_0-\gamma^2)(w+\beta \chi_+ \rme^{(\eta_\delta-\gamma)\cdot}\right) &=0 \\
    \langle (B_0-\gamma^2)(w+\beta \chi_+ \rme^{(\eta_\delta-\gamma)\cdot},\psi \rangle &=0
    \end{split}
\end{equation}
Since the exponential in the far-field term is an exact solution of the asymptotic system at $\xi=+\infty$ and the front converges exponentially, $(B_0-\gamma^2)$ maps $\chi_+\rme^{(\eta_\delta-\gamma)\cdot}$ to the weighted space $L^2_{\mathrm{exp},0,\eta_\delta+\varepsilon}$.  Consequently,   the first of these equations can be solved via the implicit function theorem to obtain $w=\beta \bar{w}(\gamma,\delta)$.  Substituting into the second equation and factoring $\beta$, we obtain the existence of a bounded solution (and hence an eigenfunction), whenever the reduced characteristic function vanishes,
\[ E(\gamma,\delta):=\langle (B_0-\gamma^2)(\bar{w}(\gamma,\delta)+\chi_+ \rme^{(\eta_\delta-\gamma)\cdot},\psi \rangle=0. \]
Note that $E(\gamma,\delta)$ is continuous in a neighborhood of $(\gamma,\delta)=(0,0)$.  When $\gamma=\delta =0$ the only bounded solution of $B_0w=0$ is $q_0'$, which however is not an element of the weighted space $L^2_{\mathrm{exp},0,\eta_\delta+\varepsilon}$ so that $E(0,0)\neq0$. Continuity of $E(\gamma,\delta)$ then precludes the existence of an eigenvalue for $\mathrm{Re}(\lambda)\geq 0$ for $\lambda$ near the origin.
Away from the essential spectrum, regular spectral perturbation theory shows upper semicontinuity of the spectrum which gives the desired marginal stability for $\delta\sim 0$.


\subsection{Singular perturbation approaches to systems of reaction-diffusion equations}

Systems of reaction-diffusion equations pose analytical challenges for the study of front propagation for a variety of reasons.  The traveling wave equations are higher dimensional rendering phase space methods more complicated to employ while some powerful PDE tools (comparison principles, variational methods, Sturm-Liouville theory) may not apply.  More systematic insight into systems is gained in singular limits where an at least formal reductions may be possible. Making such reductions rigorous at the level of the PDE can be challenging if not impossible as one seeks to reduce for \emph{all} solutions.  We sketch here, how to use geometric singular perturbation theory~\cite{fenichel}  to reduce ODEs for traveling waves and their stability. This affords a systematic and rigorous method to study pushed and pulled invasion fronts in singularly perturbed reaction-diffusion settings.  The results on nonlinear marginal stability from \S\ref{s:nlmsp} then imply that these fronts are selected in the sense that they attract open sets including compactly supported initial data in the reaction-diffusion PDEs.  The construction of traveling fronts using geometric singular perturbation theory has a long history; see for example~\cite{RottschaferWayne} and references therein. We emphasize however that the approach to stability outlined here -- in particular for pulled fronts -- constitutes a novel and convenient method for establishing the requisite properties of the traveling front required for selection; see~\cite{averyholzerscheelkellersegel} for  details of the method in a specific problem.

\textbf{Systems of equations with one fast-reaction.}
As a motivating class of examples, consider
\begin{equation}
    \begin{split}
        u_t &= u_{xx}+\frac{1}{\delta}f(u,v) \\
        v_t &= v_{xx}+g(u,v), \label{eq:GSPTgen}
    \end{split}
\end{equation}
where $0< \delta \ll 1$ is a small parameter \textcolor{black}{and we set diffusivities equal to 1 for notational ease.}  We assume that $(0,0)$ is an unstable homogeneous state while there also exists a  stable coexistence steady state $(u^*,v^*)$.  Fronts solve  the system of ODEs
\begin{equation}
\begin{displaystyle}
    \begin{array}{rlrl}
     \sqrt{\delta}\dot{u_1} &= u_2, & \qquad \qquad \dot{v_1}&=v_2, \\
       \sqrt{\delta} \dot{u_2}& = -c\sqrt{\delta}u_2  -f(u_1,v_1),&\qquad  \qquad
        \dot{v_2}&= -cv_2- g(u_1,v_1).
    \end{array}
\end{displaystyle} \label{eq:GSPTgenTW}
\end{equation}
Setting $\delta=0$, \eqref{eq:GSPTgenTW} becomes differential-algebraic and solutions are constrained to lie on slow manifolds describing branches of the zero set of $f(u,v)=0$.  Suppose that that there exist a smooth function $\phi(v)$ such that (i) $f(\phi(v),v)=0$, (ii) $\phi(0)=0$, and (iii) $\phi(v^*)=u_*$. Then we denote by
\[ \mathcal{M}_0=\left\{ (u_1,u_2,v_1,v_2) \ | \ u_1=\phi(v_1), \ u_2=0 \right\}, \]
the associated ``slow manifold''.  On $\mathcal{M}_0$, we obtain the limiting slow subsystem
\begin{equation}\label{e:slow limit}
    \begin{split}
        \dot{v_1}&= v_2 \\
        \dot{v_2}&= -cv_2-g(\phi(v_1),v_1),
    \end{split}
\end{equation}
which is exactly the traveling wave equation associated with the scalar reaction-diffusion equation
\begin{equation} v_t= v_{xx}+g(\phi(v),v). \label{eq:GSPTscalarreducedPDE} \end{equation}
We suppose that existence and stability of fronts with speed $c$ in~\eqref{eq:GSPTscalarreducedPDE} has been established, using any of the techniques discussed previously. The question then becomes under what conditions those fronts in the singular limit yield  marginally stable fronts in the full system when $\delta\gtrsim 0$?  For existence, only, the typical requirement is normal hyperbolicity of the slow manifold $\mathcal{M}_0$;  Fenichel's persistence Theorem~\cite{fenichel} then implies that this manifold will smoothly continue as an invariant manifold for  $\delta$ sufficiently small.  The perturbed manifold has the form
\begin{equation} \mathcal{M}_\delta =\left\{ (u_1,u_2,v_1,v_2) \ | \ u_1=\phi(v_1)+\delta \psi_1(v_1,v_2,\delta), u_2=\delta\psi_2(v_1,v_2,\delta) \right\}, \label{eq:Mdelta} \end{equation}
for smooth functions $\psi_{1,2}$. The manifold $\mathcal{M}_\delta$ includes all locally bounded solutions including the steady states and the heteroclinic orbit connecting them. Restricting~\eqref{eq:GSPTgenTW} to $\mathcal{M}_\delta$ then  reduces the slow existence problem to a regular perturbation problem as treated in \S\ref{s:regularpert}.

Marginal spectral stability can be studied in a similar fashion and we outline here the ideas, developed in~\cite{averyholzerscheelkellersegel} for a specific example.  Writing the  front propagating with  selected spreading speed $c_\mathrm{sel}$ (pulled or pushed) as $(U_1(\xi),U_2(\xi),V_1(\xi),V_2(\xi))$, denote perturbations as
\[ u_1=U_1(\xi)+p_1, \ u_2=U_2(\xi)+p_2, \ v_1=V_1(\xi)+q_1, \ v_2=V_2(\xi)+q_2, \]
and write the combined existence and linear stability problem as
\begin{align}
     \sqrt{\delta}\dot{U_1} &= U_2, &  \sqrt{\delta}\dot{p_1} &=p_2, \nonumber \\
     \sqrt{\delta}\dot{U_2}&=  -c_{\mathrm{sel}}\sqrt{\delta} U_2 -f(U_1,V_1) , & \sqrt{\delta}\dot{p_2}&=
     -c_{\mathrm{sel}}\sqrt{\delta} p_2 -f_u(U_1,V_1)p_1-f_v(U_1,V_1)q_1-\delta \lambda p_1, \nonumber \\
     \dot{V_1} &= V_2, &  \dot{q_1} &=q_2, \nonumber \\
     \dot{V_2}&=  -c_{\mathrm{sel}} V_2 -g(U_1,V_1) , & \dot{q_2}&=
     -c_{\mathrm{sel}} q_2 -g_u(U_1,V_1)p_1-g_v(U_1,V_1)q_1-\lambda q_1.
     \end{align}
One can again perform a singular perturbation reduction to this system obtaining the existence of a four dimensional slow manifold expressed as a graph over the slow variables $(V_1,V_2,q_1,q_2)$:
\begin{equation}
\begin{displaystyle}
    \begin{array}{rlrl}
 U_1&=\phi(V_1)+\delta \psi_1(V_1,V_2,\delta),&\qquad \quad U_2&=\delta\psi_2(V_1,V_2,\delta)  \\
  p_1&=\phi'(V_1)q_1+\delta\tilde{\psi}_1(V_1,V_2,\delta,\lambda)(q_1,q_2)^T,&\qquad \quad
p_2&=\delta\tilde{\psi}_2(V_1,V_2,\delta,\lambda)(q_1,q_2)^T,
    \end{array}
\end{displaystyle}
\end{equation}
where $\psi_{1,2}$ are as in (\ref{eq:Mdelta}), $\tilde{\psi}_{1,2}$ are linear maps and we restrict to a priori bounded $|\lambda|<\Lambda$; see~\cite[Prop 2.4]{averyholzerscheelkellersegel}.  Fenichel's theorem allows the reduction of the stability problem to a system of equations
\begin{align}
    \dot{q_1} &=q_2 \nonumber \\
    \dot{q_2}&=
    -c_{\mathrm{sel}} q_2 -\left[g_u(\phi(V_1),V_1)\phi'(V_1)+g_v(\phi(V_1),V_1)\right]q_1-\lambda q_1+\mathcal{O}(\delta),
\end{align}
which is, to leading order in $\delta$, the linearization of the reduced scalar PDE (\ref{eq:GSPTscalarreducedPDE}) near the singular front.  The only possible bounded solutions lie in this slow manifold and the stability problem is therefore reduced to a regular perturbation problem.
Finally, one needs to examine the region $|\lambda|>\Lambda$ separately, using some version of a priori bounds on spectrum; we refer to Appendix B of~\cite{averyholzerscheelkellersegel} for an example.

\textbf{Example: Model of tumor invasion.}
In the  system of reaction-diffusion equations,
\begin{equation}\label{eq:tumor}
    \begin{split}
        u_t &= u_{xx}+\frac{1}{\delta}\left((1-u-v)(u+v)-\alpha u\right) +(1-u-v)v \\
        v_t &= v_{xx}+v(1-u-v).
    \end{split}
\end{equation}
$u(t,x)$ represents a tumor cell density and $v(t,x)$ cancer stems cell density; see~\cite{avery2023growth,shyntar2022tumor} for modeling and analysis which we recount here as an illustration of the general reduction methodology introduced previously. When $\alpha>1$ and $\delta\ll 1$, the zero state is unstable (with linear spreading speed $2$) and the tumor free state $(0,1)$ is stable.  A slow manifold for the traveling wave equations is given through
\[ \phi(v_1)=-\frac{\alpha-1+2v_1}{2}+\frac{1}{2}\sqrt{(1-\alpha-2v_1)^2 -4v_1^2+4v_1}, \qquad  \phi(0)=0,  \ \phi(1)=0. \]
This leads  to a reduced scalar equation with pulled fronts and spreading speed $c_{\mathrm{lin}}=2$ of the form
\[ v_t=v_{xx}+(1-\phi(v)-v)v;\]
see~\cite{shyntar2022tumor}. From the perspective taken here, these fronts then perturb as pulled, marginally stable fronts for $\delta>0$, sufficiently small; see~\cite{avery2023growth} for details.

\textbf{Example: The Lotka-Volterra model revisited.}
We return to the population model  (\ref{eq:LV}):
\begin{equation}\label{eq:LVagain}
    \begin{split}
        u_t &= u_{xx}+u(1-u-av) \\
        v_t &= dv_{xx}+rv(1-bu-v).
    \end{split}
\end{equation}
where $0<a<1<b$ so that $(u,v)=(1,0)$ is an unstable state. Invasion dynamics across the $(d,r,a,b)$ parameter space have not been completely characterized, although there is significant progress in some regimes; see for example~\cite{HS2012,hosono98,hosono03,huang10,lewisliweinberger}.  Motivated by~\cite{hosono98}, we consider the asymptotic limit $r\to\infty$.  Setting $r=\frac{1}{\delta}$, (\ref{eq:LVagain}) fits the framework of (\ref{eq:GSPTgen}) with interchanged $u$ and $v$ species. We obtain the slow manifold as a graph over the $u$ component and a reduced scalar equation for $u$,
\begin{equation} \label{eq:LVreducsedPDE}  u_t=u_{xx}+ u\left(1-u-a\phi(u)\right),\qquad \phi(u)=\left\{\begin{array}{cc} 0 & u\geq \frac{1}{b}, \\ 1-bu & u< \frac{1}{b}. \end{array}\right.  \end{equation}
This reduced equation was analyzed  in~\cite{hosono98}, with  pulled fronts when $a=0$ and  pushed fronts for $a$ near one. The methodology presented above could then be applied in the current context with one caveat: the slow manifold loses normal hyperbolicity at $v=\frac{1}{b}$ and so  additional work would be required. We suspect that dynamical systems techniques available for this slow passage through a transcritical bifurcation~\cite{KRUPASZMOLYANTRANS} could be used to determine the location of the pushed-to-pulled transition as $r\to\infty$ rigorously; see~\cite{hosono98} for a more detailed discussion and conjectures.

    We conclude this discussion with the observation that the procedure outlined above will not always be successful. Extensions would need to for instance treat situations  where the selected front is constructed from concatenations of  several slow manifolds and fast jumps between them; see~\cite{HS2012} for an example, again in the context of the Lotka-Volterra equations.

\subsection{Practical considerations:  Efficient numerical continuation}\label{s:4.3}
The perspective presented here leads to novel numerical algorithms, substantially improving on the approach presented in \S\ref{s:3.3}.

\emph{Bistable fronts} can be computed efficiently by truncating to a bounded interval, adding suitable boundary conditions (in simple cases Dirichlet boundary conditions will work), and adding a phase condition, thus finding roots of a discretized version of
\begin{equation}
    G_L[u,c]=0\quad \Longleftrightarrow \quad \displaystyle{\left\{\begin{array}{ll}
    u''(\xi)+c u'(\xi)+f(u(\xi))=0,& \xi\in(-L,L),\\
    u(-L)-u_-=0,\\
    u(L)=0,\\
    \int_{-L}^L u(\xi)u'_\mathrm{ref}(\xi)\,\rmd \xi=0.
    \end{array}
    \right.}\label{e:bistbvp}
\end{equation}
Here, $u'_\mathrm{ref}$ mimics the derivative of the front via some reference template, which in numerical continuation could be the derivative of the solution at a previous parameter value. One finds that $(u,c)$ converge exponentially, with error $\rmO(\rme^{-\eta L})$ for some $\eta>0$ depending on spectral gaps in the linearization.

\emph{Pushed fronts} can be analyzed with the exact same boundary-value problem. However, numerical performance is improved if instead of finding $u$, one substitutes $\omega(\xi)u(\xi)=:w(\xi)$ with weight $\omega=1+\rme^{c/2 \xi}$, for instance. In fact, numerically continuing bistable fronts in the parameter $a$ when $f(u)=u(1-u)(u-a)$ from $a>0$ to $a\lesssim 0$, one finds that the bistable front connecting $1$ to $0$ continues as the pushed front, still  connecting $1$ to $0$. In both cases of bistable and pushed fronts, convergence is improved if the asymptotic eigenspace is used as a boundary condition, for instance $u_\xi=\nu_\mathrm{ss} u$ where $\nu_\mathrm{ss}=-c/2-\sqrt{c^2/4-f'(0)}$ gives the asymptotic exponential decay. In more general situations, choosing asymptotically correct boundary conditions also avoids spurious bifurcations related to non-transversality of the Dirichlet subspace to the spectral complement of the stable eigenspace.

\emph{Pulled fronts}, or pushed fronts in the vicinity of the pushed-to-pulled transition require a modified setup, where one explicitly substitutes the far-field ansatz and requires a core correction to be strongly localized. One therefore replaces the first component in~\eqref{e:bistbvp} by the first equation in~\eqref{eq:Gpulldefpersistence}, with $q$ given by the farfield-core decomposition. The strong localization of the core component can be accurately encoded in requiring $u$ to belong to the strong stable subspace, which can generally be computed directly. In the scalar case, this subspace is trivial and one would impose Dirichlet and Neumann boundary conditions. In more general situations, one can impose Dirichlet boundary conditions and an integral phase condition in the far-field that prohibits decay in the direction of the eigenspace associated with the front decay in the leading edge; see~\cite{pp}. Throughout, the pulled speed can be found by continuing a double-root of the dispersion relation at $\nu=0$ in $c$, thus solving
\[
d_c(0,\nu)=0,\qquad \partial_\nu d_c(0,\nu)=0,
\]
for the variables $\nu,c$ with a Newton method and numerical continuation. For larger systems, a generalized eigenvalue-problem formulation may be more robust; see~\cite{rss07,S23}.

\emph{Pushed-to-pulled transitions.} When continuing pulled fronts, one can easily add the condition $a=0$ in the linear term of the leading edge \eqref{eq:qdef} and continue in two parameters to find the transition; see \cite[\S6.2]{pp} for examples.

\section{Analysis: proof of the marginal stability conjecture}\label{s:nlmsp}

We summarize recent results on front selection via marginal stability. We start by formulating several sufficient conditions for establishing front selection for steep initial data in \S\ref{s:a1}, including definitions of pushed and pulled fronts as well as precise characterizations of selection. We phrase definitions and results for semilinear parabolic systems
\begin{align}
    u_t = \mathcal{P}(\partial_\xi)u + c u_\xi + f(u), \quad u(\xi,t) \in \mathbb{R}^N, \quad \xi \in \R, \quad t > 0, \label{e: proofs eqn}
\end{align}
where again $\mathcal{P}$ is a matrix-valued polynomial such that $\mathcal{P}(\partial_x)$ is an elliptic operator of order $2m$. We assume that $f(0) = 0$, so that $u \equiv 0$ is an equilibrium solution, and that the linearization $\mathcal{L}_0 = \mathcal{P}(\partial_\xi) + c \partial_\xi + f'(0)$ at $u =0$ in~\eqref{e: proofs eqn} is unstable in a translation-invariant norm. We sketch proofs in \S\ref{s:a2} and point to unresolved questions in \S\ref{s:a3}.

\subsection{Rigidly propagating fronts: rigorous results and conjectures}\label{s:a1}

\textbf{Pushed rigid fronts.} We formulate assumptions that encode the properties of pushed fronts: namely, they should travel faster than the linear spreading speed, and they should be contained in the strong stable manifold of the unstable state, giving rise to marginally stable point spectrum. \textcolor{black}{The marginal stability due to point spectrum is the key ingredient that we retain. }
\begin{definition}[Marginal spectral stability --- pushed fronts]\label{def: rigid pushed}
    A stationary solution  $u_\mathrm{ps}(\xi)$ to~\eqref{e: proofs eqn} with $c=c_\mathrm{ps}$ is a \emph{pushed front} if it is marginally spectrally stable in the following sense:
    \begin{enumerate}
        \item \emph{Selected state in the wake:} There is an exponentially stable equilibrium $u_- \in \mathbb{R}^N$ to~\eqref{e: proofs eqn}.
        \item \emph{Exponential asymptotics:} The front $u_\mathrm{ps}(\xi)$ has the asymptotics
        \begin{align*}
            u_\mathrm{ps}(\xi) = u_- + \mathrm{O}(\rme^{\eta \xi}), \quad \xi \to -\infty, \qquad u_\mathrm{ps}(\xi) = u_+^0 \rme^{-\eta_\mathrm{ps} \xi} + \mathrm{O}(\rme^{-(\eta_\mathrm{ps} + \eta) \xi}), \quad \xi \to \infty,
        \end{align*}
        for some $u_+^0 \in \R^N$ and {\color{black}$\eta_\mathrm{ps}, \eta > 0$}.
        \item \emph{Minimal marginal spectral stability:} The linearization $\mathcal{L}_\mathrm{ps} = \mathcal{P}(\partial_\xi) + c_\mathrm{ps} \partial_\xi + f'(u_\mathrm{ps}(\xi))$ at $u_\mathrm{ps}(\xi)$ as an operator
        \begin{align*}
            \mathcal{L}_\mathrm{ps} : H^{2m}_{\mathrm{exp}, 0, \eta_{\mathrm{ps}}-\tilde{\eta}} \subset L^2_{\mathrm{exp}, 0, \eta_\mathrm{ps}-\tilde{\eta}} \to L^2_{\mathrm{exp}, \tilde{\eta}, \eta_\mathrm{ps}-\tilde{\eta}},\qquad   \tilde{\eta} > 0 \text{ small, arbitrary,}
        \end{align*}
        has an algebraically simple eigenvalue at $\lambda = 0$, with eigenfunction $\partial_\xi u_\mathrm{ps}$, and its spectrum  otherwise has strictly negative real part.
    \end{enumerate}
\end{definition}

\begin{figure}
    \centering
    \includegraphics[width=.7\linewidth]{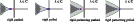}
    \caption{Typical spectra of rigid pushed and pulled invasion fronts, encoded in Defs.~\ref{def: rigid pushed} and~\ref{def: rigid pulled}, respectively, in optimal exponentially weighted spaces. Solid curves and shaded regions denote essential spectra; pink crosses denote discrete eigenvalues. Blue curves correspond to essential spectra of the unstable state $u \equiv 0$; green curves correspond to spectra of the selected state in the wake, $u_- (\xi)$.}
    \label{fig:rigid spectra}
\end{figure}
\textcolor{black}{The steep exponential decay of pulled fronts that we noticed in \S\ref{s : FKPPintro} is encoded here somewhat implicitly in spectral assumptions, where the derivative of the front contributes a zero eigenvalue; compare Fig.~\ref{f:nagumophaseportrait}. The derivative of the front always gives a pointwise solution to $\mathcal{L}_\mathrm{ps} u = 0$, but whether this solution is a genuine eigenfunction in the weighted space depends on the exponential decay rate of the front.} 
We note that this definition encodes spectral stability in a fixed exponentially weighted space, which is typical but not necessary ---  pointwise stability may hold despite unstable spectrum~\cite{FayeHolzerScheelSiemer}.
For a rigid pushed front, the only neutral spectrum is an isolated simple eigenvalue, with the essential spectrum separated by a finite spectral gap. Hence, using classical semigroup methods~\cite{Henry,Satt1977}, one obtains the following stability result.

\begin{theorem}\label{t: pushed stability}
    Suppose~\eqref{e: proofs eqn} admits a pushed front $u_\mathrm{ps}$.
    Then $u_\mathrm{ps}$ is orbitally stable with asymptotic phase against small perturbations in $L^2_{_\mathrm{exp}, 0, \eta_\mathrm{ps} - \tilde{\eta}}$.
\end{theorem}
Since $u_\mathrm{ps} (\xi) \sim \rme^{-\eta_\mathrm{ps} \xi}$ as $\xi \to \infty$, small perturbations in $L^2_{_\mathrm{exp}, 0, \eta_\mathrm{ps} - \tilde{\eta}}$ can effectively cut off the tail in the far leading edge, and so the above stability result immediately implies the following selection result.
\begin{corollary}[Selection of rigid pushed fronts]
    Assume the conditions of Thm.~\ref{t: pushed stability} hold. Then $u_\mathrm{ps}$ is a selected front in the sense of Def.~\ref{def: selected front}, with $h(t) = x_\infty(u_0) + \mathrm{O}(\rme^{-\mu t})$ for some $\mu > 0$.
\end{corollary}

\textbf{Pulled rigid fronts.} We turn our focus now to pulled fronts.
\textcolor{black}{
Pulled fronts are identified, in the simple traveling-wave equations of \S\ref{s : FKPPintro}, by the double eigenvalue at the origin and the absence of a crossing through the strong stable manifold for larger speeds. Spectrally, those assumptions translate into the presence of a pinched double root at the origin and the absence of unstable point spectrum. Our definition here relies on, similar to the case of pushed fronts, these spectral characterizations, and assumes the somewhat stronger marginal stability in a fixed exponential weight as opposed to linear pointwise characterizations. }

\begin{definition}[Marginal spectral stability --- pulled fronts]\label{def: rigid pulled}    A stationary solution  $u_\mathrm{pl}(\xi)$ to~\eqref{e: proofs eqn} with $c=c_\mathrm{lin}$ is a \emph{pulled  front} if it is marginally spectrally stable in the following sense:

\begin{enumerate}
    \item \emph{Linear spreading speed and weighted spectral stability:} Assume that the linear spreading speed $c_\mathrm{lin}$ in~\eqref{e: proofs eqn} is determined by a simple pinched double root at $(\lambda, \nu) = (0, -\eta_\mathrm{lin})$, and that the $L^2_{\exp, \eta_\mathrm{lin}}$-spectrum of the linearization $\mathcal{L}_0$ about the unstable state is contained in the open left half-plane except for a single branch touching the origin.
    \item \emph{Selected state in the wake:} There is an exponentially stable equilibrium $u_- \in \mathbb{R}^N$ to~\eqref{e: proofs eqn}.
    \item \emph{Exponential and Jordan block asymptotics:} The front $u_\mathrm{pl}(\xi)$ has the asymptotics
    \begin{align*}
        u_\mathrm{pl} (\xi) &= u_- + \mathrm{O}(\rme^{\eta \xi}), \quad \xi \to \ -\infty, \\
        u_\mathrm{pl}(\xi) &= [b(u_+^0 \xi + u_+^1) + a u_+^0] \rme^{-\eta_\mathrm{lin} \xi}, \quad \xi \to \infty
    \end{align*}
    for some $u_+^{0/1} \in \R^N$, $a, b \in \R$, $b \neq 0$, and $\eta > 0$.
    \item \emph{Minimal marginal spectral stability:} The spectrum of the linearization $\mathcal{L}_\mathrm{pl} = \mathcal{P}(\partial_\xi) + c_\mathrm{lin} \partial_\xi + f'(u_\mathrm{pl}(\xi))$ as an  operator
    \begin{align*}
        \mathcal{L}_\mathrm{pl} : H^{2m}_{\mathrm{exp}, 0, \eta_{\mathrm{lin}}} \subset L^2_{\mathrm{exp}, 0, \eta_{\mathrm{lin}}} \to L^2_{\mathrm{exp}, 0, \eta_{\mathrm{lin}}}
    \end{align*}
    has strictly negative real part asiode from the simple branch  at the origin induced by (i).
    \item \emph{No resonance at $\lambda = 0$:} There is no nontrivial solution $u \in L^\infty_{0, \eta_\mathrm{lin}}$ to $\mathcal{L}_\mathrm{pl} u = 0$.
\end{enumerate}
\end{definition}

The asymptotics captured in Def.~\ref{def: rigid pulled}(iii) are generic given a pinched double root of the dispersion relation at $(\lambda, \nu) = (0, -\eta_\mathrm{lin})$, as this double root corresponds to a length 2 Jordan block in the linearization of the traveling wave equation at the origin. The non-generic case $b = 0$ would lead to a violation of condition (v), and corresponds to a transition between pushed and pulled invasion~\cite{pp}.


Selection of pulled fronts is significantly more subtle when compared to pushed fronts since there is no exponential weight that induces a spectral gap in the linearization. Selection of rigid pulled fronts in systems without comparison principles has only recently been established~\cite{AveryScheel, avery2}.

\begin{theorem}[\cite{AveryScheel, avery2} --- selection of rigid pulled fronts]\label{t: rigid pulled unpatterned}
    Suppose $u_\mathrm{pl}$ is a pulled front in the sense of Definition~\ref{def: rigid pulled}. Then $u_\mathrm{pl}$ is a selected front in the sense of Def.~\ref{def: selected front}, with position correction $h(t) = - \frac{3}{2 \eta_\mathrm{lin}} \log t + \mathrm{O}(1)$.
\end{theorem}

The universality of the logarithmic position correction, that is, that the shift $-\frac{3}{2 \eta_\mathrm{lin}} \log t$ is independent of the steep initial condition and depends on the equation only through $\eta_\mathrm{lin}$, was predicted for pulled fronts in~\cite{EbertvanSaarloos} based on matched asymptotics, and is confirmed rigorously by Thm.~\ref{t: rigid pulled unpatterned}. 

\textbf{Rigid fronts at the pushed-to-pulled transition.} If a front $u_\mathrm{tr}$ satisfies Def.~\ref{def: rigid pulled} but with $b = 0$, and a one-dimensional solution space to $\mathcal{L}_\mathrm{pl} u = 0$ in $L^\infty_{0, \eta_\mathrm{lin}}$ spanned by $\partial_\xi u_\mathrm{tr}$, then the system is at a boundary point between pushed and pulled invasion; see \S\ref{s:4.3}. Recent work relying on comparison principles established that in scalar second order equations, these fronts at the pushed-to-pulled transition are selected in the sense of Def.~\ref{def: selected front}, with position correction $h(t) = - \frac{1}{2 \eta_\mathrm{lin}} \log t + \mathrm{O}(1)$~\cite{AHR1, AHR2, AHR3, Giletti22}. We expect that the methods of  Thm.~\ref{t: rigid pulled unpatterned} extended to  an analogous result more generally in systems. It would also be interesting to get uniform asymptotics  for nearby parameter values, where shifts change from $-\frac{1}{2\eta_\mathrm{lin}}\log(t)$ to  $-\frac{3}{2\eta_\mathrm{lin}}\log(t)$ and to $\rmO(1)$.

\begin{conjecture}[Selection at the pushed-to-pulled transition]
    Assume $u_\mathrm{tr}$ is a pushed-to-pulled transition front. Then $u_\mathrm{tr}$ is a selected front in the sense of Def.~\ref{def: selected front}, with position correction $h(t) = -\frac{1}{2 \eta_\mathrm{lin}} \log t + \mathrm{O}(1)$.
\end{conjecture}

\subsection{Rigidly propagating fronts: proof sketches}\label{s:a2}

\subsubsection{Selection of rigid pushed fronts: proof of Thm.~\ref{t: pushed stability}} Letting $u_\mathrm{ps}$ denote a rigid pushed front, one studies the dynamics of perturbed solutions $u(\xi, t) = u_\mathrm{ps}(\xi) + v(\xi, t)$. To stabilize the spectrum in the leading edge, we fix $\tilde{\eta} > 0$ small enough so that condition (iv) of Def.~\ref{def: rigid pushed} holds. We set $\eta_0 =\eta_\mathrm{ps} - \tilde{\eta}$, and define the weighted perturbation $w = \omega_{0, \eta_0} v$, which then solves the perturbed equation
\begin{align}
    w_t = \mathcal{L}^\omega_\mathrm{ps} w + \omega_{0, \eta_0} N \left( \frac{w}{\omega_{0, \eta_0}} \right), \label{e: pushed rigid perturbation equation}
\end{align}
where $\mathcal{L}^\omega_\mathrm{ps} w = \omega_{0, \eta_0} \mathcal{L}^\omega_\mathrm{ps} \left( \frac{w}{\omega_{0, \eta_0}} \right)$, and $N(v) = \mathrm{O}(v^2)$ is a quadratic nonlinear remainder. We set $\tilde{N}(w) = \omega_{0, \eta_0} N \left( \frac{w}{\omega_{0, \eta_0}} \right),$ and note that $\tilde{N}(w) = \mathrm{O}(\frac{1}{\omega_{0, \eta_0}} w^2)$ contains an exponentially localized prefactor. Solutions to~\eqref{e: pushed rigid perturbation equation} with initial data $w_0$ small may be constructed as fixed points of the variation of constants formula,
\begin{align}
    w(t) = \rme^{\mathcal{L}^\omega_\mathrm{ps} t} w_0 + \int_0^t \rme^{\mathcal{L}^\omega_\mathrm{ps} (t-s)} \tilde{N}(w(s)) \, ds. \label{e: pushed voc}
\end{align}

The spectrum of $\mathcal{L}^\omega_\mathrm{ps}$ contains an isolated eigenvalue at $\lambda = 0$, with eigenfunction $\omega_{0, \eta_0} \partial_\xi u_\mathrm{ps}$ and is otherwise stable. Using classical semigroup methods~\cite{Henry}, one finds for the linearized dynamics
\begin{align}
    \rme^{\mathcal{L}^\omega_\mathrm{ps} t} = \omega_{0, \eta_0} \partial_\xi u_\mathrm{ps} P_0 + \mathrm{O}(\rme^{-\mu t}), \label{e: pushed unpatterned linear asymptotics}
\end{align}
where $P_0$ is a projection along the adjoint eigenfunction at $\lambda = 0$. The lack of temporal decay of $\rme^{\mathcal{L}^\omega_\mathrm{ps} t}$ is an apparent obstruction to a nonlinear stability argument. One may overcome this by capturing the neutral dynamics associated to the translational mode with a \emph{temporal modulation ansatz}, $\tilde{w}(\xi, t) = \omega_{0, \eta_0}(\xi)[u(\xi + \psi(t), t) - u_\mathrm{ps} (\xi)]$. Using a fixed point argument to construct $\psi$ and $\tilde{w}$ simultaneously, one can show~\cite{Satt1977} that if $\| u_0 - u_\mathrm{ps}\|_{L^\infty_{0, \eta_0}}$ is small, there exists $\psi_\infty \in \mathbb{R}$ such that
\begin{align*}
    \left\| \omega_{0, \eta_0} [u(\cdot, t) - u_\mathrm{ps}(\cdot + \psi(t))]\right\|_{L^\infty} + | \psi(t) - \psi_\infty| \leq C\rme^{-\mu t}
\end{align*}
for some constants $C, \mu > 0$, and hence the solution indeed converges to the shifted wave $u_\mathrm{ps}(\xi + \psi_\infty)$, exponentially quickly in time; see Fig.~\ref{fig:rigid pushed convergence}, left panel. Since the initial perturbation is allowed to decay more slowly than $u_\mathrm{ps}$ itself (that is, $\eta_0 < \eta_\mathrm{ps})$, the basin of attraction of the front includes many steep initial conditions, and so the front is a selected front.

\subsubsection{Selection of rigid pulled fronts: proof of Thm.~\ref{t: rigid pulled unpatterned}} While for pushed fronts, perturbations that cut off the front tail remain localized perturbations in the  weight that stabilizes the essential spectrum, this is not the case for pulled fronts: if we conjugate with a weight $\omega_{0, \eta_\mathrm{lin}}$ that marginally stabilizes the essential spectrum, a perturbation that cuts off the front tail is large, of size $\mathrm{O}(\xi)$ as $\xi \to \infty$. As a result, selection of pulled fronts cannot be established by a direct stability argument viewing the front $u_\mathrm{pl} (\xi)$ as an equilibrium in the co-moving frame.

The asymptotic matching analysis of~\cite{vanSaarloosReview} suggests that after the front tail is cut off by a perturbation, it rebuilds itself according to the linearized dynamics near the unstable state $u \equiv 0$. The solution therefore develops a Gaussian tail, which is matched with the front profile on an intermediate length scale, thereby  determining the front location and the typical logarithmic delay.

The proof of Thm.~\ref{t: rigid pulled unpatterned} starts with this asymptotic matching procedure to construct an approximate solution describing the transition from a front profile to a diffusive tail. A nonlinear stability argument for this approximate solution establishes that nearby actual solutions also exhibit these dynamics, including solutions with steep initial data.

\begin{figure}
    \centering
    \includegraphics[width=1\linewidth]{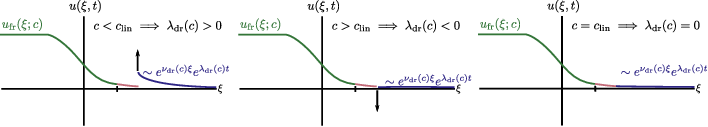}
    \caption{Schematics of attempts to construct an approximate solution which governs the propagation dynamics, matching a front with speed $c$ in the wake (green) to a tail governed by the linearized dynamics in the leading edge (navy blue), in the frame moving with the front speed, $\xi = x - ct$. If $c < c_\mathrm{lin}$ (left), the tail grows exponentially in time, and so cannot be matched with the stationary front; if $c > c_\mathrm{lin}$ (center), the tail decays exponentially in time, and so still cannot be matched with the front. If $c = c_\mathrm{lin}$, the tail is neither exponentially growing nor decaying in time, and so it is possible to match it with the front in the wake.}
    \label{fig:pulled matching}
\end{figure}

\textbf{Why marginal stability?} To execute this strategy, one determines the leading order tail behavior from the linearized equation near $u \equiv 0$, and then matches this on an intermediate length scale with the front profile. The linearized equation in the leading edge is $u_t = \mathcal{P}(\partial_\xi)u + c \partial_\xi u + f'(0)u$, and the linear analysis of \S\ref{s:lms} predicts that solutions with localized initial conditions resemble exponential profiles, $u(\xi, t) \sim \rme^{\nu_*(c) \xi} \rme^{\lambda_*(c) t}$, where $\lambda_*(c)$ is the pointwise growth mode with maximal real part, and $\nu_*(c)$ is the associated spatial eigenvalue (keeping track here only of the exponential rates, and neglecting  algebraic pre-factors). On the other hand, a rigid front with speed $c$ is an equilibrium in its co-moving frame, with spatial exponential asymptotics $u_\mathrm{fr}(\xi) \sim \rme^{\nu (c) \xi}$. In order to achieve the matching $\rme^{\nu_*(\xi) \xi} \rme^{\lambda_* (c) t} \sim \rme^{\nu(c) \xi}$, we need $\mathrm{Re} \, \lambda_*(c) = 0$: otherwise the tail will be exponentially growing or decaying in time and cannot be matched with the stationary front profile; see Fig.~\ref{fig:pulled matching}. In other words, we must require pointwise marginal stability in the leading edge.

We additionally assume in Def.~\ref{def: rigid pulled} that this marginal stability arises in the simplest fashion, determined by a simple pinched double root of the dispersion relation. Applying an exponential weight, $v(\xi, t) = \rme^{\eta_\mathrm{lin} \xi} u(\xi, t)$, the linearized equation in the leading edge at speed $c = c_\mathrm{lin}$ becomes
\begin{align}
    v_t = D_\mathrm{eff} v_{\xi \xi} + \mathrm{O}(\partial_\xi^3)v, \label{e: pulled leading edge equation}
\end{align}
where $D_\mathrm{eff} > 0$ is determined solely by the leading order expansion of the dispersion relation near the double root, and by $\mathrm{O}(\partial_\xi^3)v$ we mean terms which include three or more spatial derivatives. In self-similar variables $\tau = \log t, y = \frac{\xi}{\sqrt{t}}$ one can compute solutions to the heat-like equation~\eqref{e: pulled leading edge equation} (see e.g.~\cite{GallayWayne}) and find that localized solutions at leading order resemble the heat kernel
\begin{align*}
    v(\xi, t) \sim \frac{1}{\sqrt{t}} \rme^{-\xi^2/(4D_\mathrm{eff} t)}.
\end{align*}
Then the unweighted tail $u_\mathrm{tail} (\xi, t) = \rme^{-\eta_\mathrm{lin} \xi} v(\xi, t)$ still exhibits pointwise decay at rate $t^{-1/2}$, which obstructs matching with a stationary front. This mismatch is resolved by passing to a logarithmically drifting frame,
\begin{align}
    u_\mathrm{tail} \left(\zeta - \frac{1}{2 \eta_\mathrm{lin}} \log t, t \right) = \rme^{-\eta_\mathrm{lin} \zeta} \rme^{-(\zeta-\frac{1}{2 \eta_\mathrm{lin}} \log t)^2/(4 D_\mathrm{eff} t)}. \label{e: pulled tail matching 1}
\end{align}

\textbf{The matching procedure and the logarithmic delay.} To fully achieve pointwise marginal stability, with no algebraic decay,~\eqref{e: pulled tail matching 1} suggests that we should consider the problem in a logarithmically shifted frame, $\zeta = \xi + \alpha \log (t+t_0) - \alpha \log t_0$, so that our equation~\eqref{e: proofs eqn} reads
\begin{align*}
    u_t = \mathcal{P}(\partial_\zeta) u + \left( c_\mathrm{lin} -  \frac{\alpha}{t+t_0} \right) u_\zeta + f(u),
\end{align*}
leaving $\alpha \in \R$ free for now. The linearized equation at $u \equiv 0$ reads
\begin{align*}
    u_t = \mathcal{P}(\partial_\zeta) u + \left( c_\mathrm{lin} -  \frac{\alpha}{t+t_0} \right) u_\zeta + f'(0) u.
\end{align*}
Applying an exponential weight $v(\zeta, t) = \rme^{\eta_\mathrm{lin} \zeta} u(\zeta, t)$, one finds
\begin{align*}
    v_t = D_\mathrm{eff} v_{\zeta \zeta} - \frac{\alpha}{t+t_0} v_\zeta + \frac{\alpha \eta_\mathrm{lin}}{t+t_0} v + \mathrm{O}(\partial_\zeta^3)v.
\end{align*}
We will use this equation to describe the tail of the solution, say for $\zeta > 0$. Since the state created in the wake of the front is exponentially stable, it is natural to accompany this equation with an absorbing boundary condition, $v = 0$ at $\zeta = 0$.
To see how the logarithmic shift affects the dynamics, we then pass to self-similar variables, $V(y,\tau) = v(\zeta, t)$ with $y = \frac{1}{\sqrt{D_\mathrm{eff}}} \frac{\zeta}{\sqrt{t+t_0}}, \tau = \log (t+t_0)$, solving
\begin{align}
    V_\tau = V_{yy} + \frac{1}{2} y V_y + \alpha \eta_\mathrm{lin} V + \mathrm{O}(\rme^{-\tau/2}) V. \label{e: pulled fronts self similar equation}
\end{align}
The spectrum of the operator $L_\Delta = \partial_y^2 + \frac{y}{2} \partial_y + 1$ is well-known, for instance via its relation to the quantum harmonic oscillator~\cite{GallayWayne}. With the Dirichlet boundary condition $V = 0$ at $y = 0$, the first eigenvalue of $L_\Delta$ is located at $\lambda = 0$, with eigenfunction $V_0(y) = y \rme^{-y^2/4}$. Hence, the leading order behavior of the solution to~\eqref{e: pulled fronts self similar equation} is
\begin{align*}
    V(y, \tau) = \beta \rme^{(\alpha \eta_\mathrm{lin} - 1) \tau} V_0(y)
\end{align*}
for a constant $\beta \in \R$. In $(\zeta, t)$ variables, this corresponds to an approximate tail solution
\begin{align*}
    v_\mathrm{tail}(\zeta, t) = \beta (t+t_0)^{\alpha \eta_\mathrm{lin} - 3/2} \zeta \rme^{-\zeta^2/(4 D_\mathrm{eff} (t+t_0))},
\end{align*}
or,
\begin{align}
    u_\mathrm{tail} (\zeta, t) = \beta (t+t_0)^{\alpha \eta_\mathrm{lin} - \frac{3}{2}} \zeta \rme^{-\eta_\mathrm{lin} \zeta} \rme^{-\zeta^2/(4 D_\mathrm{eff} (t+t_0))}. \label{e: pulled proof u tail}
\end{align}%
\begin{figure}
    \centering
    \includegraphics[width=1\linewidth]{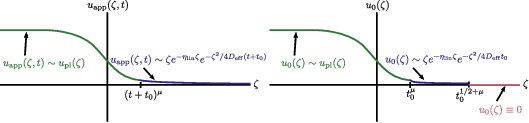}
    \caption{Left: schematic of the approximate solution constructed in the proof of Thm.~\ref{t: rigid pulled unpatterned}. Right: schematic of a steep initial datum obtained by perturbing $u_\mathrm{app}\zeta, 0)$.}
    \label{fig:approx soln}
\end{figure}
To be able to match with the pulled front, which has asymptotics $u_\mathrm{pl}(\zeta) \sim \zeta \rme^{-\eta_\mathrm{lin} \zeta}$, we choose $\alpha = \frac{3}{2 \eta_\mathrm{lin}}$ which eliminates the non-stationary prefactor $(t+t_0)^{\alpha \eta_\mathrm{lin} - 3/2}$. The coefficient $\beta$ can then be chosen to match  the coefficient of $\zeta \rme^{-\eta_\mathrm{lin} \zeta}$ in the asymptotics of $u_\mathrm{pl} (\zeta)$.

Note that the correct choice for the coefficient of the logarithmic delay is $\alpha = \frac{3}{2 \eta_\mathrm{lin}}$ rather than $\alpha = \frac{1}{2 \eta_\mathrm{lin}}$, as~\eqref{e: pulled tail matching 1} would suggest. This is an effect of the absorption provided by the front interface: if the Dirichlet boundary condition $V = 0$ were replaced by a non-absorbing, say Neumann boundary condition $V_y = 0$ at $y = 0$, the first eigenfunction of $L_\Delta$ would be $\rme^{-y^2/4}$ with eigenvalue $\lambda = \frac{1}{2}$, and $u_\mathrm{tail}$ would become
\begin{align*}
    \tilde{u}_\mathrm{tail} (\zeta, t) = \beta \rme^{-\eta_\mathrm{lin} \zeta} (t+t_0)^{\alpha \eta_\mathrm{lin} - 1/2} \rme^{-\zeta^2/(4 D_\mathrm{eff} (t+t_0))},
\end{align*}
and so choosing $\alpha = \frac{1}{2 \eta_\mathrm{lin}}$ would eliminate the non-stationary prefactor. However, note that now the $\zeta$ factor is absent in the tail asymptotics, and so this tail could not be matched with a pulled front $u_\mathrm{pl} (\zeta) \approx (b \zeta + a)\rme^{-\eta_\mathrm{lin} \zeta}$ except in the non-generic case $b = 0$, which corresponds to the pushed-to-pulled transition. Hence, for the typical pulled front $u_\mathrm{pl}$ from Def.~\ref{def: rigid pulled}, the logarithmic delay is $-\frac{3}{2 \eta_\mathrm{lin}} \log t$.

We can finally construct the approximate solution which governs the propagation dynamics by gluing the tail to the front at the length scale $\zeta = (t+t_0)^\mu$, $0 < \mu \ll 1$. More precisely, we set
\begin{align*}
    u_\mathrm{app}(\zeta, t) = \begin{cases}
        u_\mathrm{pl}(\zeta), & \zeta \leq (t+t_0)^\mu, \\
        u_\mathrm{tail}(\zeta, t), & \zeta \geq (t+t_0)^{\mu}+1,
    \end{cases}
\end{align*}
and use cutoff functions to smoothly transition across the region $(t+t_0)^\mu \leq \zeta \leq (t+t_0)^\mu +1$; see Fig.~\ref{fig:approx soln}. We also record the  weighted approximate solution, $v_\mathrm{app}(\zeta, t) = \omega_{0, \eta_\mathrm{lin}} (\zeta) u_\mathrm{app} (\zeta, t)$.
To measure the accuracy of our weighted approximate solution, we define the residual evaluator
\begin{align*}
    \tilde{F}_\mathrm{res} [u] = u_t - \mathcal{P}(\partial_\zeta) u + \left( c_\mathrm{lin} - \frac{3}{2 \eta_\mathrm{lin}(t+t_0)} \right) u_\zeta + f(u),
\end{align*}
and the weighted residual evaluator, $F_\mathrm{res}[v]= \omega_{0, \eta_\mathrm{lin}} \tilde{F}_\mathrm{res}\left[\frac{v}{\omega_{0, \eta_\mathrm{lin}}}\right]$. One finds~\cite[Proposition 2.5]{as1} the residual error estimate
\begin{align*}
    \| \rho_{0,1} F_\mathrm{res}[v_\mathrm{app}]\|_{L^1} \lesssim (t+t_0)^{-1/2 + 4 \mu},
\end{align*}
where $\rho_{0,1}$ is a one-sided algebraic weight satisfying $\rho_{0,1}(\zeta) = 1$ for $\zeta \leq -1$ and $\rho_{0,1}(\zeta) = \zeta$ for $\zeta \geq 1$. Hence, if we choose the length scale to satisfy $0 < \mu < \frac{1}{8}$ and choose $t_0 \gg 1$, $v_\mathrm{app}$ solves our weighted equation with small residual error for all times $t > 0$.

\textbf{Stability argument.} To complete the proof of Thm.~\ref{t: rigid pulled unpatterned}, one inserts an ansatz $v(\zeta, t) = v_\mathrm{app}(\zeta, t) + w(\zeta, t)$ into the weighted equation $F_\mathrm{res} [v] = 0$. The goal is to show that if the initial perturbation $w_0$ is small in an appropriate norm (but one which allows for cutting off the tail of $v_\mathrm{app}$), then $w(\zeta, t)$ remains small for all time, so that the real solution resembles the approximate solution $v(\zeta, t) \approx v_\mathrm{app}(\zeta, t)$.

The only structure which we can exploit to control $w$ in this general setting is the spectral information encoded in Def.~\ref{def: rigid pulled}, which leads to estimates on the linearized evolution $\rme^{\mathcal{L}^\omega_\mathrm{pl} t}$, where $\mathcal{L}^\omega_\mathrm{pl} u = \omega_{0, \eta_\mathrm{lin}} \mathcal{L}_\mathrm{pl} \left( \frac{u}{\omega_{0, \eta_\mathrm{lin}}} \right)$ is the weighted linearization about $u_\mathrm{pl}$. Although we are perturbing from $v_\mathrm{app} (\zeta, t)$, which encodes more detailed tail dynamics than $\omega_{0, \eta_\mathrm{lin}} (\zeta) u_\mathrm{pl}(\zeta)$, we observe that since the matching point $\zeta = (t+t_0)^\mu$ grows as $t \to \infty$, we have at least pointwise $u_\mathrm{app} (\zeta, t) \to u_\mathrm{pl} (\zeta)$ as $t \to \infty$. More precisely, one has the quantitative estimate $\| \rho_{0,1} [u_\mathrm{app} (\cdot, t) - u_\mathrm{pl}] \|_{L^\infty} \lesssim (t+t_0)^{-1/2 + 4 \mu}$~\cite[Lemma 2.8]{as1}. This allows one to replace the linearized term $f'(u_\mathrm{app}(\cdot, t)) w$ in $F_\mathrm{res}[v_\mathrm{app} + w]$ with $f'(u_\mathrm{pl})w$, at the price of an additional error term $E_1 = [f'(u_\mathrm{app}(\cdot, t))-f'(u_\mathrm{pl})]w$, hoping to view the linearized evolution near the front itself, $w_t = \mathcal{L}^\omega_\mathrm{pl} w$, as the leading part of the equation. One may then write an evolution equation for the perturbation $w$, with the form
\begin{align}
    w_t = \mathcal{L}^\omega_\mathrm{pl} w - \frac{3}{2 \eta_\mathrm{lin} (t+t_0)} \left[ w_\zeta + \omega_{0, \eta_\mathrm{lin}} \partial_\zeta \left(\frac{1}{\omega_{0, \mathrm{lin}}} \right) w \right] + R + E_1 + N(w), \label{e: pulled proof w equation}
\end{align}
where $R = F_\mathrm{res}[v_\mathrm{app}]$ and $N(w) = \mathrm{O} \left(\frac{1}{\omega_{0, \eta_\mathrm{lin}}} w^2\right)$ is a quadratic nonlinear remainder which carries a localizing factor from the exponential weight.

The linearized problem $w_t = \mathcal{L}^\omega_\mathrm{pl} w$ is comparable to the heat equation on the half-line with a Dirichlet boundary condition. Note also that for $\zeta \geq 1$, the factor multiplying $-\frac{3}{2 \eta_\mathrm{lin} (t+t_0)}$ becomes $w_\zeta - \eta_\mathrm{lin} w$. So, a good but heavily simplified model problem for~\eqref{e: pulled proof w equation}, neglecting the last three terms, is
\begin{align}
    \begin{cases}
        w_t = w_{\zeta \zeta} - \frac{3}{2 (t+t_0)} w_\zeta + \frac{3}{2(t+t_0)} w, & \zeta > 0, t > 0 \\
        w = 0, & \zeta = 0, t > 0, \\
        w = w_0, & \zeta > 0, t = 0.
    \end{cases} \label{e: selection model}
\end{align}
However, for instance computing expansions in self-similar variables, one finds that generic solutions to this equation do not decay in time, even for very localized initial conditions. This is a consequence of the term $\frac{3}{2(t+t_0)} w$, and is the exact same effect that informed our choice of the logarithmic delay, so that the tail dynamics could be marginally stable. This presents a challenge for controlling the dynamics in~\eqref{e: pulled proof w equation}: we want to view $w_t = \mathcal{L}^\omega_\mathrm{pl} w$ as the principle part of the equation, use spectral information to obtain estimates on $\rme^{\mathcal{L}^\omega_\mathrm{pl} t}$, and control all other terms perturbatively. But solutions to $w_t = \mathcal{L}^\omega_\mathrm{pl} w$ decay in time in suitable norms, while we should now expect solutions to~\eqref{e: pulled proof w equation} not to decay at all, and so we cannot hope to treat the additional terms perturbatively.

The solution to this obstacle presented in~\cite{as1} is to consider instead the variable $z = (t+t_0)^{-3/2} w$, which in the case of the model problem~\eqref{e: selection model} solves the equation
\begin{align*}
        z_t = z_{\zeta \zeta} - \frac{3}{2 (t+t_0)} z_\zeta, \quad  \zeta > 0, t > 0, \qquad \text{with boundary condition }
        z = 0, \quad \zeta = 0, t > 0,
\end{align*}
and initial condition {\color{black} $z(0,\zeta)=z_0(\zeta) = t_0^{-3/2} w_0(\zeta)$}. If one can prove that $z$ decays with rate $(t+t_0)^{-3/2}$ for large times, for small initial data, then $w$ remains bounded and small for all times, which is enough to conclude Thm.~\ref{t: rigid pulled unpatterned}. Key to this proof are the linearized estimates, valid for $t > 0$,
\begin{align*}
    \| \rho_{0, -1} \rme^{\mathcal{L}^\omega_\mathrm{pl} t} z_0 \|_{L^\infty} &\lesssim t^{-3/2} \| \rho_{0,1} z_0 \|_{L^1},\\ \| \rho_{0, 1} \partial_\zeta \rme^{\mathcal{L}^\omega_\mathrm{pl} t} z_0 \|_{L^1} &\lesssim t^{-1/2} \| \rho_{0, 2+\mu} z_0 \|_{L^\infty}, \\ \| \rho_{0, 2+\mu} \partial_\zeta \rme^{\mathcal{L}^\omega_\mathrm{pl} t} z_0 \|_{L^\infty} &\lesssim t^{\frac{1-\mu}{2}} \| \rho_{0, 2 + \mu} z_0 \|_{L^\infty},
\end{align*}
where $\rho_{0,r}(\zeta)$ is a one-sided algebraic weight satisfying $\rho_{0,r} (\zeta) = 1$, $\zeta \leq -1$ and $\rho_{0, r} (\zeta) = \zeta^r, \zeta \geq 1$. The first estimate is used to establish $(t+t_0)^{-3/2}$--decay of $z(t)$, after also perturbatively controlling the term $-\frac{3}{2 \eta_\mathrm{lin} (t+t_0)} z_\zeta$ via  the second two estimates and a bootstrap procedure. All linear estimates are inherited from the fact that the dynamics in the leading edge resemble the heat equation on the half-line (and there is no additional critical spectrum which competes with these dynamics), but proven via inverse Laplace estimates and a careful analysis of the resolvent operator near the essential spectrum; see~\cite[Sections 3-4]{as1} and~\cite{AveryScheel} for details.

Eventually, the nonlinear stability argument for $z(t)$ may be closed, establishing that $\| \rho_{0, -1} w(t) \|_{L^\infty}$ remains small for all time provided $\|\rho_{0, 2+ \mu} w_0 \|_{L^\infty}$ is small; see~\cite[Section 5]{as1} for details. This condition allows for some steep initial data for $u$ to be included, for instance
\begin{align*}
    u_0(\zeta) = \begin{cases} u_\mathrm{app} (\zeta, 0), & \zeta < t_0^{\frac{1}{2}+\mu}, \\
    0, & \zeta \geq t_0^{1/2+\mu},
    \end{cases}
\end{align*}
or perturbations thereof, small in the norm $\| \rho_{0, 2 + \mu} \cdot \|_{L^\infty}$; see Fig.~\ref{fig:approx soln} and~\cite[Section 6]{as1}.

\subsection{Practical considerations: nonlinear selection}\label{s:a3}


Defs.~\ref{def: rigid pushed} and~\ref{def: rigid pulled} give specific, verifiable conditions that are sufficient to establish selection of pushed and pulled fronts, respectively. Practically, one would then wish to establish existence of fronts that satisfy these conditions. In that light, the utility of Thms.~\ref{t: pushed stability} and~\ref{t: rigid pulled unpatterned} is to reduce the difficult PDE dynamics question of front selection to verifying  conditions in nonlinear and linearized traveling-wave equations, building for instance on the techniques outlined in \S\ref{s:robust} for perturbative or computational considerations. Beyond the open conjectures mentioned above, there are a number of other subtleties and questions.

\textbf{Basins of attraction.} Thms.~\ref{t: pushed stability} and~\ref{t: rigid pulled unpatterned} guarantee that pushed or pulled fronts attract an open class of steep initial data. The results are however perturbative and there can be many selected fronts~\cite{pnp}: it would be interesting to establish attractivity more globally, based for instance on variational structures~\cite{gallayrisler}, or to understand the basin boundaries of selected fronts.

\textbf{Orbital stability of pulled fronts.} The selection result for pushed fronts, Thm.~\ref{t: pushed stability} gives an asymptotic phase: it identifies a specific spatial translate of the front profile to which the solution converges as $t \to \infty$. By contrast, Thm.~\ref{t: rigid pulled unpatterned} only establishes asymptotic closeness to the front profile and bounds on fluctuations along translates. We suspect that with additional effort Thm.~\ref{t: rigid pulled unpatterned} could be improved to guarantee an asymptotic phase in the selection of pulled fronts.

\textbf{Beyond stability in a fixed weight.}
Defs.~\ref{def: rigid pushed} and~\ref{def: rigid pulled} are sufficient but not necessary. For instance, marginal stability need only be satisfied in a pointwise sense, rather than in an exponentially weighted space. We refer to \S\ref{s:res} for a more general characterization of marginal stability in the leading edge through resonance criteria that generalize the 1:1 resonance encoded in pinched double roots; see in particular~\cite{FayeHolzerScheel}.

\textbf{Pattern-forming fronts.} Thms.~\ref{t: pushed stability} and~\ref{t: rigid pulled unpatterned} require exponential stability in the wake, which rules out pattern-forming fronts, the subject of the next section.

\section{Pattern forming fronts --- setup, results, and conjectures}\label{sec: modulated}

Our discussion in \S\ref{s:nlms}--\S\ref{s:nlmsp} has focused on traveling waves, that is, on equilibria in a comoving frame, which connect the unstable trivial equilibrium at $x=+\infty$ to a stable equilibrium at $x=-\infty$. In this chapter, we discuss important generalizations that in particular allow us to address questions of pattern formation and selection. First, we illustrate the difficulties that arise when studying pattern forming fronts in the simplest case of rigid traveling fronts that connect to a periodic pattern at $x=-\infty$ in \S\ref{s:5.0}. We then address fronts in oscillatory invasion processes, arising most prominently since linear marginal stability is induced by pinched double roots with $\lambda=\pm\rmi\omega_*$ in \S\ref{s:5.1}, followed by a discussion of marginal stability for these fronts in \S\ref{s:5.2}. We briefly discuss the intricate question of selection of patterns in the wake, based on resonances between invasion frequencies and frequencies of wave trains in the wake in \S\ref{s:5.3}. \textcolor{black}{
We characterize pushed and pulled pattern-forming fronts, discuss the difficulties that come with attempting to establish selection of these fronts, and both formulate conjectures and point to some results on the selection of marginally stable fronts in 
\S\ref{s:5.4}. We discuss variations which all include genuinely infinite-dimensional, that is, PDE rather than ODE, existence problems in sections 
\S\ref{s:infdim} and \S\ref{s:5.5}. The latter section includes the largely unresolved question of the selection of patterns in front invasion in higher space-dimension with patterning transverse to the 
front interface.  We summarize some simple recipes for pattern prediction in \S\ref{s:modpract} and overview results in the literature on existence and stability in \S\ref{s:5:exples}.}

\textcolor{black}{Compared with the analysis and results described in previous sections, many central questions raised here remain open. We describe however a somewhat well established framework to tackle existence and stability questions, as well as predicting selected wavenumbers, but establishing the actual selection of fronts and thereby specific patterns in the wake remains a largely open question.}

\subsection{Pattern-forming fronts}\label{s:5.0}
So far, we have simply assumed that the instability of the trivial state is mediated by a front that creates a stable new equilibrium, $u(x)\equiv u_-$ in its wake. This may however not always be the case, even when linear marginal stability is associated with a pinched double root $\lambda=0$, $\nu<0$, hence with a stationary leading edge. One can easily envision examples where there simply are no equilibria other than the unstable origin. A well known example, studied in~\cite{cartersch,avery2023stabilitycoherentpatternformation}, is an array of  van der Pol oscillators, with diffusive coupling, also known as the FitzHugh-Nagumo equation in the oscillatory regime~\eqref{e:fhn},
\begin{equation}\label{e:vdp}
\begin{aligned}
u_t&=u_{xx}+u(1-u)(u-a) -v,\nonumber\\
v_t&=\eps(u-\gamma v + b),\nonumber
\end{aligned}
\end{equation}
with $a<0$, $0<\eps\ll 1$ and for instance $b=0$, $0\leq \gamma<4/(a-1)^2$. The latter condition implies that nullclines intersect only at the origin. Trajectories of the $x$-independent ODE converge to a stable limit cycle. On the other hand, the calculation of linear spreading speeds in~\eqref{e:fhnspsp} shows that dynamics of invasion fronts are expected to be stationary in the leading edge. It turns out that~\cite{cartersch,avery2023stabilitycoherentpatternformation} the invasion process is mediated by  fronts $(u_*,v_*)(x-c_* t)$ that converge to  periodic profiles as $\xi=x-c_*t\to-\infty$,
\[
(u_*,v_*)(\xi)-(u_\mathrm{per},v_\mathrm{per})(\xi)\to 0, \text{ for } \xi\to-\infty,\qquad (u_\mathrm{per},v_\mathrm{per})(\xi)=(u_\mathrm{per},v_\mathrm{per})(\xi+L) \text{ for some } L>0;
\]
see Fig.~\ref{f:cglfhn} for an illustration.
Similar examples arise in parametrically forced coupled oscillators,
\[
A_t=(1+\rmi\alpha)A_{xx}+(1+\rmi\omega)A-(1+\rmi\beta) A|A|^2+\gamma \bar{A},
\]
with $|\alpha|,|\omega|\ll 1$, $1\gg |\gamma|>|\alpha+\omega|$, so that the linear spreading speed predicts stationary behavior in the leading edge; see Fig.~\ref{f:ressliver} and the accompanying discussion. For $|\beta|\gg 1$, nonlinear oscillations are ``unlocked'' and there are no equilibria but rather a periodic orbit, so that the stationary leading edge necessarily leaves an oscillatory pattern in its wake.

In these examples, we would  like to identify  fronts that are selected in the sense of \S\ref{s:3.2}. We therefore first would like to characterize marginal stability of fronts, which then in a second step could be used to establish nonlinear selection, in analogy to our approach in \S\ref{s:nlmsp}. For the first step, we need to investigate spectral properties of a front leaving behind a periodic pattern. This yields an equation for the pointwise resolvent of the form~\eqref{e:res1stx},
\begin{equation}
U_\xi=M_\lambda(\xi)U+\delta(\xi-y)I,
\end{equation}
where $M_\lambda(\xi)\to M_{+,\lambda}$ for $\xi\to\infty$ and $M_\lambda(\xi)-M_{-,\lambda}(\xi)\to 0$ with periodic $M_{-,\lambda}(\xi)=M_{-,\lambda}(\xi+L)$ for some period $L>0$. Exponential asymptotics of solutions are dictated by eigenvalues $\nu$ of $M_{+,\lambda}$ and Floquet exponents $\nu$ of $M_{-,\lambda}(\xi)$. Since the $\xi$-derivative of the periodic solution induces a periodic solution to $U_\xi=M_{-,0}U$, $\nu=0$ is always a Floquet exponent at $\lambda=0$. In other words, the periodic pattern in the wake is at best marginally stable. In the simplest case, the multiplier $\nu=0$ is simple, and depends on $\lambda$ analytically through $\lambda=-c_\mathrm{g}\nu+\rmO(\nu^2)$, where the proportionality factor $c_\mathrm{g}\in\R$ is usually referred to as the group velocity. Introducing exponential weights $\omega(\xi)=\rme^{-\eta \xi},\xi<-1$,  then shifts the real part of the spectrum at $-\infty$ to $\Re\lambda=-c_\mathrm{g}\eta+\rmO(\eta^2)$, which is stable when $c_\mathrm{g}\eta>0$. For the fronts we are interested in here, $c_\mathrm{g}<0$, so that $\eta<0$ stabilizes the essential spectrum at $-\infty$. More intuitively, one can think of the linear relation $\lambda=-c_\mathrm{g}\nu$ as corresponding to the advection equation $u_t=-c_\mathrm{g} u_\xi$ with solution $u(\xi-c_\mathrm{g}t)$, which simply shifts initial conditions with speed $c_\mathrm{g}$, to the left when $c_\mathrm{g}<0$. The  weight, $\rme^{-\eta \xi}$, $\eta<0$,  is exponentially decreasing as $\xi\to-\infty$, thus exponentially dampening the transport to the left.
\begin{figure}
    \centering
    \includegraphics[width=0.8\textwidth]{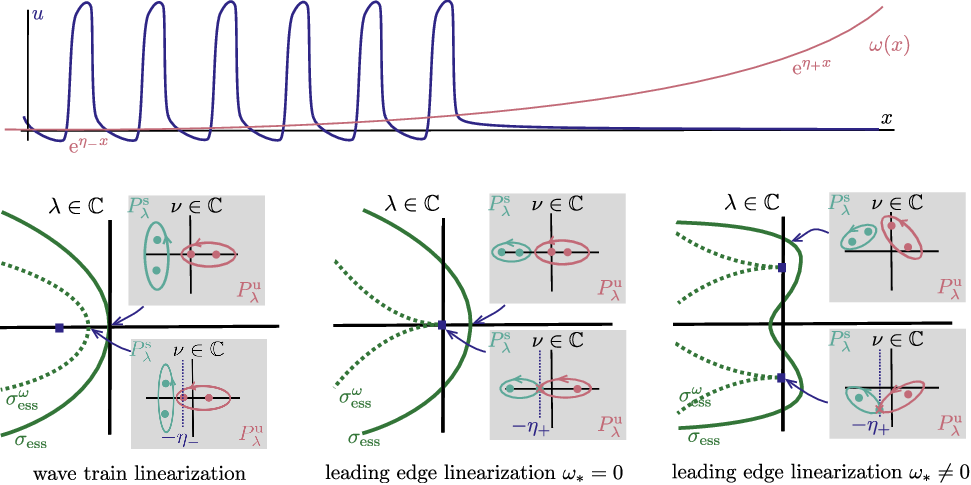}
    \caption{A traveling wave $u$ leaving a periodic pattern in its wake (top) together with weight $\omega$ that renders the linearization marginally stable. Also shown (bottom) are spatial exponential growth rates $\nu$ for solutions to the eigenvalue problem depending on $\lambda$. Exponential growth rates correspond to Floquet exponents in the wake. Note how the $\nu=0$ moves to the left as $\lambda$ is decreased, corresponding to $c_\mathrm{g}=\frac{\rmd\nu}{\rmd\lambda}=-1/c_\mathrm{g}>0$. In the leading edge, both stationary and oscillatory marginal stability are shown. }\label{f:modhalo}
\end{figure}
Unfortunately, this choice of weight allows for exponential growth of solutions and is therefore incompatible with a nonlinear analysis. A much more refined control of the effect of perturbations on the nonlinear patterns is necessary; see~\cite{avery2023stabilitycoherentpatternformation} for an in-depth discussion, references, and analysis.

We emphasize however that the \emph{marginal stability in the wake} is apparently not relevant for the selection of the front speed, but it is the \emph{marginal stability in the leading edge}, induced by pinched double roots or by point spectrum, that determines the invasion speed for compactly supported initial conditions.

\subsection{Modulated fronts: existence}\label{s:5.1}

In most examples, pattern formation already manifests itself through oscillations in the leading edge dynamics so that invasion fronts are inherently time-periodic or ``modulated'' in a comoving frame,
\begin{equation}\label{e:mf}
    u(t,x)=u_\mathrm{mtw}(x-ct,\omega t;\omega,c), \qquad u_\mathrm{mtw}(\xi,\tau;\omega,c)=u_\mathrm{mtw}(\xi,\tau+2\pi;\omega,c).
\end{equation}
The shape of the waves varies periodically in time with period $2\pi/\omega$; see Fig.~\ref{f:shch} and Fig.~\ref{f:cglfhn} (left panels) for pattern-forming modulated waves as opposed to Fig.~\ref{f:cglfhn} (right panels) for rigidly propagating fronts. Due to translation invariance in time and space, solutions come in 2-parameter families $u_\mathrm{mtw}(\xi+\xi_0,\tau+\tau_0)$. In addition, we shall find that one can vary the parameters $\omega$ and $c$, similar to the way in which one can vary the parameter $c$ in traveling-wave problems in \S\ref{s:nlms}, and find a 4-parameter family of solutions.
%
Modulated traveling waves solve a PDE rather than a traveling-wave ODE,
\begin{equation}\label{e:mtweq}
    \omega u_\tau=\mathcal{P}(\partial_\xi)u+c\partial_\xi u + f(u),\qquad u(\xi,\tau)=u(\xi,\tau+2\pi),\qquad \lim_{\xi\to\infty}u(\xi,\tau)=0.
\end{equation}
The principal part, $\partial_\tau-\mathcal{P}(\partial_\xi)$ is pseudo-elliptic rather than parabolic due to the periodic-in-time boundary conditions. This can be seen most directly when inspecting the simple linear equation
\begin{equation}\label{e:heatper}
u_\tau=u_{\xi\xi} -u +g(\tau,x), \qquad u(\xi,\tau)=u(\xi,\tau+2\pi)
\end{equation}
which, using Fourier-series in $\tau$ and Fourier transform in $\xi$, $u(\tau,\xi)=\sum_\ell \int_k \rme^{\rmi\ell\tau+\rmi k \xi} \hat{u}(\ell,k)$, gives
\[
\hat{u}(\ell,k)=(1+\rmi\ell+k^2)^{-1}\hat{g}(\ell,k),
\]
with bounded symbol $(1+\rmi\ell+k^2)^{-1} $ in $(\ell,k)\in\R^2$. The equation then is well-posed as a boundary-value problem similar to the elliptic operator $\partial_{\tau\tau}+\partial_{\xi\xi}-1$ with also bounded symbol of the inverse $(1+\ell^2+k^2)^{-1} $.

To interpret~\eqref{e:mtweq} as a dynamical system in $\xi$ on a phase space, one would write this equation as a first-order equation in $\xi$ and for instance expand in Fourier series in $\tau$ for an infinite-dimensional phase space. This infinite-dimensional initial-value problem is  strongly ill-posed in the sense of Hadamard, which can be seen in~\eqref{e:heatper} with $g=0$, and solutions $\hat{u}_\ell(\xi)=c_+\rme^{\sqrt{1+\rmi\ell}\xi}+c_-\rme^{-\sqrt{1+\rmi\ell}\xi}$, of arbitrarily fast exponential growth as $\xi\to+\infty$ or $\xi\to-\infty$. The problem is however hyperbolic for $\ell\neq 0$  in the sense that solutions either grow or decay exponentially. This is true more generally and analytic tools to pursue questions in \S\ref{s:3.1} on existence and continuation of fronts as heteroclinic orbits are available for  the infinite-dimensional modulated-wave equation~\eqref{e:mtweq}. We  briefly summarize the literature here.

\textbf{Stable manifolds.} Invariant manifolds near a given global trajectory were constructed in~\cite{PSS97},~\cite{LPSS2000}  and~\cite{SS2001}, connecting Fredholm properties of the linearized boundary-value problems to dimensions of stable and unstable subspaces. Since dimensions are infinite, one defines \emph{relative dimensions} either through homotopies or through Fredholm indices of projections onto reference subspaces.



\textbf{Counting arguments.} Multiplicities of front solutions as derived in \S\ref{s:3.1} can be analyzed for modulated waves. Key to computing relative dimensions of stable and unstable manifolds are homotopies during which one tracks spatial eigenvalues $\nu$ crossing the imaginary axis. The direction of crossing can usually be related to group velocities, which are both analytically accessible and physically relevant information; we refer to~\cite{SS2008} and~\cite{SS2004} for an in-depth discussion with many examples. In our situation, this counting yields the aforementioned multiplicity with free parameters $\xi_0, \tau_0,\omega$, and $c$.
We notice in passing that these counting arguments are also key to numerical methods based on farfield-core decompositions described in \S\ref{s:robust} that allow one to efficiently compute solutions to~\eqref{e:mtweq}; see for instance~\cite{morrisseyscheel,lloydscheel,aggms} for related implementations of this strategy.


\textbf{Large speeds.}
If one starts the homotopy at large speeds, one finds a singularly perturbed PDE rather than an ODE. As an example, consider the special case of~\eqref{e:mtweq}, the second-order system
\[
\omega u_\tau - u_{\xi\xi}+cu_\xi+f(u) = 0,
\]
with $c=1/\eps$. If we envision a front that connects to a stationary pattern in the steady frame $u_\mathrm{per}(kx)$, we find that this pattern in the fast comoving frame is of the form $u(k(\xi+ct))=u(k(\xi+\frac{c}{\omega}\tau))$, so that we need to require $\omega=ck=k/\eps$, which gives
\[
ku_\tau+u_\xi-\eps u_{\xi\xi}+\eps f(u)=0.
\]
It is now convenient to introduce a shear transformation so that the pattern $u$ is in fact independent of $\xi$, setting $k\xi+\tau=\sigma$, and rescale $\xi=\zeta/\eps$ which gives
\[
0=(\eps \partial_\zeta + k \partial_\sigma)^2u + \partial_\zeta u + f(u)=0.
\]
Formally setting $\eps=0$ then gives the backward parabolic equation
\begin{equation}\label{e:plimit}
-\partial_\zeta u= k \partial_{\sigma\sigma}u + f(u).
\end{equation}
This parallels the limit of pure backward kinetics that we found from~\eqref{e:gstp} using geometric singular perturbation theory. The limit $\eps\to 0$ here can be justified using invariant manifold constructions; see~\cite{S2006} and, for a fourth-order equation,~\cite{S2017}. Both build on earlier results for convection-dominated elliptic equations~\cite{CMS1993,S1996}. As a consequence, in a situation where the origin is unstable in the forward-dynamics of~\eqref{e:plimit}, with unstable manifold limiting on a stable pattern, periodic in $\sigma$, we find this same heteroclinic behavior in the traveling-wave equation. We note in passing that this limiting behavior also indicates  transversality of the heteroclinic and hence persistence as $\omega$ and $c$ are varied.

\subsection{Modulated fronts: stability and marginal stability}\label{s:5.2}

\textbf{Linear stability.} Stability of time-periodic solutions is determined by temporal Floquet multipliers and properties of the pointwise resolvent of the period map. An approach to stability that parallels the viewpoint taken in \S\ref{s:3.2} was laid out in~\cite{SS2001}. Paralleling the existence problem, one finds an ill-posed pseudo-elliptic generalized eigenvalue problem that gives a geometric perspective on pointwise Green's functions. We summarize briefly here the main results linking spectral properties of period maps and the geometry in an ill-posed initial-value problem. Given a modulated traveling wave solution $u_*(\xi,\tau)$ to~\eqref{e:mtweq}, we are interested in spectral properties of the period map $\Psi$,  defined through
\begin{equation}\label{e:mtweqlin}
    \omega u_\tau=\mathcal{P}(\partial_\xi)u+c\partial_\xi u + f'(u_*(\xi,\tau))v,\qquad u(\xi,0)=u_0(\xi),\qquad u(\xi,2\pi)=:(\Psi u_0)(\xi).
\end{equation}
The Floquet multiplier eigenvalue problem $\Psi u_0=\rho u_0$ can be rewritten as, setting $\rho=\rme^{2\pi\lambda/\omega} $,
\begin{equation}\label{e:mtweqlineig}
   0=\mathcal{T}_\lambda v:=-\omega v_\tau-\lambda v+\mathcal{P}(\partial_\xi)v+c\partial_\xi v + f'(u_*(\xi,\tau))v,\ \ v(\xi,0)=v(\xi,2\pi), \ \ \text{ where }  u(x,\tau)=\rme^{\lambda\tau/\omega}v(x,\tau),
   \end{equation}
It then turns out that spectral properties of  $\Psi$ and of the eigenvalue pencil $\mathcal{T}_\lambda$ coincide, and that in particular the pointwise Green's function of $\Psi-\rho$ can be constructed directly in terms of the Green's function of $\mathcal{T}_\lambda$; {\color{black} see \cite{SS2001} for proofs}.

\textbf{Marginal stability in the leading edge.} The linearization of~\eqref{e:mtweq} at the origin
can be solved explicitly using Fourier series in $\tau$. For $u=\hat{u}_\ell\rme^{\rmi\ell\tau}$, we find exponential behavior $\rme^{\nu\xi}$ when
$d(\nu):=-\omega\rmi\ell+\mathcal{P}(\nu)+c\nu+f'(0)$ has a nontrivial kernel. At the linear spreading speed $c_\mathrm{lin}$ and with $\omega=\omega_\mathrm{lin}/m$ for some $m\in\N$, we find a double root $\nu$ to $d(\nu)$ associated with the pinched double root for $\ell=m$. The pinching condition usually ensures that at this point the (strong) stable manifold with decay $\rme^{(\nu_\mathrm{lin}+\eps)\xi}$, thus including fronts with asymptotics from the double root, has relative dimension 2, that is, in a homotopy of the linearized equation with eigenvalue parameter $\lambda$ it includes two roots $\nu$ that are unstable for $\lambda\gg 1$. As a consequence,  we expect transverse intersections and robust heteroclinics; see Fig.~\ref{f:modhalo}. For nearby values of $\omega$ and $c$, decay is generically weaker: due to the square-root singularity, there is always a root with larger real part than $\nu_\mathrm{lin}$ for nearby values that limits on $\nu_\mathrm{lin}$ as $c\to c_\mathrm{lin}$ and $\omega\to\omega_\mathrm{lin}$.
The behavior of $\nu$ near this double root is at leading order a square root, mimicking diffusive-dispersive asymptotics
\[
\lambda-\rmi\omega_\mathrm{lin}=d_\mathrm{eff}(\nu-\nu_\mathrm{lin})^2 \quad \longrightarrow \quad A_t=d_\mathrm{eff}A_{xx},
\]
with effective diffusivity $\Re d_\mathrm{eff}$ and dispersion $\Im d_\mathrm{eff}$.

\textbf{Pushed fronts.}
Pushed fronts arise as intersections of stable and unstable manifolds with relative dimension 0. Since intersections occur in 2-parameter families due to $\xi$- and $\tau$-translation symmetries, they are codimension-2 heteroclinic orbits and generically necessitate 2 parameters, which are $c$ and $\omega$. The relation between transverse crossing of manifolds in these 2 parameters and geometric and algebraic multiplicity 2 of the 0 Floquet exponent was established in~\cite{SS2001}. An analysis comparable to the bifurcation analysis in \S\ref{s:robust} has not been carried out for modulated invasion fronts.

\subsection{Modulated fronts: selection in the wake}\label{s:5.3}

The discussion thus far has largely ignored the selection of a pattern in the wake. Focusing for now on pulled fronts, we have fixed $c=c_\mathrm{lin}$ and $\omega=\omega_\mathrm{lin}$.


\textbf{Periodic wave trains.} Wave trains are solutions $u=u(\omega t-kx)$, with $u(\zeta)=u(\zeta+2\pi)$ and thus lead from~\eqref{e:parnl} to the $2\pi$-periodic ODE boundary-value problem
\begin{equation}\label{e:wt}
  \mathcal{P}(k\partial_\zeta)u+\omega\partial_\zeta u + f(u)=0,\qquad u(\zeta)=u(\zeta+2\pi).
\end{equation}
The linearization $\mathcal{L}_\mathrm{wt} u= \mathcal{P}(k\partial_\zeta)u+\omega\partial_\zeta u + f'(u_\mathrm{wt})u$ at a potential solution $u_\mathrm{wt}$, has a nontrivial kernel that contains $u_\mathrm{wt}'$. We assume that the kernel is minimal, one-dimensional, and also that the zero-eigenvalue is algebraically simple so that $u'_\mathrm{wt}$ is not in the range of $\mathcal{L}_\mathrm{wt}$. Considering~\eqref{e:wt} as a functional equation, $\mathcal{F}_\mathrm{wt}[u,\omega;k]=0$, $\mathcal{F}_\mathrm{wt}:H^{2m}_\mathrm{per}\times\R^2\to L^2_\mathrm{per}$, we then find that $\partial_{u,\omega}\mathcal{F}_\mathrm{wt}$ is onto at the wave train solution, which implies the existence of a locally unique family of solutions (up to translates in $\zeta$) $u_\mathrm{wt}(kx-\omega_\mathrm{nl}(k) t;k)$ for some $\omega_\mathrm{nl}(k)$, usually referred to as the nonlinear dispersion relation; see~\cite{dsss,SS2004} for details.

\textbf{Selection through resonances.} In the frame with speed $c_\mathrm{lin}$, the frequency of wave trains includes a Doppler shift through $\omega_\mathrm{nl,co}(k):=\omega_\mathrm{nl}(k)-c_\mathrm{lin} k$. Time-periodicity of the front requires commensurate frequencies in the wake and in the leading edge,
\begin{equation}\label{e:reswake}
    p \omega_\mathrm{lin}= q \omega_\mathrm{nl,co}(k),\qquad \text{for some } p,q\in\N.
\end{equation}
One often observes $p=q=1$, which can be interpreted as \emph{node conservation} in that one oscillation in the leading edge creates one oscillation in the wake~\cite{deelanger,vanSaarloosReview}. This is clearly not true in for instance Cahn-Hilliard and FitzHugh-Nagumo examples; see Fig.~\ref{f:cglfhn} and~\ref{f:shch}.

Transitions between harmonic and subharmonic invasion, for instance $p=1$, increasing $q$ from $1$ to $2$, lead to interesting heteroclinic bifurcations. Fronts with frequency $\omega_\mathrm{lin}$ lie in an invariant subspace of the dynamics of~\eqref{e:mtweq} with $\omega=\omega_\mathrm{lin}/q$. Transitions to harmonics can then be thought of as heteroclinic bifurcations in a direction perpendicular to this invariant subspace; see~\cite{SS2007} for a technically similar analysis.
More general unlocking with temporally homoclinic limits has not been studied.

An extreme example is the case when $\omega_\mathrm{lin}=0$, which simply gives, for $m\neq 0$, $\omega_\mathrm{nl,co}(k)=0$, that is, the pattern in the wake would be stationary in the comoving frame. Choosing $p=q=0$ in~\eqref{e:reswake} however allows for \emph{any} frequency $\omega_\mathrm{nl,co}$. This transition from stationary to oscillatory front motion has been observed in the FitzHugh-Nagumo equation as one approaches the Canard regime~\cite{cartersch}.


\subsection{The marginal stability conjecture for pattern-forming fronts}\label{s:5.4}

There are currently no results that establish selection of pushed or pulled pattern-forming fronts analogous to Thms.~\ref{t: pushed stability} or~\ref{t: rigid pulled unpatterned}, with the exception of some recent progress for stability and, in the pushed case, selection for rigid pattern-forming fronts. We summarize these results and formulate conjectures for modulated fronts, here. We again consider systems of parabolic equations,
\begin{align}
    u_t = \mathcal{P}(\partial_\xi) u + c u_\xi + f(u), \quad u = u(\xi,t) \in \mathbb{R}^N, \quad \xi \in \mathbb{R}, \quad t>0. \label{e: pf proofs eqn}
\end{align}
To make spectral stability of the pattern in the wake precise, let $\mathcal{L}_- = \mathcal{P}(\partial_\xi) + c \partial_\xi + f'(u_-(\xi))$ denote the linearization of~\eqref{e: pf proofs eqn} about a periodic pattern $u_-(\xi)$, with wavenumber $k_*$ and period $L_* = \frac{2\pi}{k_*}$. By Floquet-Bloch theory \cite{ReedSimon}, $\mathcal{L}_- : H^{2m} (\mathbb{R}) \subset L^2(\mathbb{R}) \to L^2(\mathbb{R})$ is conjugate to a family of Bloch operators
\begin{align*}
    \hat{\mathcal{L}}_- (\nu) : H^{2m} (\mathbb{R}/L_*\mathbb{Z}) \subset L^2(\mathbb{R}/L_* \mathbb{Z}) \to L^2(\mathbb{R}/L_* \mathbb{Z}), \qquad \hat{\mathcal{L}}_- (\nu) = \mathcal{P}(\partial_\xi + \nu) + c (\partial_\xi + \nu) + f'(u_-(\xi)),
\end{align*}
where $\nu \in [-\rmi \frac{k_*}{2}, \rmi \frac{k_*}{2})$ is purely imaginary. In particular, the spectrum of $\mathcal{L}_-$ is purely essential spectrum, consisting of the union of the spectra of $\hat{\mathcal{L}}_- (\nu)$. The following conditions encode the typical picture for spectral stability of periodic patterns in extended domains. By translation invariance, $\partial_\xi u_-$ is in the kernel of $\hat{\mathcal{L}}(0)$, and so 0 always belongs to the essential spectrum of $\mathcal{L}_-$. So, the ``most stable'' scenario one can hope for is as follows.

\begin{definition}[Diffusive spectral stability of periodic patterns]\label{d:diffst}
    We say $u_-(\xi)$ is \emph{diffusively stable}
    \begin{enumerate}
        \item if we have $\mathrm{spec}\,(\mathcal{L}_-) \subset \{ \lambda \in \mathbb{C} : \Re \lambda < 0 \} \cup \{ 0 \}$;
        \item and if $\lambda = 0$ is an algebraically simple eigenvalue of $\hat{\mathcal{L}}_- (0)$ which continues to an eigenvalue  $\lambda(\nu)=-c_\mathrm{g} \nu +D_\mathrm{eff}\nu^2+\rmO(\nu^3)$ for some $D_\mathrm{eff}>0$ for $\hat{\mathcal{L}}_- (\nu)$.
    \end{enumerate}
\end{definition}
Here, $c_\mathrm{g}$ is the group velocity of the wave train, in the frame moving with the front speed, and $D_\mathrm{eff}$ the effective diffusivity. An analytic expansion follows from  the implicit function theorem; compare the green curves of spectra in the right two panels of Fig.~\ref{fig:rigid spectra}. Diffusive spectral stability guarantees nonlinear asymptotic stability of the periodic pattern; see e.g.~\cite{Schn1996, jnrz2, jnrz1, zumbruninventiones, uecker}.

\subsubsection{Rigid pattern-forming fronts.}
We investigate fronts that are stationary in a comoving frame but leave a periodic pattern in their wake, characterized much as in Defs.~\ref{def: rigid pushed} or~\ref{def: rigid pulled}, but with weaker adapted spectral stability in the wake.

\textbf{Pushed rigid pattern-forming fronts.} The typical spectral picture is as follows.
\begin{definition}[Marginal spectral stability --- pushed rigid pattern-forming fronts]\label{def: rigid pf pushed}
    A stationary solution to~\eqref{e: pf proofs eqn} is a \emph{rigid pushed pattern-forming front} if it is marginally spectrally stable and leaves behind a diffusively  stable periodic pattern according to Def.~\ref{d:diffst} with negative group velocity. It satisfies conditions (i)-(iii) of Def.~\ref{def: rigid pushed}, with the following modifications:
    \begin{itemize}
        \item the selected state in the wake, $u_-(\xi)$ is now periodic in $\xi$ rather than spatially uniform, and is diffusively stable as in Def.~\ref{d:diffst}, with negative group velocity;
        \item in condition (iii), we instead consider the operator $\mathcal{L}_\mathrm{ps}$ on the space $L^2_{\tilde{\eta}, \eta_\mathrm{ps}-\tilde{\eta}}$. That is, we assume that for any $\tilde{\eta} > 0$ sufficiently small, the operator
        \begin{align*}
            \mathcal{L}_\mathrm{ps} : H^{2m}_{\tilde{\eta}, \eta_\mathrm{ps}-\tilde{\eta}} \subset L^2_{\tilde{\eta}, \eta_\mathrm{ps}-\tilde{\eta}} \to L^2_{\tilde{\eta}, \eta_\mathrm{ps}-\tilde{\eta}}
        \end{align*}
        has a single geometrically and algebraically simple eigenvalue at $\lambda = 0$, with eigenfunction $\partial_\xi u_\mathrm{ps}$, and its spectrum is otherwise stable.
    \end{itemize}
\end{definition}
Notice that in contrast to Def.~\ref{def: rigid pushed} for unpatterned fronts, we use a small exponential weight on the left in condition (iii), in order to stabilize the neutral branch of spectrum arising from the periodic pattern. The space $L^2_{\mathrm{exp}, \tilde{\eta}, \eta_\mathrm{ps}-\tilde{\eta}}$ requires functions to decay exponentially with rate at least $\rme^{-(\eta_\mathrm{ps} - \eta) \xi}$ as $\xi \to \infty$, and allows for slow exponential growth with rate $\rme^{-\tilde{\eta} \xi}$ as $\xi \to - \infty$. The fact that this weight, rather than one that requires exponential decay as $\xi \to -\infty$, stabilizes the spectrum encodes that the group velocity of the periodic pattern points to the left, away from the front interface. That is, the front acts as a source of patterns. However, using this weight in the nonlinear equation for perturbations introduces coefficients to the nonlinearity which grow exponentially in space, and so this weight is incompatible with a nonlinear argument. By contrast, the localizing weight on the right introduces coefficients which decay exponentially as $\xi \to \infty$, and so only simplifies the nonlinear argument.

On the other hand, the sign of group velocities and thereby weights is inherent to observed patterns: fronts followed by patterns with positive group velocity, pointed towards the front interface, would typically not be observed; compare~\cite{SS2004}.
Again, as in \S\ref{s:nlmsp}, we also note that this definition is stronger than is strictly necessary, encoding spectral stability in a fixed weighted space.

\begin{figure}
    \centering
    \includegraphics[width=0.3\linewidth]{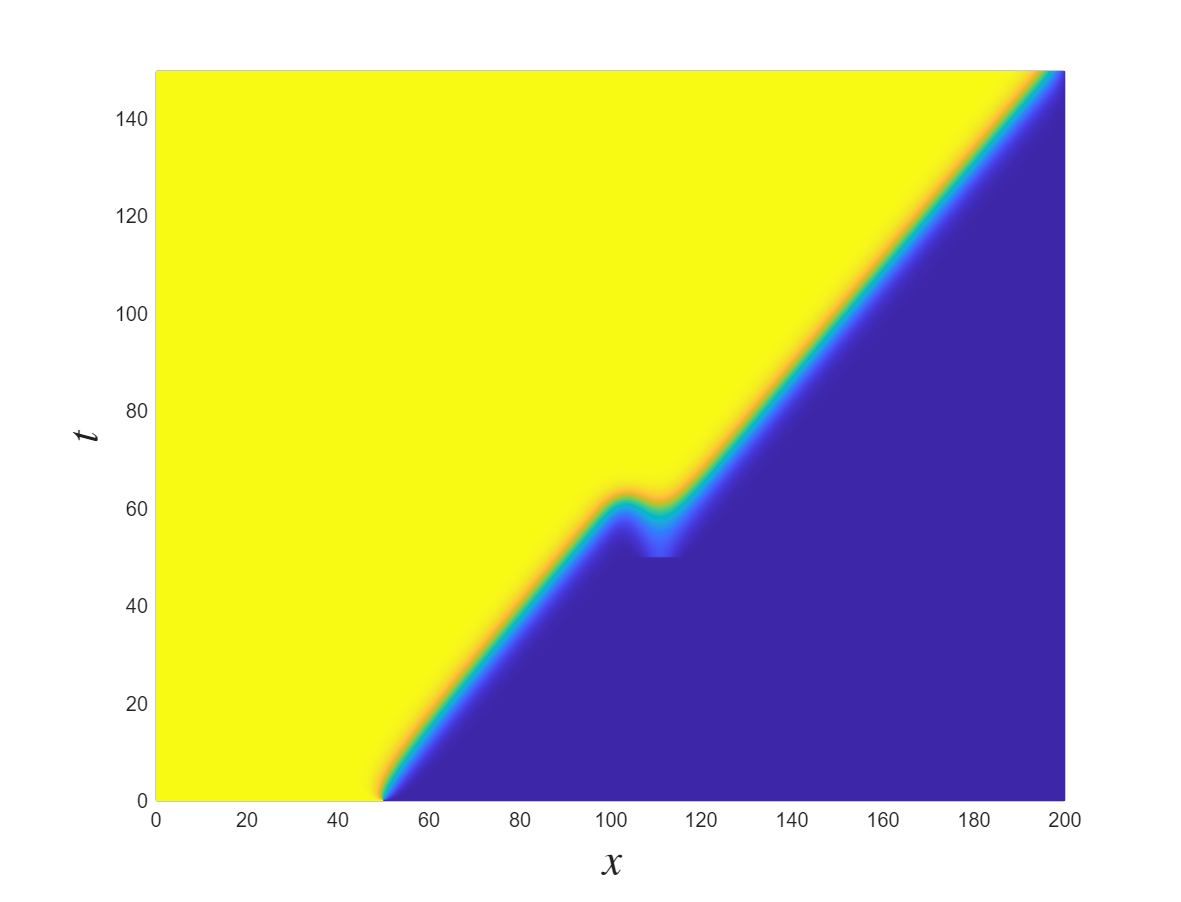}
    \includegraphics[width=0.3\linewidth]{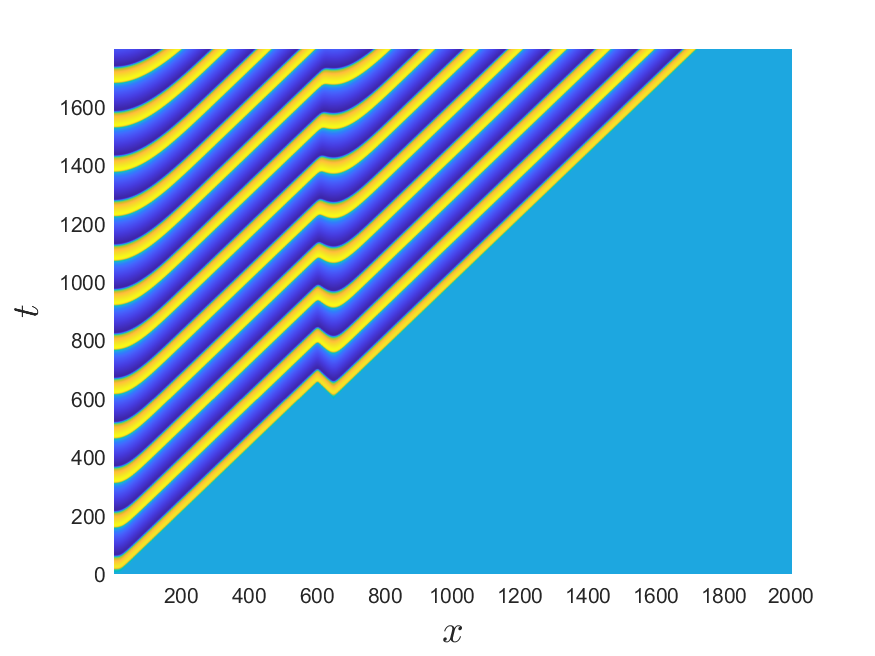}
    \includegraphics[width=0.35\linewidth]{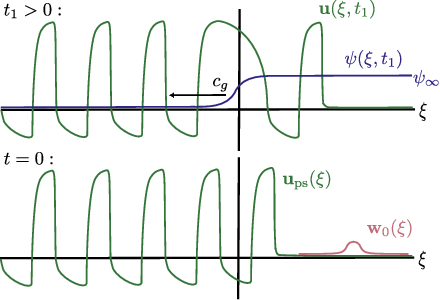}
    \caption{Left: spacetime plot depicting convergence to a pushed front in simulations of the Nagumo equation. The initial condition is a step function, which converges to a pushed front. After a moderate time, the simulation is paused and a perturbation is added ahead of the front interface. The solution then rapidly converges, uniformly in $x$, to a phase-shifted front. Center: the same experiment for a pushed pattern-forming front in the FitzHugh-Nagumo equation. The adjustment of the position of the front interface creates a phase defect in the pattern in the wake which persists for all time. Right: schematics of these dynamics in the frame co-moving with the front speed, in which the phase defect propagates to the left at the group velocity of the wave train.}
    \label{fig:rigid pushed convergence}
\end{figure}

Recall from \S\ref{s:nlmsp} that, after a perturbation which cuts off the front tail, an unpatterned pushed front rapidly adjusts its position to account for the perturbation, resulting in exponential in time convergence to a front profile with a global phase shift. When the front is pattern-forming, the front interface still rapidly adjusts its position in response to a perturbation, but this adjustment creates a defect in the phase of the pattern of the wake, and one must analyze the \emph{diffusive mixing} of the new phase created by the adjusted interface and the original phase of the front; see Fig.~\ref{fig:rigid pushed convergence}. This analysis has recently been completed in~\cite{PushedFHN}.
\begin{theorem}[\cite{PushedFHN} --- selection of rigid pushed pattern-forming fronts]\label{t: rigid pushed pattern forming selection}
    Let $u_\mathrm{ps}$ be a rigid pushed pattern-forming front.
    Then $u_\mathrm{ps}$ is a selected front in the sense of Def.~\ref{def: selected front}.
\end{theorem}
We note that the concept of selection is modified in the sense that closeness to the front is not uniform in the wake of the front due to diffusive mixing.

\textbf{Selection of rigid pushed pattern-forming fronts: proof sketch.} As in the unpatterned case, one considers perturbations $w = u - u_\mathrm{ps}$ of the front $u_\mathrm{ps}$, and derives a variation of constants formula satisfied by $w$, with the same form as~\eqref{e: pushed voc}. The difference is that in the unpatterned case, one has the linear exponential convergence $\rme^{\mathcal{L}_\mathrm{ps} t} = \omega_{0, \eta_0} \left(\partial_\xi u_\mathrm{ps}\right) u_\mathrm{ps} P_0 + \mathrm{O}(\rme^{-\mu t})$, while as a result of the diffusively stable spectrum in the wake, for pattern-forming fronts one merely has
\begin{align}
    \rme^{\mathcal{L}_{ps} t} = \omega_{0, \eta_0} \partial_\xi u_\mathrm{ps} s_{p} (t) + s_c(t) + \mathrm{O}(\rme^{-\mu t}), \label{e: pushed proof linear decomp}
\end{align}
where, roughly,
\begin{align*}
    [s_{p} (t) w_0](\xi, t) \approx \mathrm{erf}\left( \frac{\xi - c_g t}{\sqrt{4D_\mathrm{eff}^\mathrm{wt} t}} \right) P_0 w_0, \qquad\quad [s_c(t) w_0](\xi, t) \lesssim \frac{\rme^{-\frac{(\xi-c_g t)^2}{4 D_\mathrm{eff}^\mathrm{wt} t}}}{(1+t)^{1/2}},
\end{align*}
with $\mathrm{erf}(z) = \int_{-\infty}^z \rme^{-x^2} \, dx$, and $P_0$ is again a projection onto the translational mode along the adjoint eigenfunction~\cite{PushedFHN}. The leading order term $s_p(t)$ captures the interaction of the embedded eigenvalue with the essential spectrum. To resolve this at the nonlinear level, one studies the dynamics of the \emph{spatiotemporally} modulated perturbation $\tilde{w}(\xi, t) = \omega_{0, \eta_0}(\xi) [ u(\xi - \psi(\xi,t)) - u_\mathrm{ps}(\xi)]$. One can use the linearized decomposition~\eqref{e: pushed proof linear decomp} to write down a fixed-point system for $\tilde{w}$ and $\psi$, in which $\psi$ captures the leading order dynamics associated to $s_p(t)$ terms. Completing an iterative argument on this fixed-point system, one finds
\begin{align}
    \| \omega_{0, \eta_0} [u(\cdot, t) - u_\mathrm{ps}(\cdot + \psi(\cdot, t), t)] \|_{L^\infty} \lesssim (1+t)^{-1/2}, \text{ where }\ \psi(\xi, t) \approx \psi_\infty(u_0)  \mathrm{erf}\left( \frac{\xi - c_g t}{\sqrt{4D_\mathrm{eff}^\mathrm{wt} t}} \right)
    \label{e: pushed proof convergence}
\end{align}
for a constant $\psi_\infty(u_0) \in \R$. Convergence is in fact exponential on any half-line $[-K, \infty)$, with estimates
\begin{align*}
    \| \psi(\cdot, t) - \psi_\infty \|_{L^\infty[-K, \infty)} + \| \omega_{0, \eta_0} [u(\cdot, t) - u_\mathrm{ps} (\cdot +\psi_\infty)]\|_{L^\infty[-K, \infty)} \lesssim \rme^{-\mu t}.
\end{align*}
We note, however, that convergence of $u(\xi,t)$ to $u_\mathrm{ps} (\xi + \psi_\infty)$ is not uniform on the real line: the phase defect introduced by the perturbation widens but persists for all time, propagating away from the front interface according to the group velocity; see Fig.~\ref{fig:rigid pushed convergence}.

Crucial to closing the nonlinear argument is the estimate $\| P_0 w_0 \|_{L^\infty} \lesssim \| \omega_{\kappa, 0} w_0\|_{L^\infty}$, that is, the translational mode only very weakly responds to excitations to the left of the front interface. Essentially, the dynamics near and ahead of the front interface remain as in the unpatterned case, and the pattern formation mechanism in the wake simply responds to the phase adjustment introduced by the rapid adjustment of the front interface in response to the initial perturbation.

\textbf{Pulled rigid pattern-forming fronts.} The typical spectral picture now is as follows.
\begin{definition}[Marginal spectral stability --- pulled rigid pattern-forming fronts]
    A stationary solution to~\eqref{e: pf proofs eqn} is a \emph{rigid pulled pattern-forming front} if it is marginally spectrally stable and leaves behind a diffusively  stable periodic pattern according to Def.~\ref{d:diffst}, with negative group velocity. It satisfies conditions (i)-(v) of Def.~\ref{def: rigid pulled}, with the following modifications:
    \begin{itemize}
        \item the state $u_-(\xi)$ selected in the wake is now periodic in space, rather than spatially uniform;
        \item we modify condition (iv) to assume that for any $\tilde{\eta} > 0$ sufficiently small,
        \begin{align*}
            \mathcal{L}_\mathrm{pl} : H^{2m}_{\mathrm{exp}, \tilde{\eta}, \eta_\mathrm{lin}} \subset L^2_{\mathrm{exp}, \tilde{\eta}, \eta_\mathrm{lin}} \to L^2_{\mathrm{exp}, \tilde{\eta}, \eta_\mathrm{lin}}
        \end{align*}
        has stable spectrum apart from the  simple branch touching $\rmi\R$ at the origin due to (i).
    \end{itemize}
\end{definition}

As in Def.~\ref{def: rigid pf pushed}, the small exponential weight on the left is needed to stabilize the diffusive spectrum associated to the periodic pattern, and the fact that the stabilizing weight allows exponential growth on the left encodes that the group velocity associated to $u_-(\xi)$ is negative (in the frame co-moving with the front speed), so that the front acts as a source of patterns.

Rigorous selection of rigid pulled pattern-forming fronts remains an open problem. Some partial progress establishes an upper bound on the propagation speed for systems with a variational structure~\cite{CE1990}, or stability against perturbations which are small in $L^1_{\mathrm{exp},0, \eta_\mathrm{lin}}$~\cite{avery2023stabilitycoherentpatternformation} and hence cannot cut off the front tail in light of the asymptotics in Def.~\ref{def: rigid pulled}(iii). Based on the linear theory developed in~\cite{avery2023stabilitycoherentpatternformation}, we nonetheless strongly suspect that selection holds for rigid pulled pattern-forming fronts.

\begin{conjecture}[Selection of rigid pulled pattern-forming fronts]\label{c: rigid pf pulled}
    Let $u_\mathrm{pl}$ be a rigid pulled pattern-forming front. Then $u_\mathrm{pl}$ is a selected front in the sense of Def.~\ref{def: selected front}, with position correction $h(t) = -\frac{3}{2 \eta_\mathrm{lin}} \log t + \mathrm{O}(1)$.
\end{conjecture}

Intuitively, the challenge in establishing Conj.~\ref{c: rigid pf pulled} is that the logarithmic drift in the interface position introduces a slowly growing phase shift in the wake and thereby an ongoing disturbance to the phase in the wake. Technically, comparing to the proof of Thm.~\ref{t: rigid pulled unpatterned}, a key challenge is that the usefulness of the transformation from $w$ to $z = (t+t_0)^{-3/2} w$ in the proof of Thm.~\ref{t: rigid pulled unpatterned} relies on the fact that the linearized evolution decays everywhere with rate $t^{-3/2}$. Hence, one can use the linearized evolution to prove $(t+t_0)^{-3/2}$--decay for $z$, which translates to boundedness of $w$. However, the diffusive stability of the periodic pattern induces linear decay  only with rate $t^{-1/2}$ for $\xi \ll -1$. Hence, one cannot use a perturbative argument based on the linearization to prove that $z$ decays globally with the desired rate. Interestingly, the stability result~\cite{avery2023stabilitycoherentpatternformation} does establish that the $t^{-3/2}$ decay rate is preserved for the linearized evolution for $\xi \geq 1$, i.e. to the right of the front interface.

\subsubsection{Modulated fronts: conjectures}
Results on selection of modulated fronts, analogous to Thms.~\ref{t: pushed stability} and~\ref{t: rigid pulled unpatterned}, are not available.
In analogy to~\cite{avery2023stabilitycoherentpatternformation, PushedFHN} but relying on Floquet theory as in~\cite{SS2001}, we outline a framework and state conjectures relating  stability and selection of modulated fronts with patterns in their wake.

\begin{definition}[Marginal spectral stability --- modulated pushed fronts]\label{def: modulated pushed}
    A time-periodic solution $u_\mathrm{ps}(\xi, \tau)$  {\color{black} to} ~\eqref{e: pf proofs eqn} with period $T$ and speed $c = c_\mathrm{ps}$ is a \emph{modulated pushed pattern-forming front} if it is marginally stable and leaves behind a diffusively stable periodic pattern according to Def.~\ref{d:diffst}, with negative group velocity. Typical marginal stability of modulated pushed fronts is encoded as follows. 
    \begin{enumerate}
        \item \emph{Selected state in the wake:} Assume that there exists a selected wave train in the wake, $u_-(\xi, \tau) = u_\mathrm{wt}(k\xi - \omega \tau)$ to~\eqref{e: pf proofs eqn} with $c = c_\mathrm{ps}$ with $u_\mathrm{wt}(\zeta)$ periodic in its argument. This periodic solution is diffusively stable, as in Def.~\ref{d:diffst}, and the group velocity is negative in the frame co-moving with the front. 
        \item \emph{Exponential asymptotics:} Assume that $u_\mathrm{ps}(\xi, \tau)$ converges with exponential rate in $\xi$ to $u_-(\xi, \tau)$ as $\xi \to -\infty$ and to $0$ as $\xi \to \infty$. Let $\eta_\mathrm{ps}$ denote the exponential convergence rate as $\xi \to \infty$.
        \item \emph{Minimal marginal Floquet stability:} Let $\mathcal{L}_\mathrm{ps} = \mathcal{P}(\partial_\xi) + c_\mathrm{ps} \partial_\xi + f'(u_\mathrm{ps}(\xi, \tau))$ denote the linearization of~\eqref{e: pf proofs eqn} about $u_\mathrm{ps}(\xi, \tau)$. Let $Y = L^2(\mathbb{R}/T\mathbb{Z})$, and let $Y^1 = H^1(\mathbb{R}/ T\mathbb{Z})$. We assume that for any $\tilde{\eta}$ sufficiently small, the operator
        \begin{align*}
            -\partial_\tau + \mathcal{L}_\mathrm{ps}: H^{2m}_{\mathrm{exp}, \tilde{\eta}, \eta_\mathrm{ps}-\tilde{\eta}} (\mathbb{R}, Y) \cap L^2_{\mathrm{exp}, \tilde{\eta}, \eta_\mathrm{ps}-\tilde{\eta}} (\mathbb{R}, Y^1) \subset L^2_{\mathrm{exp}, \tilde{\eta}, \eta_\mathrm{ps}-\tilde{\eta}} (\mathbb{R}, Y) \to L^2_{\mathrm{exp}, \tilde{\eta}, \eta_\mathrm{ps}-\tilde{\eta}} (\mathbb{R},Y)
        \end{align*}
        has a semi-simple eigenvalue at $\lambda = 0$ of multiplicity 2, with eigenfunctions $\partial_\xi u_\mathrm{ps}$ and $\partial_\tau u_\mathrm{ps}$, and its spectrum is otherwise stable in $0\leq \Im\lambda<2\pi/T$.
    \end{enumerate}
\end{definition}
Recall from \S\ref{s:5.2} that eigenvalues of $-\partial_\tau + \mathcal{L}_\mathrm{ps}$ correspond to Floquet exponents for the periodic map, and so condition (iii) encodes marginal stability of the linearized dynamics, arising from point spectrum as expected for pushed fronts.
We expect that the theory developed in~\cite{SS2001} can be used to transfer the analysis of the linearized dynamics for rigid pattern-forming fronts in~\cite{PushedFHN} to the modulated setting (see \cite{BeckSandstedeZumbrun} for a related analysis for time-periodic viscous shocks), which would thereby establish the following general selection result for pushed modulated fronts.
\begin{conjecture}
    Let $u_\mathrm{ps}(\xi, \tau)$ be a modulated pushed pattern-forming front. Then $u_\mathrm{ps}(\xi, \tau)$ is a selected front in the sense of Def.~\ref{def: selected front}.
\end{conjecture}

We can similarly formulate definitions and conjectures regarding typical pulled modulated fronts.
\begin{definition}[Marginal spectral stability --- pulled modulated fronts]\label{def: modulated pulled}
    A time-periodic solution  $u_\mathrm{pl}(\xi, \tau)$ to~\eqref{e: pf proofs eqn} with {\color{black} period $T$ and speed} $c = c_\mathrm{lin}$ is a \emph{modulated pulled pattern-forming front} if it is marginally spectrally stable and leaves behind a diffusively  stable periodic pattern according to Def.~\ref{d:diffst} with negative group velocity.
    Typical marginal stability of modulated pulled fronts is encoded as follows.
    \begin{enumerate}
        \item \emph{Linear spreading speed via pinched double roots:} Suppose that, with $c = c_\mathrm{lin}$, the dispersion relation for linearization of~\eqref{e: pf proofs eqn} about $u = 0$ has simple pinched double roots at $\lambda = \pm i \omega_\mathrm{lin}$, with spatial eigenvalue $\nu = \nu_\mathrm{r} + \rmi \nu_\mathrm{i}$, $\nu_\mathrm{r} < 0$, where $\omega_\mathrm{lin} = \frac{2 \pi}{T}$. Let $\eta_\mathrm{lin} = -\nu_r$.
        \item \emph{Selected state in the wake:} Assume that there exists a selected wave train solution to~\eqref{e: pf proofs eqn}, $u_-(\xi, \tau) = u_\mathrm{wt}(k \xi - \omega \tau)$, with $c = c_\mathrm{lin}$, with $u_\mathrm{wt}(\zeta)$ periodic in $\zeta$. This periodic solution is diffusively stable, as in Def.~\ref{d:diffst}, and the group velocity is negative in the frame co-moving with the front. 
        \item \emph{Exponential asymptotics:} Assume that $u_\mathrm{pl}(\xi, \tau)$ converges with exponential rate in $\xi$ to $u_-(\xi, \tau)$ as $\xi \to -\infty$, and that the leading order asymptotics as $\xi \to \infty$ are $u_\mathrm{pl}(\xi, \tau) \sim A \xi \rme^{\nu \xi + \rmi \omega_\mathrm{lin} t} + \bar{A} \xi \rme^{\bar{\nu} \xi - \rmi \omega_\mathrm{lin} t}$, with $A \in \mathbb{C}^N$ nonzero.
        \item \emph{Minimal marginal Floquet stability:} Let $\mathcal{L}_\mathrm{pl} = \mathcal{P}(\partial_\xi)+c_\mathrm{lin}(\partial_\xi)+f'(u_\mathrm{pl}(\xi,\tau))$ denote the linearization of~\eqref{e: pf proofs eqn} about $u_\mathrm{pl}(\xi, \tau)$. Let $Y = L^2(\mathbb{R}/T\mathbb{Z})$ and $Y^1 = H^1(\mathbb{R} / T \mathbb{Z})$. We assume that for any $\tilde{\eta} > 0$ sufficiently small, the operator
        \begin{align*}
            -\partial_\tau + \mathcal{L}_\mathrm{pl} : H^{2m}_{\mathrm{exp}, \tilde{\eta}, \eta_\mathrm{lin}} (\mathbb{R}, Y) \cap L^2_{\mathrm{exp}, \tilde{\eta}, \eta_\mathrm{lin}} (\mathbb{R}, Y^1) \subset L^2_{\mathrm{exp}, \tilde{\eta}, \eta_\mathrm{lin}}(\mathbb{R}, Y) \to L^2_{\mathrm{exp}, \tilde{\eta}, \eta_\mathrm{lin}}(\mathbb{R}, Y)
        \end{align*}
        has two branches of essential spectrum touching the origin, corresponding to the simple pinched double roots from (i), and is otherwise stable.
        \item \emph{No resonance at $\lambda = 0$:} Assume that there are no solutions $u \in L^\infty_{\mathrm{exp}, 0, \eta_\mathrm{lin}}(\mathbb{R}, Y)$ to the equation $[-\partial_\tau + \mathcal{L}_\mathrm{ps}]u = 0$.
    \end{enumerate}
\end{definition}
\begin{conjecture}\label{c: modulated pulled selection}
    Let $u_\mathrm{pl}(\xi, \tau)$ be a modulated pulled pattern-forming front. Then $u_\mathrm{pl}(\xi, \tau)$ is a selected front in the sense of Def.~\ref{def: selected front}.
\end{conjecture}
To prove Conj.~\ref{c: modulated pulled selection}, we expect difficulties similar to Conj.~\ref{c: rigid pf pulled}. Conj.~\ref{c: modulated pulled selection} was formally resolved using matched asymptotics in~\cite{EBERT200413}. Selection throughout is defined in spaces that allow for growth in the wake due to the potential presence of non-decaying wave train dynamics.

One can envision other scenarios for modulated invasion and associated conjectures: any one of the leading edge, the front interface itself, or the state in the wake can  be oscillatory or temporally constant. The invasion will be oscillatory if any one of the three exhibits oscillations and we conjecture that all cases occur in parametrically forced Hopf bifurcation; compare for instance Fig.~\ref{f:ressliver} for some scenarios.

\subsection{Periodically forced problems,  problems in lattices, and problems in cylinders}\label{s:infdim}
The difficulty  of temporal periodicity appears in many other contexts, for instance when the equation is periodically forced in time or in space with periods $T$ or $L$, respectively,
\begin{equation}\label{e:parper}
u_t=\mathcal{P}(\partial_x,t,x)u+f, \qquad f=f(u,t,x)=f(u,t+T,x+L).
\end{equation}
Assuming, after shifting $u$ if necessary, that $u\equiv 0$ is an  equilibrium, we find a linearized equation with periodic coefficients, amenable to spatio-temporal Floquet analysis, thus expecting a linear invasion speed $c_\mathrm{lin}$ and frequency $\omega_\mathrm{lin}$. One then proceeds to identify temporally resonant nonlinear states in the wake. The situation here is simplified since the linearization at such periodic patterns would typically not possess a zero neutral eigenvalue and allow for exponential stability in the wake of the invasion. Spatial periodicity can also be due to spatial discreteness, as in equations on lattices,
\begin{equation}\label{e:lattice}
u_{j,t}=(\mathcal{P}(\delta_-,\delta_+)u)_j+f(u_j),\qquad j\in\Z,
\end{equation}
with $\mathcal{P}$ polynomial, $(\delta_-u)_j=u_j-u_{j-1}$, $(\delta_+u)_j=u_{j+1}-u_{j}$. Fronts of the form $u(j-ct)$ with temporal period $1/c$ solve a forward-backward delay equation with ill-posed initial-value problem similar to~\eqref{e:mtweq}~\cite{zinner,mallet-paret}. The methods discussed in this review could also be applied in this context to study selection of fronts; see \cite{BFRZ} for logarithmic corrections and \cite{wang2019pinned,holzer2022locked} for invasion processes where novel phenomena arises due to a combination of both temporal and spatial discreteness.

Technically similar are also problems in cylinders~\cite{heinze,berestyckinirenberg},
\begin{equation}\label{e:parcyl}
u_t=\mathcal{P}(\partial_{x_1},\nabla_{x_\perp})+f(u),\quad x=(x_1,x_\perp)\in \R\times \Omega, \qquad
\text{``+'' b.c. on } \R\times\partial\Omega.\end{equation}
The elliptic traveling-wave equation in the strip $\R\times\Omega$  accommodates a similar point of view as~\eqref{e:mtweq}. We discuss the case of  $\Omega=\T^{n-1}$,  the torus, that is, periodic boundary conditions in $x_\perp$, below.

Linear spreading speeds  are found following the recipes of analyzing pointwise resolvents through (now infinite-dimensional) stable and unstable subspaces, and identifying pinching points where stable and unstable roots collide.

\subsection{Transverse pattern selection}\label{s:5.5}

For $x\in\R^n$, $n>1$,  many new questions arise. In an isotropic system, the linear equation predicts spreading with ball-shaped level sets growing with the linear spreading speed. In the leading edge, dynamics resemble dynamics along a planar interface with small curvature corrections; see~\cite{rrr} for the FKPP case. Considering on the other hand initial conditions localized in a strip $x_1\in (-L,L)$, $x_\perp\in\R^n$, we may decompose into Fourier modes $\rme^{\rmi k_\perp\cdot x_\perp}$ and ask for spreading speeds $c_\mathrm{lin}(k_\perp)$ for each of these Fourier modes.  In an isotropic system, surprisingly this spreading speed is always maximized for $k_\perp=0$; see~\cite[\S7]{HolzerScheelPointwiseGrowth}. For pattern-forming systems, the leading edge therefore forms patterns that are independent of $x_\perp$, that is,  stripes perpendicular to the direction of spreading. Hexagonal or rectangular patterns formed in the wake of a pulled invasion front should then be interpreted as a failure of node conservation, similar to cases where the frequency of the leading edge does not predict the wavenumber correctly. Trying to understand the emergence of different patterns in the wake of fronts then naturally leads to  discussing secondary invasion fronts; see \S\ref{s:stage}.

\subsection{Practical considerations: modulated fronts}\label{s:modpract}

Many of the practical issues discussed in \S\ref{s:nlms} and \S\ref{s:robust} can in theory be immediately generalized to the case of modulated, time-periodic invasion. Very little of this has however actually been done. We mention~\cite{GLR2024} where pulled fronts have been computed in the Swift-Hohenberg equation using a Newton method for~\eqref{e:mtweq} with suitable phase conditions and an additional farfield-core decomposition in the wake. We believe that the picture emerging from analysis and computations will be much richer than described here. Hopefully, progress will emerge through a better understanding of examples, for instance establishing assumptions for marginally stable pattern-forming fronts in examples and, of course, establishing selection results.

\textbf{Wavenumber selection.} We emphasize that as a general recipe, 1:1 resonant selection is very successful, and restate here the basic strategy to find patterns selected in the wake, if the invasion is pulled:
\begin{itemize}
    \item[(i)] determine the nonlinear dispersion relation for wave trains $u_\mathrm{wt}(kx-\omega t)$ in the steady frame $\omega=\omega(k)$; see \S\ref{s:5.3};
    \item[(ii)] determine the spreading speed and linear invasion frequency $\omega_\mathrm{lin}$ in the comoving frame;
    \item[(iii)] determine the selected wavenumber $k$ solving \eqref{e:reswake}
    \[
    \omega_\mathrm{lin}=\omega_\mathrm{nl}(k)-c_\mathrm{lin} k.
    \]
\end{itemize}

\subsection{Modulated fronts: examples}\label{s:5:exples}

\textbf{The complex Ginzburg-Landau equation.} The modulation equation~\eqref{e:cgl},
\begin{equation}\label{e:cgl_a}
A_t=A_{xx}+A-A|A|^2,\qquad A\in\C,
\end{equation}
arises near the onset of instability, where a physical variable is expressed as $u(t,x)=A(\eps^2 t,\eps x)\rme^{\rmi k x}e_0+c.c.$ for some eigenmode $e_0$. The linear spreading speed is  $c_\mathrm{lin}=2$ and fronts are of the form $\rme^{\rmi\psi}u(\xi)$, $\psi\in[0,2\pi)$, and $u''+2u'+u-u^3=0$ is the Nagumo front. Selection of this front is not known but stability of the critical front was established in~\cite{EW1994} and in~\cite{as5} with sharp decay estimates. One should think of these fronts having an averaged out oscillation, which could be reintroduced by the transformation $B=\rme^{\rmi\omega t}A$, creating a pair of pinched double roots at $\pm\rmi\omega$ at $c_\mathrm{lin}=2$.
%
Variations near weakly subcritical bifurcations include the cubic-quintic version
\[
A_t=A_{xx}+A+\alpha A|A|^2-A|A|^4,\qquad A\in\C,
\]
with pulled fronts for $\alpha<2/\sqrt{3}$ and pushed fronts for $\alpha>2/\sqrt{3} $ with speed $\frac{1}{\sqrt{3}}(-a+2\sqrt{4+a^2})$.

\textbf{The Swift-Hohenberg equation.}
The simplest case where the complex Ginzburg-Landau equation actually describes the dynamics of small-amplitude solutions is the Swift-Hohenberg equation~\eqref{e:sh},
\begin{equation}\label{e:shp}
    u_t=-(\partial_{xx} +1 )^2 + \eps^2 u - u^3,
\end{equation}
with $\eps\gtrsim 0$ and amplitude approximation $u(t,x)=\eps A(\eps^2 t,\eps x)+c.c.$, $A_T=4A_{XX}+A-3A|A|^2$, with spreading speed $4$, that is, spreading speed $4\eps$ in~\eqref{e:shp}. In fact, the spreading speed is explicit from~\eqref{e:4cth4} with $a=-2,b=-1+\eps^2$, with expansion
\begin{equation}\label{e:shspeed}
c_\mathrm{lin}=4\eps+\eps^3-\frac{9}{8}\eps^5+\rmO(\eps^7),\
\omega_\mathrm{lin}=4\eps+\frac{3}{2}\eps^3-\frac{45}{32}\eps^5+\rmO(\eps^7),\  k=\omega_\mathrm{lin}/c_\mathrm{lin}=1+\frac{1}{8}\eps^2-\frac{13}{128}\eps^4+\rmO(\eps^6).
\end{equation}
Existence and linear stability of  fronts with $c>c_\mathrm{lin}$ has been established in a series of papers exploiting the approximation by~\eqref{e:cgl_a}; see~\cite{CE1986,CE1987,CE1990,EW1991,es1}.
Existence can be understood as a center-manifold reduction in the modulated traveling-wave equation with frequency and speed of order $\eps$,
\[
\eps\omega_1 u_\tau=-(\partial_{\xi\xi}+1)^2u+\eps^2 u + c_1\eps u_\xi.
\]
Decomposing into Fourier modes, $u=\sum\hat{u}\rme^{\rmi\ell\tau}$, we find spatial eigenvalues as roots $\nu_\ell$ of
\[
\eps\omega_1\rmi\ell=-(\nu^2+1)^2+\eps^2+c_1\eps\nu.
\]
At $\eps=0$, $\nu_\ell=\pm\rmi$ for all $\ell$. Expanding in $\eps$, and substituting the $\omega_1$ and $c_1$ from~\eqref{e:shspeed}, we find
\[
\nu=\rmi-2\eps+\rmO(\eps^2), \text{ when } \ell=1, \qquad \nu\sim \pm\sqrt{\eps\rmi(\pm 1-\ell)}, \text{ otherwise}.
\]
The gap in the real part of the spectrum between $\nu={\pm 1}+\rmO(\eps)$ and the comparatively large $\rmO(\sqrt{\eps})$ allows one to reduce to a 4-dimensional center-manifold parameterized by the modes $|\ell|=1$ and recover at leading order the Ginzburg-Landau traveling-wave equation~\cite{EW1991}. Marginal stability is not known.







\textbf{Other pattern-forming systems.}
The Swift-Hohenberg equation was designed to mimic instabilities in fluids~\cite{SH1977},  appeared in a slightly modified form much earlier in Turing's notes on morphogenesis~\cite{DAWES201649}, and has been used across the sciences to model the emergence of spatial patterns. Experimental examples where quantitative comparisons to theory are feasible include fluid experiments, such as the motivating example of Rayleigh-B\'enard convection in~\cite{SH1977} and the Taylor-Couette experiment. Mimicking the ideas for the analysis of fronts in the Swift-Hohenberg equation, existence and stability of fronts with $c>c_\mathrm{lin}$ was established for the Navier-Stokes equations modeling the Taylor-Couette experiment~\cite{HSchn1999,ES2000}. In a reaction-diffusion pattern-forming context, the technique was adapted to construct pattern-forming fronts in     FKPP with non-local reaction terms in~\cite{FH2015}.

\textbf{Hopf bifurcations.}
Bifurcations discussed thus far form spatially periodic patterns. Other instabilities may generate patterns periodic in time, or in both space and time; see~\cite{crosshohenberg,SRadial} for some background. In case of a Hopf bifurcation, the modulation equation is the complex coefficient Ginzburg-Landau equation~\eqref{e:ccgl}. We discussed linear invasion speeds in \S\ref{s:2ex}. Existence, stability, or selection of invasion fronts are not known. One can however find pulled invasion fronts of the form $A=\rme^{\rmi\Omega t} B(x-ct)$, with $\Omega=\alpha$ and $c=2\sqrt{1+\alpha^2}$, as  equilibria of
\begin{equation}\label{e:cglfr}
   B_t= (1+\rmi\alpha)B_{\xi\xi}+cB_\xi + (1-\rmi\alpha)B-(1+\rmi\beta)B|B|^2 \stackrel{!}{=}0.
\end{equation}
For the resulting ODE in $(B,B_\xi)\in\C^2\sim \R^4$, the origin is stable with algebraically quadruple, geometrically double eigenvalue $\nu=-(1-\rmi\alpha)/(1+\alpha^2)$, reflecting that we chose $c$ and $\Omega$ according to the marginal stability condition. The ODE also possesses spatially periodic solutions of the form $B_\mathrm{wt}(k)=r(k)\rme^{\rmi k\xi}$, $r(k)^2=1-k^2$,
\begin{equation}\label{e:cglks}
k_\mathrm{s}^\pm = \left(\sqrt{1+\alpha^2}\pm\sqrt{1+\beta^2}\right)/\left(\alpha-\beta\right).
\end{equation}
The dispersion relation in the steady frame $\omega_\mathrm{nl}(k)=\beta+(\alpha-\beta)k^2$ gives $c_\mathrm{g}(k)=2(\alpha-\beta)k$ so that in the frame with speed $c_\mathrm{lin}$, precisely the solution with $k=k_\mathrm{s}^-$ has negative group velocity. A short calculation gives that the associated periodic solution in~\eqref{e:cglfr} has a two-dimensional unstable manifold. Invasion fronts exist precisely when solutions in the unstable manifold converge to the asymptotically stable origin, a fact that is easy to verify numerically.

More information can be obtained when linear and nonlinear dispersion are balanced, $|\alpha-\beta|\ll 1$. At $\alpha=\beta$, we set $B=\rme^{\rmi k \xi}b$ with $k=\alpha/\sqrt{1+\alpha^2}$ and find real fronts in
\[
b''+\frac{2}{\sqrt{1+\alpha^2}}b'+\frac{1}{1+\alpha^2} b -b|b|^2.
\]
A perturbative analysis gives existence for $|\alpha-\beta|\ll 1$ and possibly information on marginal stability.

Front dynamics can however be rather complex since the pattern in the wake can be unstable, leading to secondary invasion that we shall discuss in \S\ref{s:stage}; see~\cite{VANSAARLOOS1992303} and~\cite{NozBek1983,sherratt2009locating}. In the cubic-quintic Ginzburg-Landau, modeling weakly subcritical bifurcations, one can find pushed fronts, sometimes even explicitly~\cite{VANSAARLOOS1992303}; see also
\cite{aransonkramer} for more context on the Ginzburg-Landau equation.
Beyond modulation equations, one can recover the traveling-wave equation on a center-manifold for the ill-posed modulated-wave equation~\eqref{e:mtweq}; see~\cite{SRadial}. In this case, the temporal frequency is close to the Hopf frequency and therefore not small, simplifying the reduction process.  Finally, little seems to be known about fronts in the case of instabilities that create wave trains $u(\omega t - kx)$ with $\omega,k=\rmO(1)$ near onset.

\textbf{Phase separation.} Simple models include the Allen-Cahn and Cahn-Hilliard equations, phase-field systems, and even variants of the Keller-Segel model. One is usually concerned with the formation of clusters and subsequent dynamics of clusters through motion, merging, or parasitic growth. Periodic arrangements of clusters are known to exist in all the models mentioned, via a reduction to simple scalar nonlinear pendulum equations, yet they are typically unstable. Spatially constant states can be unstable against what is usually referred to as spinodal decomposition, or cluster formation, most easily understood in the Cahn-Hilliard equation with linearization $u_t=-(u_{xx}+(1-3m^2)u)_{xx}$. One can find linear spreading speeds and frequencies in this equation and establish existence of invasion fronts that leave behind periodic patterns~\cite{S2017}. One can also find fronts that invade the unstable periodic patterns and leave behind new periodic patterns with larger wavelength, thus mediating a coarsening process~\cite{S2006,S2017}. The existence proofs rely on Conley index arguments, specifically connection matrices, which rely on topological forcing of heteroclinic connections based on a priori bounds and detailed knowledge of invariant sets and their Morse indices~\cite{Conley1978,FiMi1992}. The dynamical systems tools are applied to a finite-dimensional approximation with truncated Fourier series in time $\tau$ for~\eqref{e:mtweq}, together with a priori bounds when passing to the limit of infinite series.




Numerical studies indicate that invasion is always pulled~\cite[Fig.~5\&7]{S2017}, but that the linearly predicted and the observed wavenumber differ for weak instabilities. Similar phenomena were observed in the Keller-Segel model of chemotaxis
\cite{bose2013invasion} and in a phase-field model
\cite{GMS1,kotza2012}. The change in wavenumber can be understood in the framework of \S\ref{s:5.3} as changes in harmonics in~\eqref{e:reswake}; see
~\cite[Fig.~8]{S2017} 
  for a numerical exploration indicating preferences for $p=1$, $q=1,2,3,\ldots$ in~\eqref{e:reswake}.



\begin{figure}
    \includegraphics[height=1.5in]{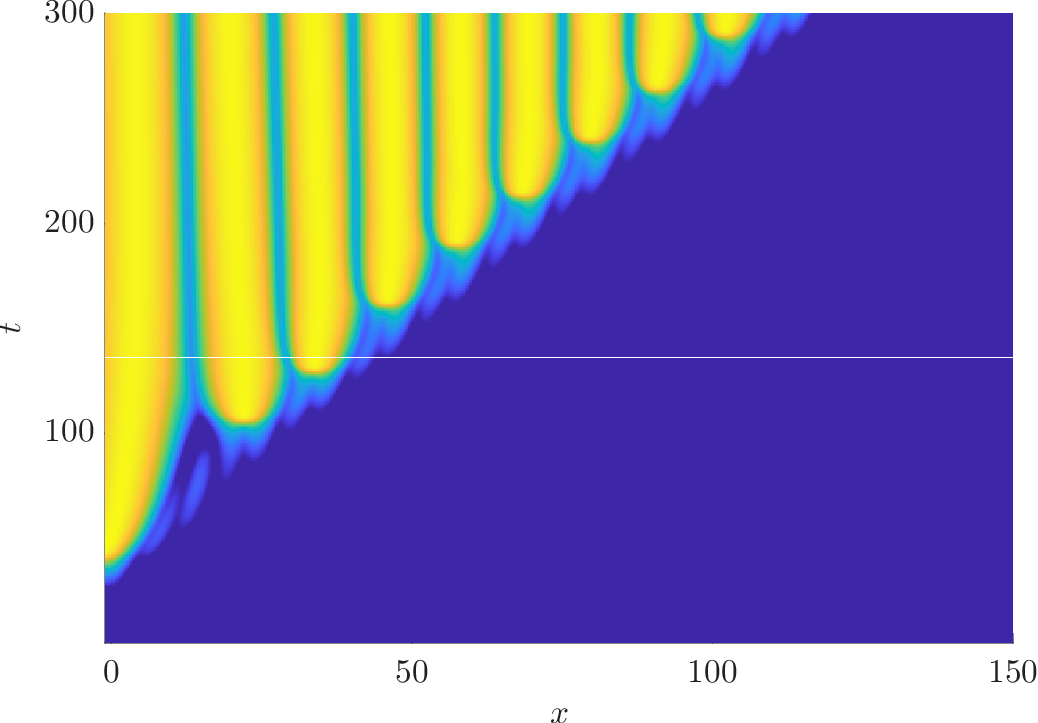}\hfill
    \includegraphics[height=1.5in]{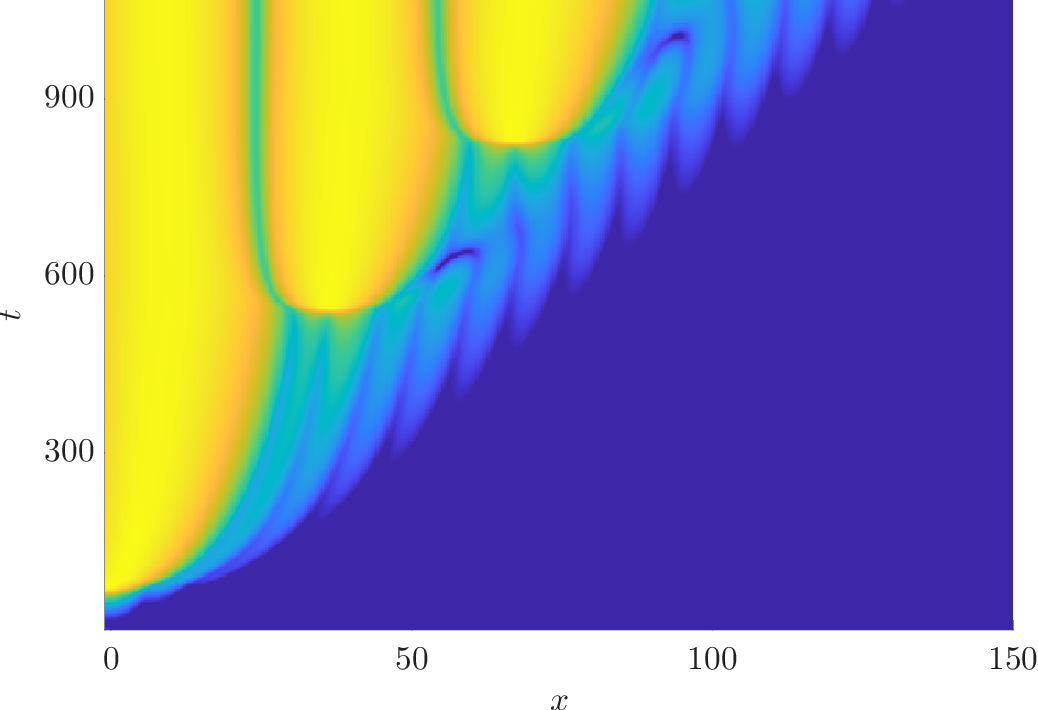}\hfill
    \includegraphics[height=1.5in]{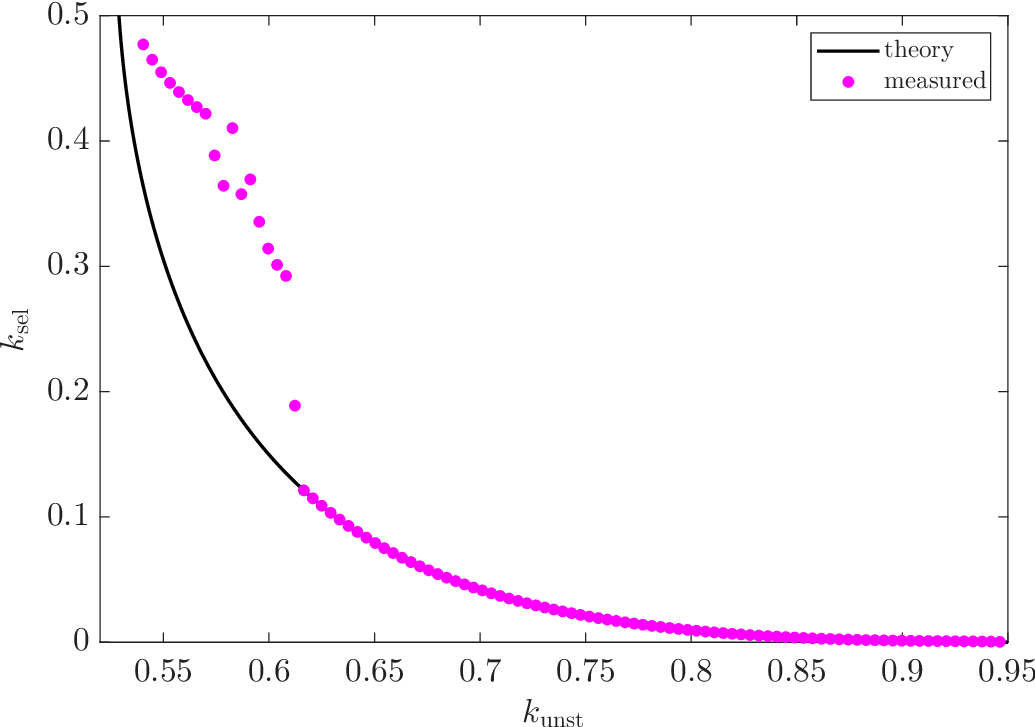}
    \caption{Dynamics of fronts invading $A=r\rme^{\rmi k_\mathrm{unst}x}$, $r^2=1=k_\mathrm{unst}^2$, in CGL~\eqref{e:cgl}. Shown are the perturbation with $k_\mathrm{unst}=0.67498$ (left) and $k_\mathrm{unst}=0.61158$ (center), and the selected wavenumber in the wake versus linear prediction (right). While speeds (not shown) agree with the linear theory for all $k_\mathrm{unst}$, wavenumbers differ due to an intricate interaction between small oscillations in the leading edge and large nonlinear structures (center).}
    \label{f:eck}
\end{figure}

\textbf{Sideband instabilities.} Instability of periodic patterns is often due to sideband instabilities~\cite{crosshohenberg}, which can at small amplitudes be described by Cahn-Hilliard equations~\cite{vanHarten1995}. The resulting invasion process has been observed in many contexts, including dispersive Hamiltonian equations~\cite{PhysRevLett.116.043902,mossman2025nonlinearstagemodulationalinstability} and studied theoretically to some extent in the Ginzburg-Landau equation~\cite{EckGal1993}. Again, fronts, appear to be pulled and one observes higher resonances for weak instabilities; compare Fig.~\ref{f:eck}.

\section{Discussion and more open problems}\label{s:dis}

\subsection{Instabilities of fronts and  wakes, and staged invasion}\label{s:stage}

\begin{figure}
\centering
 \includegraphics[width=0.36\textwidth]{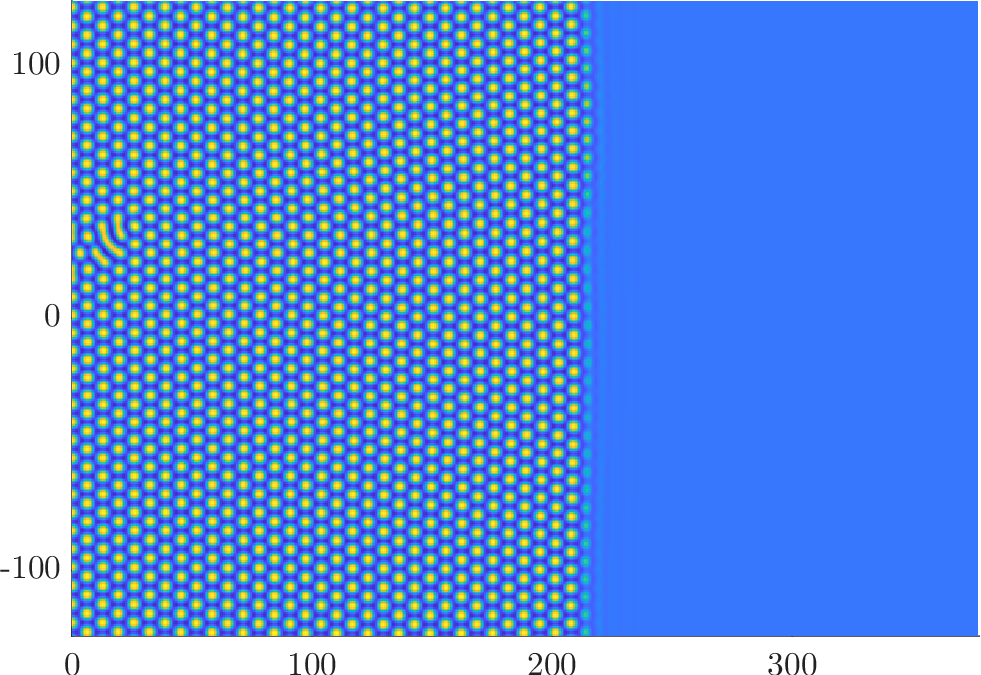}\qquad
 \includegraphics[width=0.36\textwidth]{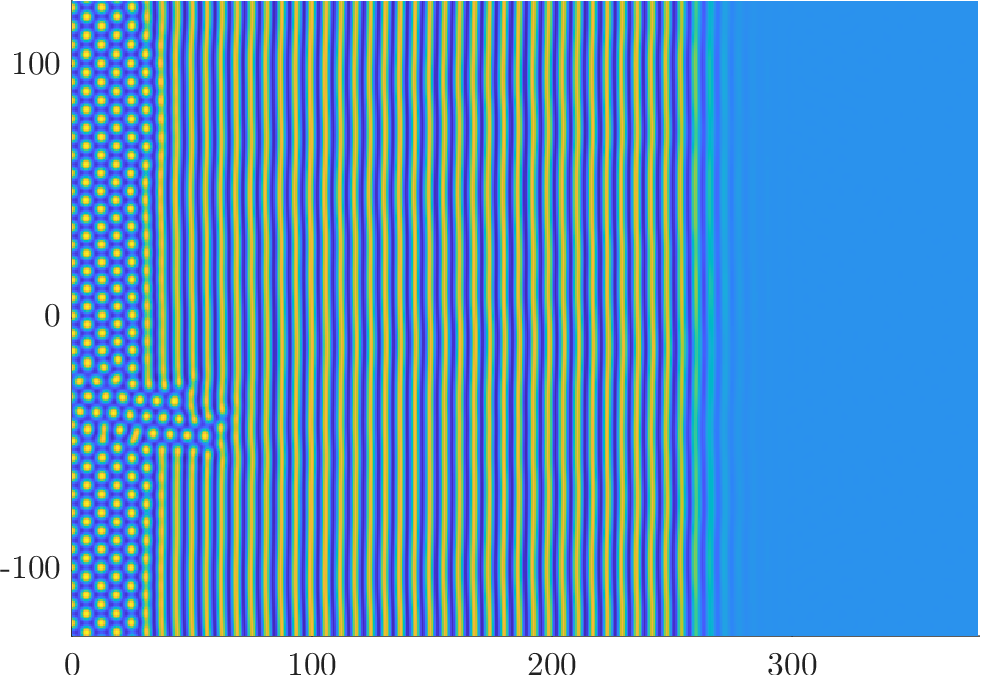}
    \caption{Invasion from a small vertical strip of random initial data in the planar Swift-Hohenberg equation~\eqref{e:sh} with quadratic coefficient $\gamma=2$ and $\mu=0.05$ (left) and $\mu=0.2$ (right). }\label{f:second}
\end{figure}

Somewhat ironically, the propagation via invasion front into an unstable state can  itself be unstable, despite the marginal stability in the leading edge. Such instabilities, usually observed in the wake, are in fact quite common and can lead to cascades of unstable states visited at any fixed location in space.  We pointed to some of those bifurcations in \cite{pnp} and  in \S\ref{s:3.1}.
The skew-product coupling presented there in \eqref{e:frontmorsechange} mimics the presence of invariant subspaces. A Hopf bifurcation in the wake is for instance simply not visible in the traveling-wave equation. In the modulated traveling-wave equation where the Hopf bifurcation is visible and gives rise to invasion fronts, the pure traveling-wave equation occurs as an equation on the invariant subspace of time-independent functions. More generally, instabilities affect the existence of a primary front only if they are resonant with the primary frequency. Instabilities of patterns in the wake such as in the Cahn-Hilliard to patterns with say double the wavelength are not resonant in this sense. They give rise to secondary fronts, so that the invasion process proceeds in stages. In fact, many of the examples discussed thus far involve such staged invasion processes, yet fairly little is understood about them. Most analytical results focus on situations with some explicit decoupling, or with a stable leading edge; see for instance
\cite{SS2001b,SS2005}, and  
\cite{GSU2004,GS2007}.  
The situation is complicated by the fact that instabilities are  at times generated at the front interface, yet they may impact the dynamics in the wake significantly~\cite{HSaccelerated}. Two somewhat eminent challenges in the context described here are:
\begin{enumerate}
    \item Explain the selection of harmonic wavenumbers through fronts; see ~\cite[Fig.~8--10]{S2017} and Fig.~\ref{f:eck}!
    \item Find selected fronts that leave  hexagonal patterns in the wake; see Fig.~\ref{f:second}!
\end{enumerate}
We conjecture that both cases can be understood through secondary fronts invading a primary pattern created in the leading edge, a short-wavelength pattern in the first case, and a stripe pattern parallel to the front interface in the second situation. Secondary invasion by hexagons in a Swift-Hohenberg equation is shown in Fig.~\ref{f:second}, both when hexagons form immediately after invasion, and when they form in a secondary, slower front. In both cases, the leading edge creates a striped pattern.
The predominant selection of stripes in the leading edge was first observed in~\cite{bose2013invasion} and generally proven for linearized problems in~\cite{HolzerScheelPointwiseGrowth}. Multi-dimensional invasion fronts were constructed in a Swift-Hohenberg context in~\cite{DSSS2003} using methods similar to~\cite{HSchn1999}, but not systematically studying selection mechanisms
 ~\cite[Fig.~9\&10]{S2017}. Staged invasion was studied in~\cite{pisnepo,csahok,PhysRevE.61.4835}.

Technically, the questions here bear some resemblance to questions of front invasion in shifting environments, when thinking of the invasion by a primary species as  preparing an environment for a second species to invade; see for instance~\cite{girardin2019invasion,WLFQ2022,LY2022}.


\subsection{Linear prediction through resonances}\label{s:res}
Our main results on selection rely on marginal stability in weighted function spaces, adapted to the pinched double roots. Such an approach can fail as seen in \eqref{e:cmc}, or, in the simplest situation, in
\begin{equation}\label{e:anom}
u_t=u_{xx}+u-u^2+v,\qquad v_t=d v_{xx}-\delta v
\end{equation}
with $d>2$, $0\leq \delta\ll 1$, in somewhat subtle ways.  The Gaussian $v$-profile effectively selects a front in the FKPP equation for $u$ with a speed $c>2$ that corresponds to a \emph{double pinched double root}, where spatial roots from the $u$- and the $v$-equation coalesce~\cite{HolzerAnomalous,HolzerAnomalous2}. Furthermore, replacing the $v$-source term in the $u$-equation by $v^2$ or $v^3$, acceleration still occurs for large enough $d$, although the linear equation does not predict this acceleration. This specific observation led to generally identifying resonances in the linear dispersion relation beyond the pinched double roots that have dominated this review. Criteria for speeds based on such generalized resonances were derived and corroborated in
\cite{FayeHolzerScheel}, although selection proofs beyond very simple skew-coupled systems with comparison principles do not exist. An example that illustrates the relevance in pattern-forming problems with competing modes was studied in
\cite{FayeHolzerScheelSiemer}, where nonlinear asymptotic stability (although not selection) of the apparently selected front was shown despite the fact that the front is spectrally unstable in any exponentially weighted space. It is remarkable that in this special situation \eqref{e:anom}, the selected front is not the steepest front, nor the slowest stable front. It is however marginally stable.

\subsection{More open problems and conclusions}
Somewhat tangential to the main focus of this review, there are many open questions.

\textbf{Invading periodic patterns.} We discussed instabilities of patterns created in the wake of invasion fronts, which lead to secondary fronts invading periodic, rather than constant, states. Investigating speed, frequencies, and selection of these fronts leads to somewhat intricate phenomena related to resonances, discussed before. Such fronts were studied in the phase separation context but most of the methods described in this review have not been adapted to this setting.

\textbf{Quasiperiodic and chaotic invasion.} When mentioning the relevance of resonances in the selection of states in the wake, one naturally asks if quasi-periodic invasion might be prevalent more generally in dynamical systems. We are not aware of systematic studies that would explore this question. Numerical studies in lattice-differential equations~\cite{Ben-Naim_2015} appear to suggest the absence of frequency locking. More dramatically, fronts in the complex coefficient Ginzburg-Landau equation may leave a chaotic state in their wake and there do not appear to be tools available that would characterize marginal stability and pushed versus pulled behavior in this situation.

\textbf{Inhomogeneous media.} In the presence of comparison principles, there is a large literature on front propagation in inhomogeneous media; see for example \cite{freidlingartner,nolen2012existence,shigesada86,xin2000front}. The dynamical systems perspective taken here seems well suited to address perturbative questions, such as slowly or rapidly varying media, from periodic to quasiperiodic, and to ergodic.

\textbf{Bounded domains.}  Instabilities in large bounded domains in the presence of transport are naturally mediated by invasion fronts and can lead to steep bifurcation diagrams; see~\cite{ssbasin,stegemerten} for  scalar case studies and~\cite{ADSS21} for conceptual and quantitative results more generally.

\textbf{Geometry and higher space dimensions.} Spreading of instabilities in higher-dimensional systems leads to novel curvature corrections to the propagation speed \cite{rrr} and, more dramatically, to intricate changes in selected patterns in the wake. It would clearly be interesting how manipulating the geometry of the growth process can be exploited to design patterns formed in the wake; see  for instance~\cite{growing}.


\textbf{Conclusion and more conjectures.}
We presented a dynamical systems perspective on the spreading  of instabilities in spatially extended systems that complements the robust techniques rooted in comparison principles. From this dynamics perspective, one of the key features of the invasion process is its ability to select states in its wake, even when the system intrinsically supports a multitude or even a continuum of stable states. We focused throughout on marginal stability as a selection criterion. When discussing generalizations of the pinched double root criterion to higher resonances in the dispersion relation and how they can lead to different speeds and selected states, we also realized that marginal stability is intrinsically quite subtle and far from completely understood. It is then interesting to return to other criteria, that happen to coincide with marginal stability in simple models: marginal positivity and steepest leading edge. In this light, we conclude this review with  a question that connects these concepts:
\begin{itemize}
    \item Are there examples of systems where marginal positivity, steepest decay in the leading edge, and marginal stability give different predictions (and which one is correct)?
\end{itemize}
Resonances, as discussed in \S\ref{s:res} provide intriguing case studies.
We realize that this review asks more questions than it answers, but we believe that the view point taken here can help organize our understanding of instabilities across many scientific areas, and at the same time stimulate research by pointing to natural dynamical systems inspired questions in this field. The flexibility and universality of the dynamical systems approach we have taken here is best summarized quoting van Saarloos~\cite{vanSaarloosReview}:
\begin{center}
\emph{\ldots the approach we have advanced here is ready to be put on rigorous footing. If this will be done, it will undoubtedly allow one to approach large classes of equations at one fell swoop.}
\end{center}

\fontsize{9}{9}\selectfont
\bibliographystyle{abbrv}

\bibliography{references}

\end{document}